\documentclass{sl}     
\usepackage{slsec}     
\usepackage{stud_log} 
\usepackage{slfoot_arx}
\usepackage{slthm}     
                       
\usepackage[latin9]{inputenc}
\usepackage{longtable}
\usepackage{float}
\usepackage{pmboxdraw}
\usepackage{mathrsfs}
\usepackage{amsmath}
\usepackage{amssymb}

\usepackage{enumitem}
\usepackage{graphicx}
\usepackage{wrapfig}
\usepackage{arydshln}
\usepackage{xcolor}

\makeatletter

\providecommand{\tabularnewline}{\\}

\usepackage{enumitem}		
\newlength{\lyxlabelwidth}      

\newenvironment{elabeling}[2][]%
{\settowidth{\lyxlabelwidth}{#2}
\begin{description}[font=\normalfont,style=sameline,leftmargin=\lyxlabelwidth,#1]}
{\end{description}}

\newcommand{\cml}{[\hspace{-1.5pt}[}
\newcommand{\cmr}{]\hspace{-1.5pt}]}
\usepackage{enumitem}

\AtBeginDocument{
  
}

\theoremstyle{definition}

\HeadingsInfo{Adrian Soncodi}{A Characterization of Non-Iterative Normal Modal Logics}

\begin{document}

\setcounter{page}{1}     

\raggedbottom{}

 

\AuthorTitle{Adrian Soncodi}{A Characterization of Non-Iterative Normal Modal Logics}

\begin{abstract}
Non-iterative normal modal logics are defined by axioms of modal degree 1.
In this paper we use calculations with normal forms
to determine the set of all possible non-iterative normal modal logics,
unimodal propositional extensions of $\mathbf{K}$.
To characterize them, we construct the complete set (lattice) of such logics
and we provide the generic specification of their axioms and their semantics.
\end{abstract}


\section{Introduction}

\paragraph{Scope.}\vspace{-2mm}
The normal modal logic $\mathbf{K}$ and its extensions are defined
in classical text books such as \cite{Blackburn} and \cite{Hughes}.
We use the terms \emph{modal logic}, \emph{modal logic system},
or simply \emph{logic} or \emph{system} interchangeably.
Note that, as a general rule, in this paper by ``system" we mean
a finitely-axiomatizable unimodal propositional extension of $\mathbf{K}$;
exceptions shall be pointed out explicitly.

A non-iterative system is one that can be axiomatized by formulas
with no iterated (i.e. no nested) modalities.
These systems have been studied before and several of their attributes
(e.g. having the finite model property or being canonical)
have been proven in papers like \cite{Lewis} and \cite{Surendonk}.
However, non-iterative systems are not determined solely from these properties,
because other systems may have them too.

So in this paper we attempt to provide a comprehensive characterization
of all non-iterative systems.
We basically ask the following question:

\vspace{-1mm}
\begin{itemize}[label=\textbf{(Q):}, leftmargin=*]
\item \textit{What are all the possible non-iterative systems,
i.e. how can we  define every single instance from this set explicitly?}
\end{itemize}
	
\vspace{-1mm}
To answer, we apply some of the methods described in \cite{Soncodi}.
We use normal forms,
specifically the \textit{disjunctive normal forms} (DNF) described in \cite{Fine},
which are a special case of the canonical forms defined in \cite{Soncodi}.
The results that we provide consist of:\vspace{-2mm}

\begin{itemize}[label=\textbullet, leftmargin=*]
\item A method to define/construct the lattice of non-iterative systems (partially ordered by inclusion),
showing the relationships between systems.
\item For each such system, the explicit formula for a defining axiom,
as well as the characterizing class of Kripke frames.
\end{itemize}

\pagebreak{}

\paragraph{Methodology.}
To begin with,
in Sections 2 and 3 we introduce some prerequisite concepts and basic properties.
Then the main results are in Sections 4-6.
The following is an outline of our strategy for answering \textbf{(Q)}:

\vspace{-2mm}
\begin{enumerate}[label=(\arabic*), ref=\arabic*]
\item First, in Section 4, using the Lindenbaum-Tarski algebra,
we partition all (classes of equiprovable) non-iterative modal formulas
into a countable set of \emph{modal contexts},
where each context includes only a finite number of formulas.
Any such formula can be used as an axiom for a non-iterative system.
But a number of different axioms may end up generating the same system.
So we prune the formulas to a minimal set of candidate axioms
that tentatively yield all the non-iterative systems.
We show that in each context
there is only a finite number of systems determined by our axioms,
and that they form a lattice.
This lattice is remarkably regular,
such that by using some ad-hoc coordinates
we can conveniently label each system by its unique position in the lattice.
\item Next, in Section 5,
we assemble the finite lattices from all contexts into a single candidate lattice
that should include all the non-iterative systems.
Here, each system is again determined by its axiom
and its unique position in the assembled, infinite lattice.
Still at this point we don't know yet
whether or not the systems in the lattice are all distinct.
\item Finally, in Section 6, for each system in the candidate lattice
we explicitly determine the corresponding Kripke semantics,
i.e. its characterizing class of Kripke frames.
And since we are able to prove
that each system is characterized by a different class of frames,
we conclude that
\textit{all the non-iterative systems in our candidate lattice are distinct}.
\end{enumerate}

\vspace{-2mm}

We achieve all this by
identifying the DNF \textit{patterns} that yield suitable axioms.
The relevant patterns turn out to be what we call
\textit{prime orbits of normal forms},
which are like building blocks for the systems axioms.

\section{Modal Contexts and their Systems}

\paragraph{Notations.}

\setlength{\jot}{3mm}
The notations to use in this paper
have been a significant challenge and are thus the result of an
extensive consideration.
The problem is that in our case the traditional notations
often result in very long formulas.
Since our paper consists mainly of algebraic calculations with normal forms,
which are essentially sums of products,
we opted for a shorter representation,
in the style of Boolean algebra.
For example we write $pq+!rs$
instead of $(p \land q) \lor (\sim r \land s)$.
So we kindly ask the reader to bear with us
in terms of the notations introduced further on,
because we believe that in the end it is worth it,
for the sake of compactness and readability of the proofs.

For formulas we use the Greek letters
$\varphi,\psi\ldots$ as well as specific notations
detailed further on. Propositional variables are $\left\{ p,q,r\ldots\right\} $
or $\{p_{i}\}$ and the logical constants are $\mathit{0}$ and $\mathit{1}$ (note the italics). The operators
are listed in the table below, in descending order of precedence:\vspace{-1mm}

{\renewcommand{\arraystretch}{1.2}
\begin{center}
\begin{tabular}{|c|l|}
\hline 
{\scriptsize{$\square$, $\lozenge$ and $!$}} & \emph{\scriptsize{necessity}}{\scriptsize{, }}\emph{\scriptsize{possibility}}{\scriptsize{
and }}\emph{\scriptsize{negation}}\tabularnewline[1pt]
\hline 
{\scriptsize{$\circ$}} & \emph{\scriptsize{uniform substitution (US) application }}{\scriptsize{(see
further on)}}\tabularnewline[1pt]
\hline 
{\scriptsize{$\cdot$}} & \emph{\scriptsize{conjunction}}{\scriptsize{, as well as }}\emph{\scriptsize{US}}{\scriptsize{
}}\emph{\scriptsize{composition}}{\scriptsize{ (but the symbol is typically omitted)}}\tabularnewline[1pt]
\hline 
{\scriptsize{$+$}} & \emph{\scriptsize{disjunction}}\tabularnewline[1pt]
\hline 
{\scriptsize{$\rightarrow$ and $\leftrightarrow$}} & \emph{\scriptsize{logical implication}}{\scriptsize{ and }}\emph{\scriptsize{logical
equivalence}}\tabularnewline[1pt]
\hline 
{\scriptsize{$=$ and $\approx_{\mathbf{S}}$}} & \emph{\scriptsize{identity}}{\scriptsize{ and}}\emph{\scriptsize{
equivalence (equiprovability)}}{\scriptsize{ }}\emph{\scriptsize{in}}{\scriptsize{
$\mathbf{S}$ of formulas; $\mathbf{S}$ omitted when $\mathbf{K}$ }}\tabularnewline[1pt]
\hline 
\end{tabular}\vspace{-0.5mm}
\par\end{center}}

Whenever necessary, parentheses are used for disambiguation.
Examples of well-formed formulas
(wff) are $p!q+qr\rightarrow\mbox{!}qr$ and $(\square p\leftrightarrow\square\lozenge p)\rightarrow\lozenge\square p$.

The \emph{modal degree} of a formula is the largest number
of nested modal operators found in any sub-formula of the given formula.
We also refer to it as the \emph{level} of the formula.
In the above examples, the modal degrees (levels) are 0 and 2 respectively.

\paragraph{Modal contexts.}
We denote systems by boldface and axioms by Roman letters.
\textbf{F} is reserved for the inconsistent system (consisting of all formulas).

Let $\mathcal{F}(v,d)$ be the set of all unimodal wff
in a number of variables not exceeding $v$ and of modal degree not exceeding $d$.
Consider a modal logic system $\mathbf{B}$ as a basis for the discussion.
In this paper $\mathbf{B}$ is typically $\mathbf{K}$, the weakest normal modal logic.
We define a \emph{modal context} $\mathbf{B}[v,d]$ as the quotient
$\mathcal{F}(v,d)/\approx_{\hspace{-1pt}_{\mathbf{B}}}$, the Lindenbaum-Tarski
algebra of classes of formulas equiprovable in $\mathbf{B}$.
But we tacitly equate a class with a representative from it,
so that we can still refer to these classes as formulas.
This is similar to writing
$\mbox{1+1=0}$ instead of $\hat{1}+\hat{1}=\hat{2}=\hat{0}$ in $\mathbb{Z}_{2}$,
which is often done and has the advantage of avoiding hats over large formulas.

The reason for this \emph{countable-contextualization} of formulas is as follows.
Let $\mathbf{S}$ be the system that extends $\mathbf{K}$ by a finite set of axioms $\{\alpha_{i}\}$.
Then there is an equivalent axiomatization for $\mathbf{S}$
using a single axiom $\alpha=\prod\alpha_{i}$, with $\alpha\in\mathbf{K}[v,d]$,
where $v$ and $d$ are the largest among all $\alpha_{i}$.
But as we shall see, $\mathbf{K}[v,d]$ is a finite set of (classes of equiprovable) formulas.
So for \textbf{(Q)} we ask what are all $\alpha\in\mathbf{K}[v,d]$ that yield distinct extensions of $\mathbf{K}$.
It turns out that for non-iterative contexts this question can be answered,
and we shall provide a constructive way to determine all the corresponding systems, axioms and semantics.

\paragraph*{Boolean contexts.}

A level 0 context $\mathbf{K}[v,0]$ consists of non-modal formulas with up to $v$ variables.
These can also be written in the well-known,
equivalent \emph{disjunctive normal form} (DNF).
For DNF, the \emph{factors} (conjuncts) in every \emph{term} (disjunct)
must consist of all $v$ variables,
where each variable is either complemented or not.
Such terms are usually referred to as \emph{full normal forms},
but we call them \emph{minterms},
as by ``normal form'' we typically mean the DNF representation of formulas in general.

Note that for $v<v'$, $\mathbf{K}[v,0]\subset\mathbf{K}[v',0]$,
since prior to (re)normalization in $\mathbf{K}[v',0]$
any missing variable $p_{i}$ can be reintroduced in the formula
by conjunction with $(p_{i}+\mbox{!}p_{i})$.

The following example shows a $\mathbf{K}[v,0]$ formula $\varphi$
and its DNF equivalent:\setlength{\jot}{1mm} \vspace{-12mm}

\begin{eqnarray*}
\varphi(p,q,r) & = & p+q+r\rightarrow(p\rightarrow q)r\\
 & \approx & pqr+\mbox{!}pqr+\mbox{!}p!qr+\mbox{!}p\mbox{!}q\mbox{!}r
\end{eqnarray*}

\vspace{-2mm}The above DNF formula can be displayed in matrix form as follows:\vspace{-2mm}

\[
\varphi(p,q,r)\approx\begin{array}{c|cccc|}
p & 1 & 0 & 0 & 0\\
q & 1 & 1 & 0 & 0\\
r & 1 & 1 & 1 & 0
\end{array}
\]
\vspace{-3mm}

We call this the minterms-matrix (short \emph{minmatrix}) representation of the formula.
The labels on the left denote the DNF factors, which in this case are all the propositional variables.
There is one column per minterm
and the minmatrix entries represent the \emph{state} of the row's factor in the corresponding minterm:
0 if the factor is complemented, otherwise 1.

Barring the reordering of rows and columns the minmatrix representation of a formula is unique.
However, we prefer to use a \emph{standard order} for rows and columns,
with the variables ordered alphabetically
or by some index and the minterms left-to-right in descending order of their state-tuples,
as shown in the example above.

For $v$ propositional variables $p_{0},\ldots,p_{v-1}$
there are $n=2^{v}$ Boolean minterms
denoted as $m_{0},\ldots m_{n-1}$.
Then $\mathbf{K}[v,0]$ is the set of $2^{n}$ distinct formulas
(including $\mathit{0}$ and $\mathit{1}$)
obtained by adding the minterms from any subset of $\{m_{i}\}$.
A formula is called \emph{positive} if its DNF includes $m_{n-1}$;
otherwise it is called a \emph{negative} formula.
(The reason is that any positive formula has a non-DNF equivalent formula
that does not use the negation operator or the constant $\mathit{0}$;
but we do not need to detail this property here.)

A notation that we shall use further on is $\mathrm{E}_{v}(k)$
for the set of $\mathbf{K}[v,0]$ formulas whose DNF contains precisely $k$ minterms.
Also, $\mathrm{E}_{v}^{+}(k)$ and $\mathrm{E}_{v}^{-}(k)$
denote the subsets of positive and negative formulas from $\mathrm{E}_{v}(k)$ respectively.
Note that if $e\in\mathrm{E}_{v}^{-}(k)$ then $\mbox{!}e\in\mathrm{E}_{v}^{+}(n-k)$, $0\leq k<n$.

\paragraph{Modal disjunctive normal form.}

In \cite{Fine}, the DNF for $\mathbf{K}[v,d]$ formulas with $d\geq1$ is defined recursively.
The corresponding minmatrices include:
\begin{itemize}
\item\vspace{-2mm}
\emph{Non-modal factors} consisting of all the $v$ propositional variables $p_{i}$
(the same as in non-modal minmatrices).
\item
\emph{Modal factors} of the form $\lozenge\mu_{i}$,
where $\mu_{i}$ is every DNF minterm from the previous context $\mathbf{K}[v,d-1]$.
\end{itemize}

\vspace{-1.5mm}For example, in $\mathbf{K}[1,1]$ the factors are $p$, $\lozenge p$
and $\lozenge\mbox{!}p$. In $\mathbf{K}[2,1]$ the factors are $p$,
$q$, $\lozenge(pq)$, $\lozenge(p\mbox{!}q)$, $\lozenge(\mbox{!}pq)$
and $\lozenge(\mbox{!}p\mbox{!}q)$, etc.
The following is an example minmatrix from $\mathbf{K}[2,1]$:\vspace{-2mm}

\begin{center}
$\varphi(p,q)=
\begin{array}{c|cccc:cc:ccc:cc|}
p & 1 & 1 & 1 & 1 & 1 & 1 & 0 & 0 & 0 & 0 & 0\\
q & 1 & 1 & 1 & 1 & 0 & 0 & 1 & 1 & 1 & 0 & 0\\
\hline \lozenge(pq) & 1 & 1 & 0 & 0 & 1 & 0 & 1 & 0 & 0 & 1 & 1\\
\lozenge(p!q) & 1 & 0 & 1 & 0 & 0 & 0 & 1 & 1 & 1 & 0 & 0\\
\lozenge(!pq) & 0 & 1 & 1 & 0 & 0 & 1 & 1 & 1 & 0 & 1 & 1\\
\lozenge(!p!q) & 1 & 1 & 0 & 1 & 1 & 1 & 1 & 1 & 1 & 1 & 0
\end{array}$\vspace{-1mm}

\par\end{center}

The horizontal line that separates the modal and non-modal factors is optional,
for convenience only.
The dashed vertical lines are also optional,
separating \emph{minmatrix sections},
which are subsets of minterms whose non-modal factors have the same states.
The product of non-modal factors in their respective states is called
the \emph{(non-modal) prefix} of the minterm.
The minmatrix section corresponding to minterms with prefix $m_{n-1}=\prod_{i=0}^{n-1}p_{i}$
is called its \emph{positive section}.
Observe that a standard order also applies to the modal factors,
which we denote by $\nu_{i}$, $0\leq i<n$; namely $\nu_{i}=\lozenge \mu_{i}$.

We wish to point out that the above ``matrix''
is only a schematic representation of a normalized formula.
This is useful to quickly show what is
relevant for its DNF, namely the states of the minterm factors.
But this representation can be equivalently replaced at any time
by the full algebraic formula.
Other than that, a minmatrix \emph{is}, in fact, a normalized formula.

Our notation for the minmatrix of a formula $\varphi$ is $[\varphi]$,
or $_{v}^{d}[\varphi]$ to emphasize the context.
This is still the same formula in DNF, therefore $\varphi\approx[\varphi]$.
The formula $\mathit{0}$ is represented by the empty minmatrix $[\:]=[\mathit{0}]$
and the formula $_v^d[\sum_i\mu_i]$ containing all the minterms from the context
by $[\mathit{1}]$
(notation not to be confused with the reference \cite{Blackburn}).

Per our convention to freely refer to $\mathbf{K}[v,d]$ classes as formulas,
$\mathbf{K}[v,d]$ is the (obviously finite) set of all minmatrices from the context.
Then the Boolean operations on formulas
can be performed as set operations on the corresponding minmatrix minterms:
union for disjunction, intersection for conjunction,
complementation (with respect to $[\mathit{1}]$) for negation, etc.

For this reason we may also interpret $[\varphi]$ as a set (rather than a sum) of minterms.
This notation overloading allows us to avoid
the constant use of back-and-forth conversion operators
between DCF formulas and their sets of minterms,
while all can still be disambiguated from the surrounding text
(as in, for example, $f:[\mathit{1}]\rightarrow[\mathit{1}])$.
And this also allows us to be brief by writing ``minmatrix intersection''
instead of ``the minmatrix that is the sum of the minterms
from the intersection of the sets of minterms from ...''

With this notation we can also write, for example, $[\varphi]\subset[\psi]$.
This partial order relationship then determines a lattice structure on $\mathbf{K}[v,d]$,
where each $\mathbf{K}[v,d]$ is a Boolean lattice
isomorphic to the power set of the set $\{\mu_{i}\}$ of all level $d$ minterms.

A context $\mathbf{K}[v',d']$ where $v'\geq v$ and $d'\geq d$
with at least one inequality being strict,
is called a \emph{successor} of $\mathbf{K}[v,d]$.
When $v=v'$ and $d'=d+1$, or when $d=d'$ and $v'=v+1$,
it is an \emph{immediate successor} context.
The reversed relationship says that $\mathbf{K}[v,d]$
is a \emph{predecessor} of $\mathbf{K}[v',d']$,
or an \emph{immediate predecessor} respectively.
Other contexts are incomparable.

Let $[\varphi] \in \mathbf{K}[v,d]$ have variables $p_{i}$, $0\leq i<v$.
Then $[\varphi] \approx [\varphi](p_v+!p_v)$,
and after re-normalization every $\mathbf{K}[v,d]$ minterm
yields several $\mathbf{K}[v+1,d]$ minterms
of an equiprovable formula in the immediate successor context.
We say that $_{v}^{d}[\varphi]$ was \emph{promoted} to $_{v+1}^{d}[\varphi]$.
Although we don't need to do it in this paper,
we can also promote $_{v}^{d}[\varphi]$ to an equiprovable $_{v}^{d+1}[\varphi]$
by artificially increasing the level,
e.g. by replacing the innermost level 0 formulas $\pi$ under the modalities,
as in $\lozenge\pi\approx\lozenge(_{v}^{1}[\mathit{1}]\,\pi)$,
and then re-normalizing the formula.
Thus, by induction we get $\mathbf{K}[v,d] \subset \mathbf{K}[v',d']$
for any successor context.

\paragraph*{Characteristic minmatrix and $\mathbf{K}[v,d]$-systems.}

From \cite{Fine} we have that a formula (or minmatrix) $\varphi$ is a theorem of $\mathbf{K}$
iff $_{v}^{d}[\varphi]=[\mathit{1}]$,
i.e.\ the DNF of $\varphi$ contains all the minterms from the corresponding context.
Obviously, the sum of all the minterms in a context is a Boolean tautology,
so if $\varphi,\psi$ are $\mathbf{K}[v,d]$ formulas then:\vspace{-5mm}

\begin{equation}
\vdash\mathbf{_{\hspace{-2pt}_{K}}\,}\varphi\rightarrow\psi\quad\text{iff}\quad[\varphi\rightarrow\psi]=[\mathit{1}]\quad\text{iff}\quad[\varphi]\subseteq[\psi]\label{eq:basic-1}
\end{equation}

\vspace{-5mm}
\begin{equation}
\vdash\mathbf{_{\hspace{-2pt}_{K}}}\,\varphi\leftrightarrow\psi\quad\text{iff}\quad[\varphi\leftrightarrow\psi]=[\mathit{1}]\quad\text{iff}\quad[\varphi]=[\psi]\label{eq:basic-2}
\end{equation}
\vspace{-5mm}

Let $\mathbf{S\in\mathrm{NExt}\mathbf{K}}$.
The \emph{characteristic minmatrix} (CMM) of $\mathbf{S}$ for a context $\mathbf{K}[v,d]$,
denoted as $\cml\mathbf{S}\cmr$ or $_{v}^{d}\cml\mathbf{S}\cmr$,
is defined as the conjunction (i.e. intersection of the minmatrices)
of all the $\mathbf{S}$-theorems from the context.
This definition is sound,
since under equivalence there is only a finite number
of $\mathbf{S}$-theorems per context.
Note that the minmatrices used in the intersection must belong to the context,
but the formal proof of the corresponding $\mathbf{S}$-theorems
may involve formulas from other contexts (and in fact this is often necessary).

The definition of a CMM combined with (\ref{eq:basic-1})
implies that a formula $_{v}^{d}\varphi$ is a theorem of $\mathbf{S}$
iff $\,{}_{v}^{d}\cml\mathbf{S}\cmr\subseteq\,_v^d[\varphi]$.
When this is the case,
one can construct a formal proof for $\varphi$ from $\cml\mathbf{S}\cmr$
using mainly the rule of \emph{substitution of equivalents} (EQ), defined
for example in \cite{Hughes}, as well as the propositional calculus
(PC) \emph{monotony} rule: from $\vdash p$ infer $\vdash p+q$. In
this sense the CMM is the ``strongest'' theorem that a system can
prove within a context.

For a formula $_{v}^{d}[\mathrm{\varphi}]$,
the system $\mathbf{S}$ that extends $\mathbf{K}$ with axiom $\mathrm{\varphi}$
is called \emph{the} $\mathbf{K}[v,d]$-\emph{system} corresponding to $\mathrm{\varphi}$,
and we denote it as $\mathbf{K}\mathrm{sys}(\mathrm{\varphi})$.
We also define $\cml\varphi\cmr\triangleq\cml\mathbf{S}\cmr$.
If $\varphi$ happens to be the strongest theorem that $\mathbf{S}$ can prove in the context,
then $\cml\varphi\cmr=[\varphi]$.
But it may be the case that $\cml\varphi\cmr\subsetneq[\varphi]$,
i.e. not every $[\varphi]\in\mathbf{K}[v,d]$ is necessarily the CMM of some system.
Indeed, if we can derive from $\mathrm{\varphi}$
another formula $\psi$ such that $[\varphi\psi]\subsetneq[\varphi]$,
then $[\varphi]$ cannot be $\cml\mathrm{\mathbf{S}}\cmr$ for any system $\mathbf{S}$,
since by definition $\cml\mathrm{\mathbf{S}}\cmr\subseteq[\varphi\psi]$.

As an example, below are the axiom of the well-known normal system $\mathbf{T}$ and its axiom T,
where in $\mathbf{K}[1,1]$ we have $\cml\mathbf{T}\cmr=\cml\mathrm{T}\cmr$ yet $\cml\mathbf{T}\cmr\subsetneq[\mathrm{T}]$:
\vspace{-2mm}

\begin{center}
$\mathrm{T}=\square p\rightarrow p=
\begin{array}{c|cccc:cc|}
p & 1 & 1 & 1 & 1 & 0 & 0\\
\hline
\lozenge p & 1 & 1 & 0 & 0 & 1 & 0 \\
\lozenge!p & 1 & 0 & 1 & 0 & 1 & 1
\end{array}$\,;
\quad
$\cml \mathbf{T}\cmr=
\begin{array}{c|cc:cc|}
p & 1 & 1 & 0 & 0\\
\hline
\lozenge p & 1 & 1 & 1 & 0 \\
\lozenge!p & 1 & 0 & 1 & 1
\end{array}\,;$\vspace{-1mm}

\par\end{center}

Let $\mathbf{K}\cml v,d\cmr$ be the set of CMMs from $\mathbf{K}[v,d]$.
Then $\mathbf{K}\cml v,d\cmr\subset\mathbf{K}[v,d]$ and $\mathbf{K}\cml v,d\cmr$
is typically not a sublattice of $\mathbf{K}[v,d]$.
But since it is finite and includes $\cml\mathbf{F}\cmr=[\mathit{0}]$ and $\cml\mathbf{K}\cmr=[\mathit{1}]$,
the set-based partial order between formulas can still be used to define
a lattice structure on $\mathbf{K}\cml v,d\cmr$.
However, the resulting $\mathbf{K}\cml v,d\cmr$ operations $\vee$ and $\wedge$
differ from the set operations inherited from $\mathbf{K}[v,d]$ as follows:
\bgroup \renewcommand*\theenumi{\alph{enumi}}\renewcommand*\labelenumi{\theenumi)}

\begin{elabeling}{00.0000.0000.000}
\item [{\textbf{\textsc{Theorem~1}}}] \noindent
\emph{Let} $\cml\mathbf{S'}\cmr$ \emph{and} $\cml\mathbf{S''}\cmr$
\emph{be CMMs from} $\mathbf{K}\cml v,d\cmr$\emph{. Then:}\vspace{-2mm}
\end{elabeling}

\noindent\vspace{0mm}a)
$\>\cml\mathbf{S'}\cmr\vee\cml\mathbf{S''}\cmr
=
\cml\mathbf{S'}\cmr\cup\cml\mathbf{S''}\cmr\>$,
\emph{(a CMM union is a CMM), and}

\noindent\vspace{0mm}b)
$\>\cml\mathbf{S'}\cmr\wedge\cml\mathbf{S''}\cmr
\subseteq
\cml\mathbf{S'}\cmr\cap\cml\mathbf{S''}\cmr\>$,
\emph{(a CMM intersection may not be a CMM)}.\vspace{2mm}

\noindent\textbf{\textsc{Proof.}}\qquad{}
Let $[\mathrm{S'}]\triangleq\cml\mathbf{S'}\cmr$
and $[\mathrm{S}'']\triangleq\cml\mathbf{S''}\cmr$
be the axioms of $\mathbf{S'}$ and $\mathbf{S''}$ from the context,
and $[\mathrm{S}]\triangleq\cml\mathbf{S'}\cmr\cup\cml\mathbf{S''}\cmr$.
For a), the inclusion $\supseteq$ is obvious,
but for $\subseteq$ we must show that $[\mathrm{S}]\triangleq\cml\mathbf{S'}\cmr\cup\cml\mathbf{S''}\cmr$ is indeed a CMM.
Let $\mathbf{S}=\mathbf{K}\mathrm{sys(S)}$.
If $[\mathrm{S}]$ is not a CMM,
then there is a theorem $\varphi$ of $\mathbf{S}$
such that $[\mathrm{S\,\varphi}]\subsetneq[\mathrm{S}]$,
so at least one of $[\mathrm{S'\,\varphi}]\subsetneq[\mathrm{S'}]$
or $[\mathrm{S''\,\varphi}]\subsetneq[\mathrm{S''}]$ holds.
But by (\ref{eq:basic-1}), both $\mathbf{S'}$ and $\mathbf{S''}$ already prove $\mathrm{S}$,
hence any theorem of $\mathbf{S}$, including $\varphi$,
so at least one of $\cml\mathbf{S'}\cmr$ and\emph{ }$\cml\mathbf{S''}\cmr$ is not a CMM,
contradicting our assumption.
Next, b) holds because any system that proves $[\mathrm{S'}]$ and $[\mathrm{S}'']$
proves at least $\cml\mathbf{S'}\cmr\cap\cml\mathbf{S''}\cmr$.
But here the inclusion may be strict,
since axioms $[\mathrm{S'}]$ and $[\mathrm{S}'']$ together
may prove a CMM that is stronger than this intersection.
\textsc{\scriptsize{\hfill{}$\blacksquare$}}\vspace{3mm}

It turns out, as we shall prove further on,
that for non-iterative contexts b) is also an equality.
But when $d>1$ the inclusion can actually be strict.
See an example in Section 3 of \cite{Soncodi},
where $\cml\mathbf{S5}\cmr\subsetneq\cml\mathbf{S4}\cmr\cap\cml\mathbf{B}\cmr$
in $\mathbf{K}[1,2]$.

As defined, a $\mathbf{K}[v,d]$ system is finitely-axiomatizable,
and it is also the weakest extension of $\mathbf{K}$ that has that CMM in the given context.
Obviously, every finitely-axiomatizable system is a $\mathbf{K}[v,d]$-system in some context,
as well as in all successor contexts thereafter.
A system that is not finitely-axiomatizable is not a $\mathbf{K}[v,d]$-system in any context,
but it still has a CMM in every context.

\vspace{-1mm}A caveat, however, is that $_{v}^{d}\mathbf{\cml}\mathbf{S}\cmr={}_{v}^{d}\mathbf{\cml}\mathbf{S'}\cmr$
does not imply $\mathbf{S}=\mathbf{S'}$,
unless of course they are both $\mathbf{K}[v,d]$-systems.
But systems that originate in different contexts may share a CMM within a context
when they prove the same theorems there.
Let $\mathbf{S'}\not=\mathbf{S}$, where $\mathbf{S}$ is a $\mathbf{K}[v,d]$-system and $\mathbf{S'}$ is a (stronger) system from a successor context,
yet still  $_{v}^{d}\mathbf{\cml}\mathbf{S}\cmr={}_{v}^{d}\mathbf{\cml}\mathbf{S'}\cmr$.
Then we say that in $\mathbf{K}[v,d]$ $\mathbf{S'}$ \emph{sinks into} $\mathbf{S}$.
For example, in $\mathbf{K}[1,1]$
the well-known level 2 systems $\mathbf{S4}$ and $\mathbf{S5}$ sink into $\mathbf{T}$.
But even though they share the level 1 CMM with $\mathbf{T}$,
$\mathbf{S4}$ and $\mathbf{S5}$ are not $\mathbf{K}[1,1]$-systems
(they are $\mathbf{K}[1,2]$-systems).

Let $\mathbf{K}\mathrm{sys}\cml v,d\cmr =
\{\mathbf{K}\mathrm{sys}(\varphi)\!:\!\varphi\in\mathbf{K}[v,d]\}$
be the set of $\mathbf{K}[v,d]$-systems.
The next statement says that,
when we limit our view to the finite number of systems from a context,
we can just map them to their CMMs.

\begin{elabeling}{00.0000.0000.00}
\item [{\textbf{\textsc{Claim~2}}}] \noindent
$\mathbf{K}\mathrm{sys}\cml v,d\cmr$
\emph{is a lattice isomorphic to} $\mathbf{K}\cml v,d\cmr$.
\end{elabeling}
\textbf{\textsc{Proof.}}\qquad{}
For any given context,
we have already shown that each system (including $\mathbf{K}$ and $\mathbf{F}$)
has a CMM that uniquely defines the system.
So we can induce a lattice structure on $\mathbf{K}\mathrm{sys}\cml v,d\cmr$
using a canonical isomorphism:
to each CMM $\cml\mathbf{S}\cmr\in\mathbf{K}\cml v,d\cmr$
we associate $\mathbf{K}\mathrm{sys}(\cml\mathbf{S}\cmr)$,
and:
\vspace{-0.5mm}\begin{eqnarray*}
\mathbf{K}\mathrm{sys(}\cml\mathbf{S'}\cmr)\vee\mathbf{K}\mathrm{sys}(\cml\mathbf{S''}\cmr)
& \triangleq &
\mathbf{K}\mathrm{sys}(\cml\mathbf{S'}\cmr\vee\cml\mathbf{S''}\cmr)
\\
\mathbf{K}\mathrm{sys(}\cml\mathbf{S'}\cmr)\wedge\mathbf{K}\mathrm{sys}(\cml\mathbf{S''}\cmr)
& \triangleq &
\mathbf{K}\mathrm{sys}(\cml\mathbf{S'}\cmr\wedge\cml\mathbf{S''}\cmr)
\end{eqnarray*}

\vspace{-1mm}\noindent can be used to define the operations in the new lattice.
\textsc{\scriptsize{\hfill{}$\blacksquare$}}

\vspace{3mm}
But we need to point out that $\mathbf{K}\mathrm{sys}\cml v,d\cmr$
is not a sublattice of $\mathrm{Next}\mathbf{K}$.
First, in $\mathbf{K}\mathrm{sys}\cml v,d\cmr$
the systems are "upside-down" with respect to
$\mathrm{Next}\mathbf{K}$,
since the partial order between CMMs is the reverse of the order
that results from defining the systems as sets of theorems.
Then the lattice operations from $\mathrm{Next}\mathbf{K}$,
denoted as $\oplus$ and $\odot$,
are not the same as $\vee$ and $\wedge$ from $\mathbf{K}\mathrm{sys}\cml v,d\cmr$.
In this paper we do not need to work with sublattices of $\mathrm{Next}\mathbf{K}$,
but the reader can verify that
for $\mathbf{S'},\mathbf{S''}\in\mathbf{K}\mathrm{sys}\cml v,d\cmr$:
\vspace{-1mm}\begin{eqnarray*}
\mathbf{\cml}\mathbf{S'}\odot\mathbf{S''}\cmr & = & \mathbf{\cml}\mathbf{S'}\cmr\vee\cml\mathbf{S''}\cmr\\
\mathbf{\cml}\mathbf{S'}\oplus\mathbf{S''}\cmr & = & \mathbf{\cml}\mathbf{S'}\cmr\wedge\cml\mathbf{S''}\cmr
\end{eqnarray*}

\pagebreak{}

\section{Prime Orbits of Minterms}

In this section we derive an initial \emph{necessary} condition
for a minmatrix to be the CMM of some system.

\paragraph*{Uniform substitutions.}

In the following we shall assume that the working context $\mathbf{K}[v,d]$
can accommodate all the formulas involved.
We write $\langle\alpha_{i}\rangle$
as a shorthand for $(\alpha_{0},\ldots,\alpha_{v-1})$, for example
$\varphi(p_{0},\ldots,p_{v-1})=\varphi\langle p_{i}\rangle$ and
$(\sigma_{0}\langle p_{j}\rangle,\ldots,\sigma_{v-1}\langle p_{j}\rangle)
=
\langle\sigma_{i}\langle p_{j}\rangle\rangle$,
even though not all these formulas
may effectively reference all the propositional variables $p_{i}$, $0\leq i<v$.
(But when they don't, if needed,
we can always promote them to equiprovables that do reference all the $v$ variables.)

Let $\sigma$ be the uniform substitution
$\langle p_{i}\rangle\mapsto\langle\sigma_{i}\langle p_{j}\rangle\rangle$,
where $\sigma_{i}$ are formulas, $0\leq i<v$.
The result, denoted as $\varphi\circ\sigma$,
of applying the substitution $\sigma$
to a formula $\varphi=\varphi\langle p_{i}\rangle$
is the formula $(\varphi\circ\sigma)\langle p_{i}\rangle$
obtained by consistently replacing every occurrence
of every propositional variable $p_{i}$ in $\varphi$
by the corresponding $\sigma_{i}\langle p_{j}\rangle$.
We write this as
$\varphi\circ\sigma
=
(\varphi\circ\sigma)\langle p_{i}\rangle
=
\varphi\langle\sigma_{i}\langle p_{j}\rangle\rangle$.
Formally, this operation is defined by the following rules,
applied recursively to the sub-formulas $\psi$, $\theta$, $\ldots$ that occur in $\varphi$:
\begin{elabeling}{00.00.0000}
\item [{(US-1)}] $0\circ\sigma\triangleq0$ and $1\circ\sigma\triangleq1$
\item [{(US-2)}] $p_{i}\circ\sigma\triangleq\sigma_{i}\langle p_{j}\rangle$
\item [{(US-3)}] $(\psi+\theta)\circ\sigma\triangleq\psi\circ\sigma+\theta\circ\sigma$
\item [{(US-4)}] $(!\psi)\circ\sigma\triangleq\mbox{!}(\psi\circ\sigma)$
\item [{(US-5)}] $(\lozenge\psi)\circ\sigma\triangleq\lozenge(\psi\circ\sigma)$
\end{elabeling}

Then the following are immediate consequences:
\begin{elabeling}{00.00.0000}
\item [{(US-6)}] $(\psi\,\theta)\circ\sigma\approx(\psi\circ\sigma)(\theta\circ\sigma)$
\item [{(US-7)}] $(\psi\rightarrow\theta)\circ\sigma\approx\psi\circ\sigma\rightarrow\theta\circ\sigma$
\item [{(US-8)}] $(\psi\leftrightarrow\theta)\circ\sigma\approx\psi\circ\sigma\leftrightarrow\theta\circ\sigma$
\item [{(US-9)}] $(\psi\not\leftrightarrow\theta)\circ\sigma\approx\psi\circ\sigma\not\leftrightarrow\theta\circ\sigma$
\item [{(US-10)}] $(\square\psi)\circ\sigma\approx\square(\psi\circ\sigma)$
\end{elabeling}

Also, using US-1 to US-10 above and the fact that $\varphi\approx\psi$ is defined as
$\vdash\mathbf{_{\hspace{-2pt}_{\mathbf{K}}}\,}\varphi\leftrightarrow\psi$,
one can easily prove the following additional properties:
\begin{elabeling}{00.00.0000}
\item [{(US-11)}] If $\varphi\approx\psi$ then $\varphi\circ\sigma\approx\psi\circ\sigma$
\item [{(US-12)}] If $\varphi\,\psi\approx0$ then $(\varphi\circ\sigma)(\psi\circ\sigma)\approx0$
\end{elabeling}

Our immediate interest is in level 0 substitutions,
i.e.~substitutions where all $\sigma_{i}$ are level 0 formulas.
We also call them \emph{context-preserving} substitutions,
because by applying them to any formula
neither the number of variables nor the modal degree increase.
If $v$ or $d$ actually decrease for a sub-formula,
it can always be promoted back to $\mathbf{K}[v,d]$.

There are $2^{2^{v}}$ formulas in $\mathbf{K}[v,0]$,
therefore $2^{v\cdot2^{v}}$ context-preserving substitutions $\sigma$
that can be applied to any $\varphi\in\mathbf{K}[v,d]$.
These substitutions are defined
independently of the formulas $\varphi$ and their modal degree.
Denote by $\mathcal{S}(v,0)$ the set of all level 0 substitutions $\sigma$ in $v$ variables.
The \emph{composition} $\sigma\sigma'$ of substitutions
$\sigma=\langle\sigma_{i}\langle p_{j}\rangle\rangle\in\mathcal{S}(v,0)$ and
$\sigma'=\langle\sigma_{i}'\langle p_{j}\rangle\rangle\in\mathcal{S}(v,0)$
is defined as follows:\vspace{-3mm}

\[
(\sigma\sigma')\langle p_{i}\rangle\triangleq\langle\sigma_{i}\langle\sigma_{j}'\langle p_{k}\rangle\rangle\rangle
\]

The composition of level 0 substitutions is obviously well-defined
(i.e. the result is context-preserving)
and its associativity is straightforward to verify.
With this operation $\mathcal{S}(v,0)$ is a monoid,
whose unit is the identical substitution
$\varsigma_{0}=\varsigma_{0}\langle p_{i}\rangle\triangleq\langle p_{i}\rangle$.
Then $\varphi\circ\sigma$ actually defines
a (right) monoid action of $\mathcal{S}(v,0)$ on $\mathbf{K}[v,d]$,
compatibility being ensured
since for any formula $\varphi\in\mathbf{K}[v,d]$ we have:\vspace{-3mm}

\[
\varphi\circ(\sigma\sigma')=\varphi\langle\sigma_{i}\langle\sigma_{j}'\langle p_{k}\rangle\rangle\rangle=\varphi\langle\sigma_{i}\langle p_{j}\rangle\rangle\circ\sigma'=(\varphi\circ\sigma)\circ\sigma'
\]
\vspace{-8mm}

\paragraph{Prime substitutions and prime orbits.}

Assume $\mathbf{S}$ extends $\mathbf{K}$ with axiom $\mathrm{S}$.
From $\vdash\mathbf{_{\hspace{-2pt}_{S}}\,}\mathrm{S}$
infer $\vdash\mathbf{_{\hspace{-2pt}_{S}}\,}\varphi$ for some $\varphi$,
then $\vdash\mathbf{_{\hspace{-2pt}_{S}}\,}\mathrm{S\,\varphi}$,
with $[\mathrm{S\,\varphi}]=[\mathrm{S}]\cap[\varphi]$.
If $[\mathrm{S}\,\varphi]\subsetneq[\mathrm{S}]$,
we say that the candidate CMM $[\mathrm{S}]$ \emph{collapses}
(by intersection with the minmatrix of some other theorem).
Let $\mathrm{\varphi=S\circ\sigma}$ for a context-preserving $\sigma$.
If $[\mathrm{S\,}\varphi]\subsetneq[\mathrm{S}]$
we say that $[\mathrm{S}]$ \emph{collapses under} $\sigma$;
otherwise if $[\mathrm{S}\,\varphi]=[\mathrm{S}]$,
i.e.~$[\mathrm{S}]\subseteq[\varphi]$,
we say that $[\mathrm{S}]$ is \emph{immune} to $\sigma$.

As we have seen, if $[\mathrm{S}]$ collapses under $\sigma$ then it can not be a CMM,
since this would imply $\mathrm{[S}]\subseteq[\mathrm{S}\,\varphi]$.
Thus, a candidate CMM must first of all be immune to all the context-preserving substitutions.
This is a necessary, although not sufficient condition for a minmatrix to be a CMM.

We therefore analyze which minmatrices are immune to level 0 substitutions.
We begin by considering a subset of level 0 substitutions
that we call \emph{prime substitutions}.
They are defined as the invertible elements of the monoid $\mathcal{S}(v,0)$,
hence they form a group that we denote by $\mathcal{S}_{p}(v,0)$.

It turns out that prime substitutions are precisely those level 0 substitutions $\varsigma$
that always transform a minterm into a single minterm.
As such, they generate automorphisms of the lattice $\mathbf{K}[v,d]$.
The following theorems establish the result.
\begin{elabeling}{00.0000.0000.000}
\item [{\textbf{\textsc{Theorem~3}}}]
\emph{For any given} $v$
\emph{, the group} $\mathcal{S}_{p}(v,0)$
\emph{is isomorphic to the symmetric group} $\mathrm{S}_{2^{v}}$\emph{.}
\end{elabeling}
\textbf{\textsc{Proof.}}\qquad{}For every $\varsigma\in\mathcal{S}_{p}(v,0)$,
let $f_{\varsigma}:\mathbf{K}[v,0]\rightarrow\mathbf{K}[v,0]$
be defined as $f_{\varsigma}(\varphi)\triangleq\varphi\circ\varsigma$.
We prove that $f_{\varsigma}$ is an automorphism of the lattice $\mathbf{K}[v,0]$.
First, we show that it is injective.
Assume $\varphi,\psi\in\mathbf{K}[v,0]$
and $f_{\varsigma}(\varphi)\approx f_{\varsigma}(\psi)$, i.e.~$\varphi\circ\varsigma\approx\psi\circ\varsigma$.
Since $\varsigma$ is invertible, we have
$\varphi\circ\varsigma\circ\varsigma^{-1}\approx\psi\circ\varsigma\circ\varsigma^{-1}$
and by compatibility $\varphi\circ(\varsigma\,\varsigma^{-1})\approx\psi\circ(\varsigma\,\varsigma^{-1})$,
i.e.~$\varphi\approx\psi$.
Next, being injective on the finite set $\mathbf{K}[v,0]$,
$f_{\varsigma}$ must be a bijection.
Properties US-1, US-3 and US-6 show that $f_{\varsigma}$
is also compatible with the lattice operations in $\mathbf{K}[v,0]$,
hence it is an automorphism.

We observe that if the prime substitution
$\varsigma$ is $\langle p_{i}\rangle\mapsto\langle\varsigma_{i}\langle p_{j}\rangle\rangle$,
then we have $f_{\varsigma}(p_{i})=p_{i}\circ\varsigma=\varsigma_{i}\langle p_{j}\rangle$,
therefore we can also write $\varsigma$ as
$\langle p_{i}\rangle\mapsto\langle f_{\varsigma}(p_{i})\rangle$.
Conversely, for every automorphism $f_{\varsigma}:\mathbf{K}[v,0]\rightarrow\mathbf{K}[v,0]$,
define $\varsigma$ to be the substitution
$\langle p_{i}\rangle\mapsto\langle f_{\varsigma}(p_{i})\rangle$.
Obviously $\varsigma\in\mathcal{S}(v,0)$ and it follows immediately that
$\langle p_{i}\rangle\mapsto\langle f_{\varsigma}^{-1}(p_{i})\rangle$ is its inverse,
also in $\mathcal{S}(v,0)$. \bgroup \renewcommand*\theenumi{\alph{enumi}}\renewcommand*\labelenumi{\theenumi)}

Thus, \emph{$\mathcal{S}_{p}(v,0)$} is isomorphic
to the group of automorphisms of $\mathbf{K}[v,0]$.
However, from lattice theory
any automorphism of a finite Boolean lattice 
is determined by its values on the atoms of the lattice,
and that these automorphisms correspond to the permutations of the atoms.
In our case the atoms are the $2^{v}$ minterms of $\mathbf{K}[v,0]$,
which proves our claim.
\textsc{\scriptsize{\hfill{}$\blacksquare$}}\vspace{3mm}

\begin{elabeling}{00.0000.0000.000}
\item [{\textbf{\textsc{Theorem~4}}}] \noindent \emph{Let $\varsigma\in\mathcal{S}_{p}(v,0)$.
Then for every $\mathbf{K}[v,d]$ context:}\end{elabeling}
\begin{enumerate}
\item \vspace{-2mm}\emph{The function}
$f_{\varsigma}:[\mathit{1}]\rightarrow[\mathit{1}]$\emph{ defined
as }$f_{\varsigma}(\mu)=\mu\circ\varsigma$\emph{ is a bijection on
the set }$_{v}^{d}[\mathit{1}]$\emph{ of minterms.}
\item \emph{The function }$f_{\varsigma}:\mathbf{K}[v,d]\rightarrow\mathbf{K}[v,d]$\emph{
defined as }$f_{\varsigma}(\varphi)=\varphi\circ\varsigma$\emph{
is a lattice automorphism.}
\end{enumerate}
\textbf{\textsc{Proof.}}\qquad{}Since $\mathbf{K}[v,d]=\wp([\mathit{1}])$
when minmatrices are viewed as sets, b) is a corollary of a).
We prove a) by induction on the modal level $d$.

For $d=0$ the result follows directly from Theorem 3.
Assume a) holds up to some level $d$.
Let $\mu=m\,\pi$ be a level $d+1$ minterm, where:
\begin{itemize}
\item \vspace{-2mm}$m$, the non-modal prefix of $\mu,$ is a level 0 minterm,
hence so is $m\circ\varsigma$.
\item $\pi$ is a product of all the level $d+1$ modal factors from the
set $\{\lozenge\mu_{k}\}$ in their respective states,
with $\mu_{k}$ being all the level $d$ minterms.
But by the induction hypotheses,
$\varsigma$ permutes the level $d$ minterms, therefore $\{\lozenge(\mu_{k}\circ\varsigma)\}=\{\lozenge\mu_{k}\}$.
Then $\pi\circ\varsigma$ is again
a product of all the level $d+1$ modal factors from $\{\lozenge\mu_{k}\}$
in correspondingly permuted states.
\end{itemize}

We apply property US-6 to conclude that $\mu\circ\varsigma\approx(m\circ\varsigma)(\pi\circ\varsigma)$
is a level $d+1$ minterm, so $f_{\varsigma}$ is well-defined. Also,
if $\mu_{1},\mu_{2}\in[\mathit{1}]$, then $\mu_{1}\circ\varsigma\approx\mu_{2}\circ\varsigma$
implies $\mu_{1}\circ\varsigma\circ\varsigma^{-1}\approx\mu_{2}\circ\varsigma\circ\varsigma^{-1}$
and $\mu_{1}\approx\mu_{2}$.
Thus, $f_{\varsigma}$ is injective on the finite set $_{v}^{d+1}[\mathit{1}]$,
hence it is a bijection.\textsc{\scriptsize{\hfill{}$\blacksquare$}}\vspace{3mm}

Theorem 4 implies that for every context $\mathbf{K}[v,d]$, $\mathcal{S}_{p}(v,0)$
determines a group action on the set $_{v}^{d}[\mathit{1}]$ of minterms. We
then define the (context-dependent) \emph{prime orbits} of minterms
as the orbits of this group action.\egroup

The \emph{order} (or \emph{count}) of a prime orbit is the number of its minterms.
Applying all the prime substitutions to a minterm $\mu$
does not necessarily yield only distinct minterms.
However, as we shall see shortly,
this yields an equal number of minterms per level 0 prefix.
Thus, in general, the prime orbits from a $\mathbf{K}[v,d]$ context
do not all have the same order.
But the order is always a multiple of $2^{v}$, its smallest value being $2^{v}$.
(The largest possible value would be $(2^{v})!$,
but this is typically not reached for $\mathbf{K}$.)

\vspace{2mm}

\begin{elabeling}{00.0000.0000.0000}
\item [{\textbf{\textsc{Theorem~5}}}] \emph{If a minmatrix} $\varphi$
\emph{includes some, but not all, the minterms of a prime orbit} $\omega$\emph{,}
\emph{then it collapses under some prime substitution} $\varsigma$\emph{.}
\end{elabeling}
\textbf{\textsc{Proof.}}\qquad{}This follows from Theorem 4 and the properties of orbits.
Let $\omega^{*}=\{\mu_{1},\ldots,\mu_{k}\}$
be a non-empty set of minterms $\omega^{*}\subsetneq\omega$ and $\mu_{x}\in\omega\setminus\omega^{*}$.
Since prime orbits are orbits
and the restriction of a group action to an orbit is transitive,
there exists a prime substitution $\varsigma$
such that $\mu_{1}\circ\varsigma\approx\mu_{x}$.
Then if $f_{\varsigma}:[\mathit{1}]\rightarrow[\mathit{1}]$
is $f_{\varsigma}(\mu)=\mu\circ\varsigma$,
we have $f_{\varsigma}(\mu_{1})\notin\omega^{*}$,
so $f_{\varsigma}(\omega^{*})\not\subset\omega^{*}$.
But since $\omega^{*}$ is finite,
there must be a $\mu_{i}\in\omega^{*}$, $1\leq i\leq k$
such that $\mu_{i}\notin f_{\varsigma}(\omega^{*})$,
hence for every $\mu_{y}\in\omega^{*}$, $\mu_{y}\circ\varsigma\not\approx\mu_{i}$.
Also, given any $\mu_{y}\in\varphi\setminus\omega^{*}$
with $\mu_{y}$ in some prime orbit $\omega'\not=\omega$
we have $\mu_{y}\circ\varsigma\in\omega'$,
and because prime orbits, as all group orbits, are disjoint, $\mu_{y}\circ\varsigma\notin\omega^{*}$
and once again $\mu_{y}\circ\varsigma\not\approx\mu_{i}$.
Thus, overall, there is no $\mu_{y}\in\varphi$
such that $\mu_{y}\circ\varsigma\approx\mu_{i}$,
therefore $\mu_{i}\notin\varphi\circ\varsigma$,
which means that $\varphi$ collapses under the substitution $\varsigma$.
\textsc{\scriptsize{\hfill{}$\blacksquare$}}\vspace{2mm}

\begin{elabeling}{00.0000.0000.0000}
\item [{\textbf{\textsc{Corollary~6}}}] \noindent
\emph{Any CMM must consist only of complete prime orbits}.\vspace{2mm}

\end{elabeling}

For a $\mathbf{K}[v,d]$-system $\mathbf{S}$, or for its CMM $\cml\mathrm{\mathbf{S}}\cmr$,
let $\Omega(\mathbf{S})=\Omega(_{v}^{d}\mathbf{S})$
be the set of prime orbits in $_{v}^{d}\cml\mathbf{S}\cmr$.
Obviously, $\Omega(\mathbf{K})=\Omega([\mathit{1}])$ in all contexts.

The reader may note that prime orbits from $\mathbf{K}[v,d]$ contexts
are a special case of the prime orbits from $\mathbf{E}[v,d]$ contexts
that we defined in \cite{Soncodi} using canonical forms,
and that results 3-6 are quite similar for both bases.

The results we have proved so far are general,
as they apply to all $\mathbf{K}[v,d]$ contexts.
In the next sections we shall focus on non-iterative contexts.

\pagebreak{}

\section{The $\mathbf{K}[\![v,1]\!]$ Lattices}

\vspace{-1mm}We now determine all the CMMs from a level 1 context.
To avoid collapses under prime substitutions,
we already know that they must be sums of complete prime orbits.
But we are able to impose stricter conditions
by requiring that they do not collapse under non-prime substitutions either.

Note that there are other methods to collapse a candidate CMM,
for example using theorems from other contexts.
But our findings imply that for level 1 systems
other methods yield nothing new,
i.e. if a level 1 minmatrix does not collapse under any level 0 substitution
then it is indeed a CMM.

\paragraph{Definitions}
Consider the formulas $\sum_{i=0}^{n-1}\varphi_{i}$
and $\prod_{i=0}^{n-1}\varphi_{i}$.
We will need to work with variations of them,
specifically some in which a fixed number of the $\varphi_i$ are complemented,
or are absent.

For any subformula $\varphi_{i}$ that may occur in a certain state
(complemented or not) in a sum or product, we define a \emph{pseudo-coefficient}
$\varepsilon_{i}$ to reflect its state:
$\varepsilon_{i}=1$ for $\varphi_{i}$ and
$\varepsilon_{i}=0$ for $\mbox{!}\varphi_{i}$.
Then we can write sums and products of mixed complementation uniformly as
$\sum_{i=0}^{n-1}\varepsilon_{i}\varphi_{i}$ or $\prod_{i=0}^{n-1}\varepsilon_{i}\varphi_{i}$.

When we need several $\varepsilon$-tuples in nested sums or products
we write $\varepsilon_{i}=\langle\varepsilon_{i,j}\rangle$ or $\varepsilon_{i,j}=\langle\varepsilon_{i,j,k}\rangle$,
with indices in the nesting order.

Let $\varepsilon=(\varepsilon_{0},\varepsilon_{1}
,\ldots,
\varepsilon_{n-1})=\langle\varepsilon_{i}\rangle$
be an $\varepsilon$-\emph{tuple} of pseudo-coefficients.
The \emph{count} function $\chi(\varepsilon)$
is defined as the number, range $0\ldots n$,
of $\varepsilon_{i}$ from $\varepsilon$ that are in state 1,
i.e. their $\varphi_i$ are not complemented.

\vspace{3mm}
Let $\mathbf{K}[v,1]$ be a non-iterative context and $n=2^{v}$.
We assume $v\geq1,$ as the degenerate context $\mathbf{K}[0,1]$
can be analyzed separately.
Recall that for the given context:
\begin{itemize}\vspace{-2mm}
\item $n$ is the number of level 0 (Boolean) minterms $m_{i}$, $0\leq i<n$,
\item $n$ is also the number of sections $s$ in a minmatrix, $0\leq s<n$,
\item $n$ is also the number of
level 1 modal DNF factors $\nu_{i}=\lozenge m_i$, $0\leq i<n$.
\end{itemize}

\vspace{-1mm}With the above notations,
an arbitrary minterm $\mu_j\in\mathbf{K}[v,1]$ is:\vspace{-2mm}

\[
\mu_j=m_{s}\prod_{i=0}^{n-1}\varepsilon_{j,i}\nu_{i}
\]

\noindent where $m_{s}$ is the Boolean prefix of $\mu$,
$s$ is the minmatrix section containing $\mu$, $\nu_{i}$ is the $i$-th modal
factor from $\mathbf{K}[v,1]$, namely $\nu_{i}=\lozenge m_{i}$,
and $\varepsilon_j=\langle\varepsilon_{j,i}\rangle$ is a suitable tuple
of pseudo-coefficients.

\vspace{1mm}For this minterm expression, we also define $\chi(\mu_j)\triangleq\chi(\varepsilon_j)$.

\pagebreak{}

\paragraph{Determining $\mathbf{K}[\![v,1]\!]$}
The next lemma shows that the minterms of a prime orbit follow a certain pattern:
in every minterm from a given section $s$:
\begin{itemize}\vspace{-2mm}
\item the modal factor $\nu_{s}$ has the same state (0 or 1), and
\item of the other modal factors, precisely $k$ ($0\leq k \leq n-1$) have state 1.
\end{itemize}

\vspace{-2mm}Note that this also implies that prime orbits
have an equal number ${n-1 \choose k}$ of minterms per section,
where $k$ depends on the prime orbit.

\begin{elabeling}{00.0000.0000.000}
\item [{\textbf{\textsc{Lemma~7}}}] \emph{Let $\omega\in\Omega([\mathit{1}])$
be a prime orbit and} $\mu=m_{s}\prod_{i=0}^{n-1}\varepsilon_{i}\nu_{i}$\emph{
a minterm from $\omega$. Then $\omega$ consists of all the minterms
}$\mu'=m_{s'}\prod_{i=0}^{n-1}\varepsilon'_{i}\nu_{i}$ \emph{for which
$\varepsilon'_{s'}=\varepsilon_{s}$ and
$\chi(\mu')=\chi(\mu)$.}\vspace{1mm}

\end{elabeling}
\textbf{\textsc{Proof.}}\qquad{}
Since prime orbits are group orbits,
by applying all \emph{$\varsigma\in\mathcal{S}_{p}(v,0)$}
to $\mu$ we must get all the minterms from $\omega$.
Fix \emph{$\varsigma\in\mathcal{S}_{p}(v,0)$}.
By Theorem 4, $\varsigma$ permutes the atoms of the lattice $\mathbf{K}[v,0]$,
which are all the $m_{i}$, therefore:\vspace{-5mm}

\begin{eqnarray*}
\mu\circ\varsigma & \approx & (m_{s}\prod_{i=0}^{n-1}\varepsilon_{i}\nu_{i})\circ\varsigma\;
\approx
\;(m_{s}\circ\varsigma)\prod_{i=0}^{n-1}\varepsilon_{i}(\nu_{i}\circ\varsigma)
\\
& \approx & m_{s'}
\prod_{i=0}^{n-1}\varepsilon_{i}\lozenge(m_{i}\circ\varsigma)\;
\approx
\;m_{s'}\prod_{i=0}^{n-1}\varepsilon'_{i}\nu_{i}
\end{eqnarray*}\vspace{-3mm}

\noindent where $s'$ is the new section of the transformed minterm
and $\varepsilon'_{i}$ are all the $\varepsilon{}_{i}$ permuted,
hence $\chi(\varepsilon')=\chi(\varepsilon)$.
Furthermore, since this permutation maps $m_{s}$ to $m_{s'}$ we have
$\nu_{s}\circ\varsigma
\approx
\lozenge(m_{s}\circ\varsigma)
\approx
\lozenge m_{s'}
\approx
\nu_{s'}$
and so $\varepsilon'_{s'}=\varepsilon_{s}$.

Next, by Theorem 3, prime substitutions correspond to \emph{all} permutations
of the atoms of the lattice $\mathbf{K}[v,0]$.
Then conversely,
if $\mu'$ is any minterm for which $\chi(\varepsilon')=\chi(\varepsilon)$
and $\varepsilon'_{s'}=\varepsilon_{s}$,
we can choose a $\varsigma'\in\mathcal{S}_{p}(v,0)$
that corresponds to a permutation mapping $m_{s'}$ to $m_{s}$ and
any $m_{j}\ne m_{s'}$ for which $\varepsilon'_{j}=0$ in $\mu'$
to an $m_{i}\ne m_{s}$ for which $\varepsilon_{i}=0$ in $\mu$.
This is possible since $\chi(\varepsilon')=\chi(\varepsilon)$,
and then similar calculations yield $\mu'\circ\varsigma'\approx\mu$,
hence all the minterms $\mu'$ satisfying the stated conditions belong to $\omega$.
\textsc{\scriptsize{\hfill{}$\blacksquare$}}\vspace{4mm}

With the above property
we can now determine all the prime orbits from $\mathbf{K}[v,1]$.
This will help us to find all the CMMs in the following way:
Suppose any $\varphi \in \mathbf{K}[v,1]$ could be a CMM.
Then, since there are $2^{v+n}$ minterms in the context,
we would need to check $2^{2^{v+n}}$ candidate CMMs to see whether or not they collapse.
And most of them will indeed collapse.
But from Corollary 6, any CMM is a sum of (complete) prime orbits.
As we now show there are only $2n$ prime orbits in $\mathbf{K}[v,1]$,
we will have only $2^{2n}$ candidate CMMs to check.

\pagebreak{}

\begin{elabeling}{00.0000.0000.00}
\item [{\textbf{\textsc{Theorem~8}}}]
\emph{The context }$\mathbf{K}[v,1]$
\emph{has $2n=2\cdot 2^v$ prime orbits}
\emph{denoted as}
$\mathrm{Vv}_{0}$, $\mathrm{Dd}_{0}$, $\mathrm{Dc}_{1}\ldots\mathrm{Dc}_{n-1}$
\emph{and} $\mathrm{Dw}_{1}\ldots\mathrm{Dw}_{n-1}$,
\emph{with formulas:}\vspace{-1mm}
\begin{eqnarray*}
\mathrm{Vv}_{0}\: & = & \sum_{s=0}^{n-1}m_{s}\mbox{!}\nu_{s}
\prod_{\substack{i=0\\i\ne s}}^{n-1}\mbox{!}\nu_{i} \quad = \quad 
\sum_{s=0}^{n-1}m_{s}
\prod_{i=0}^{n-1}\mbox{!}\nu_{i},
\\
\mathrm{Dd}_{0}\: & = & \sum_{s=0}^{n-1}m_{s}\phantom{\mbox{!}}\nu_{s}
\prod_{\substack{i=0\\i\ne s}}^{n-1}\mbox{!}\nu_{i},
\\
\mathrm{Dc}_{k}\: & = & \sum_{s=0}^{n-1}m_{s}\mbox{!}\nu_{s}\!\!\!\!
\sum_{\substack{j=1\\\chi(\varepsilon_{s,j})=k}}^{N(s)}
\prod_{\substack{i=0\\i\ne s}}^{n-1}\varepsilon_{s,j,i}\,\nu_{i},\qquad1\leq k<n,
\\
\mathrm{Dw}_{k} & = & \sum_{s=0}^{n-1}m_{s}\phantom{\mbox{!}}\nu_{s}\!\!\!\!
\sum_{\substack{j=1\\\chi(\varepsilon_{s,j})=k}}^{N(s)}
\prod_{\substack{i=0\\i\ne s}}^{n-1}\varepsilon_{s,j,i}\,\nu_{i},\qquad1\leq k<n,
\end{eqnarray*}
\end{elabeling}
\vspace{-2mm}
\emph{where} $N(s) = {n-1 \choose k}$ \emph{is the number of minterms in section} $s$
\emph{and} $\nu_{i}=\lozenge m_{i}$.
\setlength{\jot}{1mm}\vspace{3mm}

\noindent\textbf{\textsc{Proof.}}\qquad{}This is a consequence of the restrictions
imposed on the prime orbits by Lemma 7.
For simplicity, we can look only at the minterms from a single section $s$.
Minterms from the other sections then follow suit,
and for each prime orbit we add the minterms from all sections, $0$ to $n-1$.

Write a section $s$ minterm as $\mu_j=m_{s}\prod_{i=0}^{n-1}\varepsilon_{j,i}\nu_{i}$.
The choices for $\chi(\varepsilon_j)$ are $0,1,\ldots,n$.
But the state of modal factor $\nu_{s}$
must be the same across all section $s$ minterms.
When $\chi(\varepsilon_j)=0$ (all $\nu_{i}$ complemented) we obtain $\mathrm{Vv}_{0}$.
When $\chi(\varepsilon_j)=1$
we can have either $\varepsilon_{j,s}=1$, which yields $\mathrm{Dd}_{0}$,
or $\varepsilon_{j,s}=0$, in which case $\chi(\varepsilon_j)=1$
over the other $\varepsilon_{j,i}$
and the matching minterms yield $\mathrm{Dc}_{1}$.
Similar reasoning shows that when $\chi(\varepsilon_j)=k$,
with $2\leq k\leq n-1$,
the resulting prime orbits are either $\mathrm{Dw}_{k-1}$
(if $\varepsilon_{j,s}=1$) or $\mathrm{Dc}_{k}$ (if $\varepsilon_{j,s}=0$).
Finally, when $\chi(\varepsilon_j)=n$
we can only obtain $\mathrm{Dw}_{n-1}$.
Since these choices exhaust all possibilities for $\mu_j$,
there are no other $\mathbf{K}[v,1]$ prime orbits.
\textsc{\scriptsize{\hfill{}$\blacksquare$}}

\vspace{3mm}The reader may note that
$\mathrm{Vv}_{0}$ and $\mathrm{Dd}_{0}$
are quite similar to $\mathrm{Dc}_{k}$ and $\mathrm{Dw}_{k}$
and we could have labeled them $\mathrm{Dc}_{0}$ and $\mathrm{Dw}_{0}$,
for $k=\chi(\varepsilon_{s,j})=0$.
However, $\mathrm{Vv}_{0}$ and $\mathrm{Dd}_{0}$ are special,
since they are associated
with systems $\mathbf{Ver}$ and $\mathbf{Triv}$ respectively,
so we prefer to distinguish them in this way.

As per Corollary 6, a minmatrix that includes only complete prime
orbits is immune to prime substitutions.
Yet it may still collapse under non-prime substitutions.
In the worst case it may collapse to $\cml\mathbf{F}\cmr=[\mathit{0}]$,
which is the bottom element in $\mathbf{K}\cml v,1\cmr$.

\pagebreak{}

We now show that, starting from minmatrix $[\mathit{0}]$,
in order to construct (non-collapsing) level 1 CMMs,
prime orbits must be added in a particular order,
determined by certain dependencies between prime orbits.

To state these dependencies we introduce the concept of
\emph{prime orbit coverage} under a (typically non-prime) substitution $\sigma$.
To begin with, where $\Omega(_v^1[\mathit{1}])=\{\omega_{i}:0\leq i<2n\}$
we observe that:\vspace{-2mm}
\[
\sum_{i=0}^{2n-1}(\omega_{i}\circ\sigma)
\approx
(\sum_{i=0}^{2n-1}\omega_{i})\circ\sigma
\approx
[\mathit{1}]\circ\sigma
\approx
[\mathit{1}]\vspace{-3mm}
\]
and since all $\omega_{i}$ are pairwise disjoint, by US-6 and US-12
all $\omega_{i}\circ\sigma$ are also pairwise disjoint.
Comparing each $\omega_{i}\circ\sigma$ with the corresponding $\omega_{i}$
we observe that under $\sigma$ some of these minmatrices lose (potentially
all) their minterms, while others gain minterms, yet overall we still
find all the minterms from the context redistributed among several $\omega_{i}\circ\sigma$.

Any minterm $\mu\in\omega_{i}\circ\sigma$ is said to be \emph{covered}
by $\omega_{i}$ under $\sigma$.
Whenever $(\omega_{i}\circ\sigma)\cap\omega_{j}\ne\varnothing$
and $(\omega_{i}\circ\sigma)\cap\omega_{j}\subsetneq\omega_{j}$
we say that $\omega_{i}$ \emph{partially covers} $\omega_{j}$ under
$\sigma$. And when $\omega_{j}\subseteq(\omega_{i}\circ\sigma)$
we say that $\omega_{i}$ \emph{fully covers} $\omega_{j}$ under
$\sigma$.
Note that self-coverage (when $i=j$) is included here.

\begin{elabeling}{00.0000.0000.00}
\item [{\textbf{\textsc{Claim~9}}}] \emph{For level 1 prime orbits, if}
$\omega_{i}$ \emph{covers} $\omega_{j}$
\emph{(partially or fully), then a CMM can include} $\omega_{j}$
\emph{only if it also includes} $\omega_{i}$.
\end{elabeling}

\noindent\textbf{\textsc{Proof.}}\qquad{}Since all $\omega_{j}\circ\sigma$ are pairwise disjoint,
for any minterm $\mu\in\omega_{j}$ not covered by $\omega_{j}$ under $\sigma$
there must be precisely one other $\omega_{i}$ that covers $\mu$ under that $\sigma$.
So if a candidate CMM $\xi$ includes $\omega_{j}$ but not $\omega_{i}$,
then $\xi\circ\sigma$ includes at best the incomplete prime orbit $\omega_{j}$,
i.e. $\xi$ collapses under $\sigma$ to a minmatrix that does not include $\omega_{j}$,
therefore it cannot be a CMM.
\textsc{\scriptsize{\hfill{}$\blacksquare$}}\vspace{3mm}

Using this argument,
we can infer some \emph{necessary} dependencies between prime orbits in our CMMs.
The calculations below may seem a little complex,
but in fact they are just basic Boolean algebra with DNF formulas,
where we need to keep track of the prime orbits minterms
\emph{based on their definition from Theorem 8.}

\begin{elabeling}{00.0000.0000.00}
\item [{\textbf{\textsc{Theorem~10}}}] \emph{Let $\xi\in\mathbf{K}\cml v,1\cmr$ be a $\mathrm{CMM}$
with $\xi\ne\cml\mathbf{F}\cmr$. Then:}\vspace{-2mm}

\begin{elabeling}{0000000}
\item [{(DR1):}] \emph{$\xi$ includes at least one of $\mathrm{Vv}_{0}$
or $\mathrm{Dd}_{0}$.}
\item [{(DR2):}] \emph{$\xi$ includes $\mathrm{Dw}_{k}$, $1\leq k<n$,
only if it includes $\mathrm{Dd}_{0}$ and all $\mathrm{Dw}_{l}$, $1\leq l<k$.}
\item [{(DR3):}] \emph{$\xi$ includes $\mathrm{Dc}_{k}$, $1\leq k<n$,
only if it includes $\mathrm{Dd}_{0}$
as well as all $\mathrm{Dw}_{l}$,
$1\leq l<k$, and all $\mathrm{Dc}_{l}$, $1\leq l<k$.}
\end{elabeling}
\end{elabeling}

\noindent\textbf{\textsc{Proof.}}\qquad{}Let $\sigma$ be the non-prime substitution where
$p_{0}\mapsto p_{0}\mbox{!}m_{n-1}\mbox{!}m_{n-2}$
and $p_{1},p_{2},\ldots,p_{v-1}$ are unchanged.
This is the same as $p_{0}\mapsto p_{0}\sum_{i=0}^{n-3}m_{i}$,
or $p_{0}\mapsto\mathit{0}$ if $v=1$ and the context has only $n=2^{v}=2$ minterms.
Then:
\vspace{-1mm}\[
m_{n-1}\circ\sigma \,\approx\,
(\prod_{i=0}^{n-1}p_{i})\circ\sigma \approx
p_{0}\mbox{!}m_{n-1}\mbox{!}m_{n-2}\prod_{i=1}^{n-1}p_{i} \approx
\mbox{!}m_{n-1}\mbox{!}m_{n-2}m_{n-1} \approx \mathit{0}
\]
\vspace{-4mm}\begin{eqnarray*}
m_{n-2}\circ\sigma
&\!\!\!\approx\!\!\!&
(\mbox{!}p_{0}\prod_{i=1}^{n-1}p_{i})\circ\sigma \approx
(!p_{0}+m_{n-1}+m_{n-2})(m_{n-1}+m_{n-2})\qquad\,\,
\\
&\!\!\!\approx\!\!\!&
m_{n-1}+m_{n-2}
\end{eqnarray*}

\noindent and for the other $m_{j}$ (if any, i.e if $v>1$)
we have $m_{j}\circ\sigma \approx m_{j}$.
Indeed, when $j$ is even, $j+1<n-2$ and the state of $p_{0}$ in $m_{j}$ is 0, therefore:
\vspace{-2mm}\[
m_{j}\circ\sigma \approx
(\mbox{!}p_{0}\prod_{i=1}^{n-1}\varepsilon_{i}p_{i})\circ\sigma \approx
(!p_{0}+m_{n-1}+m_{n-2})(m_{j}+m_{j+1}) \approx
m_{j}
\]
\noindent\vspace{-2mm}and when $j$ is odd, $j<n-2$ and the state of $p_{0}$ in $m_{j}$ is 1, therefore:
\[
m_{j}\circ\sigma \approx
(p_{0}\prod_{i=1}^{n-1}\varepsilon_{i}p_{i})\circ\sigma \approx
p_{0}(\sum_{i=0}^{n-3}m_{i})(m_{j}+m_{j-1}) \approx m_{j}
\]

On level 1, for all $\mu$ in section $s=n-1$, $\mu\circ\sigma\approx\mathit{0}$;
and for $s<n-1$:
\begin{eqnarray*}
\mu\circ\sigma & \approx &
(m_{s}\circ\sigma) \prod_{i=0}^{n-1}\varepsilon_{i}\lozenge(m_{i}\circ\sigma)
\\
& \approx &
(m_{s}\circ\sigma)\,
(\varepsilon_{n-1}\lozenge\mathit{0})\,
(\varepsilon_{n-2}\lozenge(m_{n-1}+m_{n-2}))\!
\prod_{i=0}^{n-3}\varepsilon_{i}\lozenge m_{i}
\end{eqnarray*}

(Note that when $n=2^{v}=2$, the last product is $\mathit{1}$; otherwise $n\geq 4$.)
\vspace{1mm}

\textbf{For} $\mathbf{Vv}_{0}$, there is one minterm per section, and if
$\mu$ is in section $n-2$:\vspace{-6mm}

\begin{eqnarray*}
\mu\circ\sigma
& \approx &
(m_{n-1}+m_{n-2})\,\mbox{!}\lozenge\mathit{0}\,\mbox{!}\lozenge(m_{n-1}+m_{n-2})\!
\prod_{i=0}^{n-3}\mbox{!}\lozenge m_{i}
\\
& \approx & (m_{n-1}+m_{n-2})\,
\mbox{!}\lozenge m_{n-1}\,
\mbox{!}\lozenge m_{n-2}\,
\prod_{i=0}^{n-3}\mbox{!}\lozenge m_{i} \approx \mu_{v}+\mu
\end{eqnarray*}
which covers $\mu$ and the vanished section $n-1$ minterm $\mu_{v}$.
Minterms from remaining sections, if any, are unchanged by $\sigma$,
as for $s < n-2$ we have:
\[
\mu\circ\sigma \approx
m_{s}\,\mbox{!}\lozenge\mathit{0}\,\mbox{!}\lozenge(m_{n-1}+m_{n-2})
\prod_{i=0}^{n-3}\mbox{!}\lozenge m_{i} \approx \mu
\]
so overall $\mathrm{Vv}_{0}$ covers only itself, fully.
Thus $\mathrm{Vv}_{0}$ does not depend on other prime orbits
and it may be added to any CMM without causing a collapse.

\pagebreak{}

\textbf{For} $\mathbf{Dd}_{0}$, also with one minterm per section, and
for $\mu$ in section $n-2$:\vspace{-6mm}

\begin{eqnarray*}
\mu\circ\sigma & \approx &
(m_{n-1}+m_{n-2})\,\mbox{!}\lozenge \mathit{0}\,\lozenge(m_{n-1}+m_{n-2})\!
\prod_{i=0}^{n-3}\mbox{!}\lozenge(m_{i}\circ\sigma)
\\
& \approx & (m_{n-1}+m_{n-2})(\lozenge m_{n-1}+\lozenge m_{n-2})
\prod_{i=0}^{n-3}\mbox{!}\lozenge m_{i}
\end{eqnarray*}
\vspace{-3mm}

Here, prefix $m_{n-1}$ covers the vanished $\mathrm{Dd}_{0}$ minterm from section $n-1$
and also a $\mathrm{Dw}_{1}$ and a $\mathrm{Dc}_{1}$ minterm,
while prefix $m_{n-2}$ covers $\mu$ and more.
$\mathrm{Dd}_{0}$ minterms from remaining sections, if any,
are unchanged by $\sigma$,
so overall $\mathrm{Dd}_{0}$ covers itself fully
and covers $\mathrm{Dw}_{1}$ and $\mathrm{Dc}_{1}$ partially.
Thus, $\mathrm{Dd}_{0}$ does not depend on any other prime orbit,
but both $\mathrm{Dw}_{1}$ and $\mathrm{Dc}_{1}$ depend on $\mathrm{Dd}_{0}$.

\textbf{For} $\mathbf{Dw}_{k}$, we just saw that $\mathrm{Dw}_{1}$ depends on $\mathrm{Dd}_{0}$.
If $n=2$ there is no other $\mathrm{Dw}_{k}$.
Otherwise we pick a $\mathrm{Dw}_{k}$ minterm $\mu$ from section $n-2$
with $\varepsilon_{n-1}=0$ and
$\chi(\varepsilon)=k$ (excluding $\varepsilon_{n-2}=1$), $1\leq k \leq n-2$.
Then:\vspace{-6mm}

\begin{eqnarray*}
\mu\circ\sigma & \approx & (m_{n-1}+m_{n-2})\,
\mbox{!}\lozenge \mathit{0}\,\lozenge(m_{n-1}+m_{n-2})\!
\prod_{i=0}^{n-3}\varepsilon_{i}\lozenge m_{i}
\\
& \approx & m_{n-2}
(\lozenge m_{n-1}\lozenge m_{n-2}+\lozenge m_{n-1}\mbox{!}\lozenge m_{n-2})
\prod_{i=0}^{n-3}\varepsilon_{i}\lozenge m_{i}+\ldots
\end{eqnarray*}
\vspace{-3mm}

Here, because of $\lozenge m_{n-1}$ instead of $!\lozenge m_{n-1}$,
prefix $m_{n-2}$ adds 1 to $\chi(\varepsilon)$
and it includes both states for $\lozenge m_{n-2}$,
so it covers $\mathrm{Dw}_{k+1}$ as well as $\mathrm{Dc}_{k+1}$ minterms.
Therefore $\mathrm{Dw}_{k+1}$ and $\mathrm{Dc}_{k+1}$
depend on $\mathrm{Dw}_{k}$,
thus transitively on all previous $\mathrm{Dw}_{i}$, $1\leq i \leq k$,
as well as on $\mathrm{Dd}_{0}$.

\textbf{For} $\mathbf{Dc}_{k}$ we saw that $\mathrm{Dc}_{1}$ depends on both $\mathrm{Dd}_{0}$ and $\mathrm{Dw}_{1}$.
If $n=2$ there is no other $\mathrm{Dc}_{k}$.
Otherwise we pick a $\mathrm{Dc}_{k}$ minterm $\mu$ from section $n-3$
with $\varepsilon_{n-1}=0$, $\varepsilon_{n-2}=1$, $\varepsilon_{n-3}=0$
and $\chi(\varepsilon)=k$ (overall), $1\leq k \leq n-2$. 
Then:\vspace{-6mm}

\begin{eqnarray*}
\mu\circ\sigma & \approx & m_{n-3}\,
\mbox{!}\lozenge\mathit{0}\,\lozenge(m_{n-1}+m_{n-2})\,\mbox{!}\lozenge m_{n-3}\!
\prod_{i=0}^{n-4}\varepsilon_{i}\lozenge m_{i}
\\
& \approx & m_{n-3}(\lozenge m_{n-1}\,\lozenge m_{n-2}\,\mbox{!}\lozenge m_{n-3}
\prod_{i=0}^{n-4}\varepsilon_{i}\lozenge m_{i}+\ldots)
\end{eqnarray*}
\vspace{-3mm}

Here, because of $\lozenge m_{n-1}\lozenge m_{n-2}$
instead of $\mbox{!}\lozenge m_{n-1}\lozenge m_{n-2}$,
prefix $m_{n-3}$ adds 1 to $\chi(\varepsilon)$,
so it covers a $\mathrm{Dc}_{k+1}$ minterm plus more.
Therefore $\mathrm{Dc}_{k+1}$,
which already depends on all $\mathrm{Dw}_{i}$, $1\leq i \leq k$
and on $\mathrm{Dd}_{0}$, also depends on $\mathrm{Dc}_{k}$,
thus transitively on all previous $\mathrm{Dc}_{i}$, $1\leq i \leq k$.

Finally, by aggregating all the dependencies that we found above,
we get precisely the rules stated by this theorem.
\textsc{\scriptsize{\hfill{}$\blacksquare$}}

\pagebreak{}

The choice for $\sigma$ in the above proof
is due to a result that we don't actually need to prove here.
In fact,
every non-prime substitution reveals a set of dependencies between prime orbits.
But many substitutions yield the same set of dependencies,
so we can partition all of them into equivalence classes.
One particular class of \emph{critical substitutions} happens to determine
the strongest set of dependencies between prime orbits.
The $\sigma$ from Theorem 10
is just a simple and convenient substitution from this class.

From rules DR1-DR3 we build a candidate lattice $\mathbf{K}\cml v,1\cmr$,
which is as shown in Figure \ref{fig:Kv1-CMMs}.
In the diagram the CMMs are represented by dots.
Prime orbits are marked along the edges
such that one can determine what sets are included in any particular CMM.
A simple calculation shows that
the total number of CMMs from $\mathbf{K}\cml v,1\cmr$ is $n(n+3)$.

\begin{figure}[H]
\begin{centering}
\includegraphics[scale=0.64]{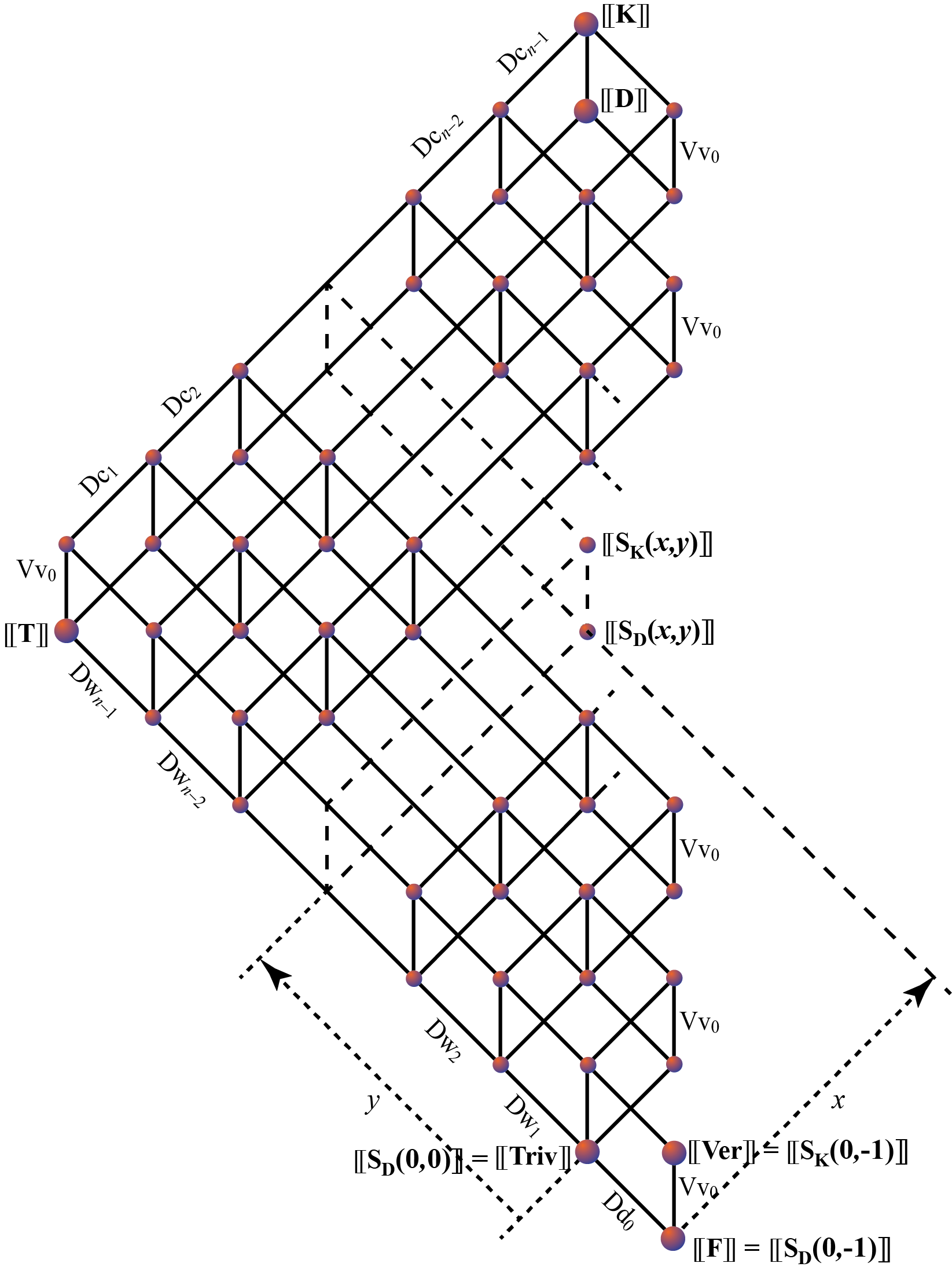}
\par\end{centering}

\caption{\label{fig:Kv1-CMMs}The candidate lattice $\mathbf{K}\cml v,1\cmr$}
\end{figure}

\pagebreak{}

By Theorem 10, there are no other CMMs in $\mathbf{K}\cml v,1\cmr$,
since any other minmatrix from the context
collapses under a critical substitution.
Yet in principle, any of these candidate CMMs may still collapse by other methods.
But note that adjacent $\mathbf{K}\cml v,1\cmr$ CMMs are a single prime orbit apart.
So if we can show that they all correspond to distinct systems,
then they are indeed all CMMs,
as there are no other suitable minmatrices to collapse to.
This would also imply that level 0 prime substitutions
plus a single critical substitution
suffice to determine the CMM of a level 1 system.

We observe that a number of $\mathbf{K}\cml v,1\cmr$ CMMs
include prime orbit $\mathrm{Vv}_{0}$ and an equal number do not.
These are referred to as \emph{K-plane} and \emph{D-plane}
CMMs respectively and the corresponding $\mathbf{K}\mathrm{sys}\cml v,1\cmr$
systems as \emph{K-plane} and \emph{D-plane} systems respectively.
The diagram does not provide names for them, except in a few cases
where we anticipate the position of some well-known systems. But due
to the regularity of this lattice we can denote the CMMs by the $x$
and $y$ coordinates defined as illustrated.

Figure \ref{fig:K01-and-K11} shows the lattices for the first two contexts.
Note that we chose generic names for the $\mathbf{K}[v,1]$ prime orbits,
but in fact they are context-dependent;
so in principle they must be renamed in each particular context. 

\begin{figure}[h]
\begin{centering}
\includegraphics[scale=0.66]{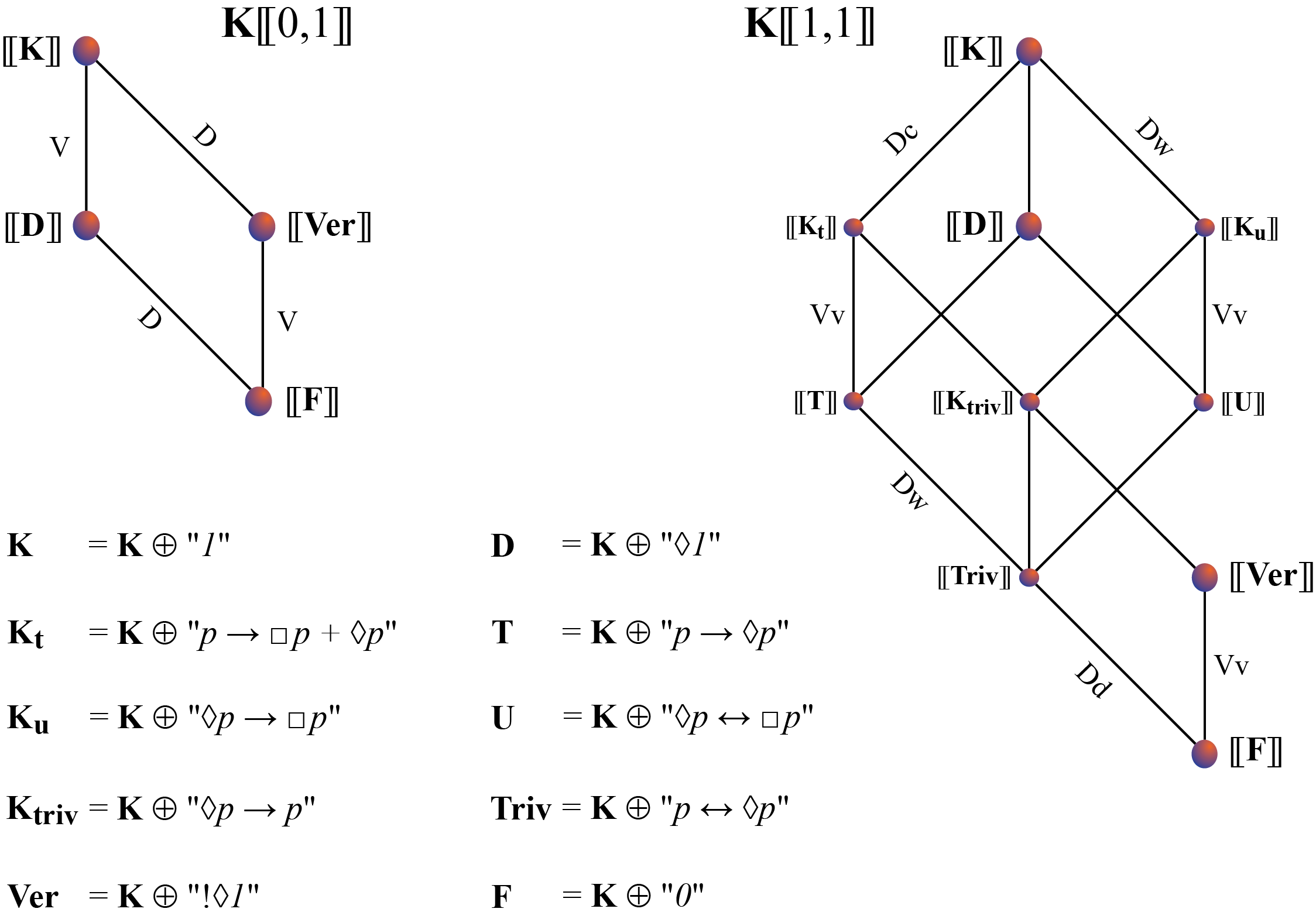}
\par\end{centering}

\caption{\label{fig:K01-and-K11}The lattices $\mathbf{K}\cml0,1\cmr$ and
$\mathbf{K}\cml1,1\cmr$ }
\end{figure}

\pagebreak{}

\section{The Candidate Lattice $\mathbf{Ksys}\mathbf{[\![*,1]\!]}$}

Having all the $\mathbf{K}\cml v,1\cmr$ lattices,
we note that it doesn't make sense to assemble them into a single lattice,
since a system has different CMMs in different contexts.
However, it does make sense to assemble all the corresponding
(as per Claim 2) lattices $\mathbf{\mathbf{K}\mathrm{sys}}\cml v,1\cmr$
into a single one.

We denote the assembled lattice by $\mathbf{\mathbf{K}\mathrm{sys}}\cml*,1\cmr$.
For now, it is a candidate lattice,
as we still need to prove that all the systems in it are distinct.
And we recall from Section 2 that this is not meant to be a sublattice of $\mathrm{Next}\mathbf{K}$.
It is just our way to present all the finitely-axiomatizable, non-iterative systems.

We construct $\mathbf{\mathbf{K}\mathrm{sys}}\cml*,1\cmr$
according to the pattern observed in Section 4,
as the fractal-like limit of all $\mathbf{\mathbf{K}\mathrm{sys}}\cml v,1\cmr$.
By definition the generic $\mathbf{\mathbf{K}\mathrm{sys}}\cml v,1\cmr$
includes all the systems from context $\mathbf{K}[v,1]$.
But then it must obvioulsly include all the systems from all predecessor contexts too,
only we do not know how exactly the latter
are positioned in $\mathbf{\mathbf{K}\mathrm{sys}}\cml v,1\cmr$.
This will become clear when we determine the semantics for these systems,
i.e. their characterizing classes of Kripke frames.

\paragraph{Construction of $\mathbf{\mathbf{K}\mathrm{sys}}\cml*,1\cmr$.}\vspace{-1mm}
For both the K-plane and the D-plane we define an $x$-axis and a
$y$-axis, so that we can refer to their $\mathbf{K}\mathrm{sys}\cml v,1\cmr$
systems as $\mathbf{S_{K}}(x,y)$ and $\mathbf{S_{D}}(x,y)$ respectively,
where $0\leq x\leq n-1$, $-1\leq y\leq n-1$ and $x\leq y+1$. For
$\mathbf{K}\mathrm{sys}\cml*,1\cmr$ we add an extra coordinate $*$
defined such that $*>n$ for all integers $n$.
Then we map $\mathbf{\mathbf{K}\mathrm{sys}}\cml v,1\cmr$
to $\mathbf{K}\mathrm{sys}\cml*,1\cmr$ as follows:\setlength{\jot}{1mm}
\begin{enumerate}
\item For $0\leq x\leq n-1$ and $-1\leq y\leq n-2$, $\mathbf{S_{K}}(x,y)$
and $\mathbf{S_{D}}(x,y)$ from $\mathbf{K}\mathrm{sys}\cml v,1\cmr$
map to position $(x,y)$ in the K-plane and D-plane of $\mathbf{\mathbf{K}\mathrm{sys}}\cml*,1\cmr$
respectively.
\item For $0\leq x\leq n-2$, $\mathbf{S_{K}}(x,n-1)$ and $\mathbf{S_{D}}(x,n-1)$
from $\mathbf{K}\mathrm{sys}\cml v,1\cmr$ map to position $(x,*)$
in the K-plane and D-plane of $\mathbf{K}\mathrm{sys}\cml*,1\cmr$
respectively.
\item Finally, $\mathbf{S_{K}}(n-1,n-1)$ and $\mathbf{S_{D}}(n-1,n-1)$
from $\mathbf{K}\mathrm{sys}\cml v,1\cmr$ map to position $(*,*)$
in the K-plane and D-plane of $\mathbf{K}\mathrm{sys}\cml*,1\cmr$
respectively.
\end{enumerate}
The resulting lattice $\mathbf{\mathbf{K}\mathrm{sys}}\cml*,1\cmr$
is illustrated in Figure \ref{fig:K[vd]-systems-all}.
The solid lines show 
how the generic $\mathbf{K}\mathrm{sys}\cml v,1\cmr$ lattice is embedded in it,
therefore we can add up the prime orbits to get an axiom for each system.
Also note how for every successor context, new systems will appear in the shaded area only.

From Section 4, we know that all non-iterative system must be in this lattice,
as the minmatrix of any other axiom would collapse in some context.
Next, by showing that every $(x,y)$ position from $\mathbf{K}\mathrm{sys}\cml *,1\cmr$
is associated with a distinct semantics
will prove that the corresponding systems are also distinct,
i.e. that $\mathbf{K}\mathrm{sys}\cml *,1\cmr$
contains all the systems we are looking for.

\pagebreak{}

\begin{figure}[H]
\begin{centering}
\includegraphics[scale=0.56]{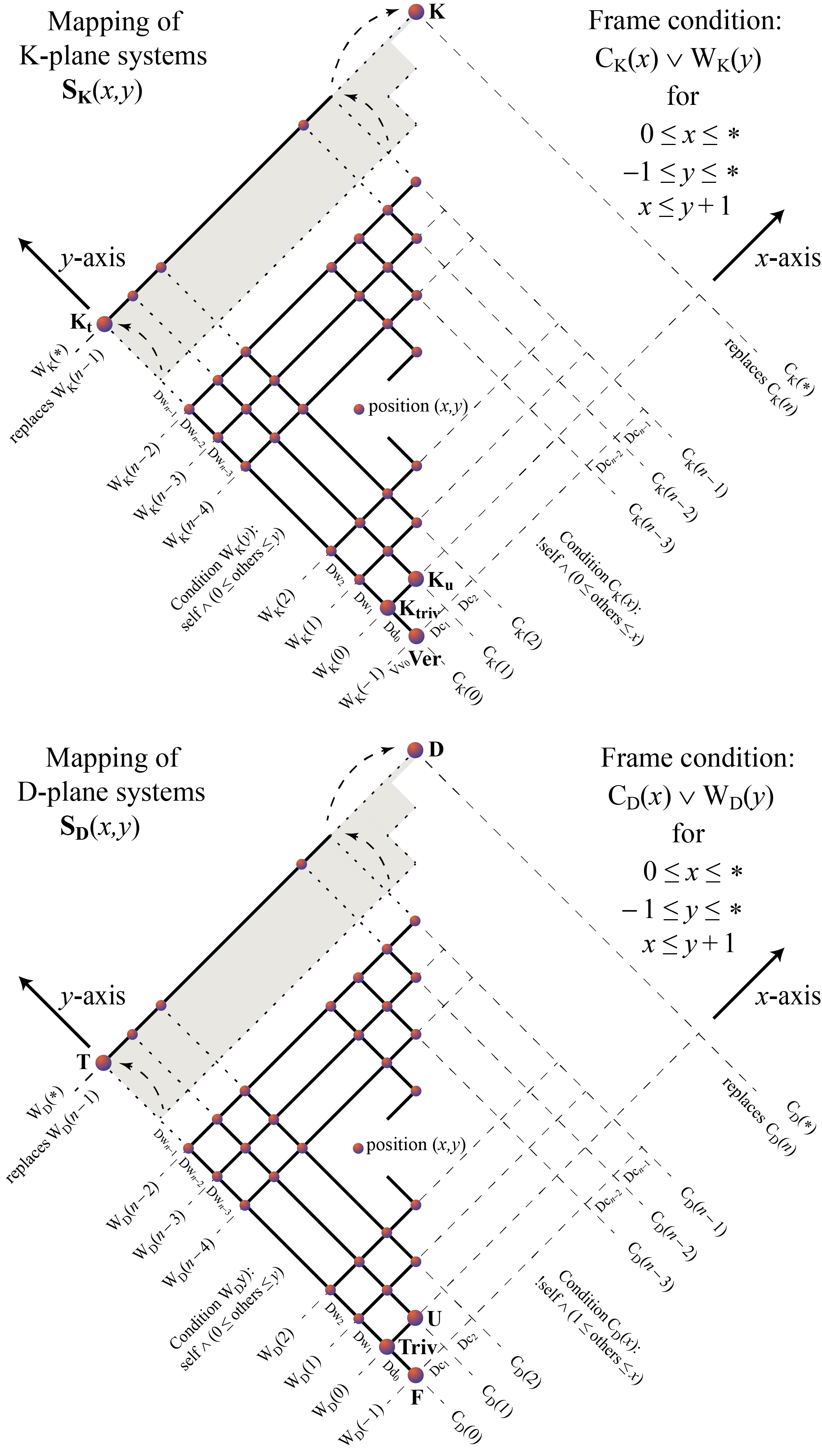}
\par\end{centering}
	
\caption{\label{fig:K[vd]-systems-all}The lattice $\mathbf{K}\mathrm{sys[*,1]}$ $-$
Mapping of K-plane and D-plane systems}
\end{figure}

\pagebreak{}

\paragraph{Axioms.}In the next section we will need the defining axioms
$\alpha_{\mathbf{S_{K}}(x,y)}$ and $\alpha_{\mathbf{S_{D}}(x,y)}$
of the $\mathbf{K}\mathrm{sys}\cml *,1\cmr$ systems.
Since all $\mathbf{K}\mathrm{sys}\cml v,1\cmr$
are embedded in $\mathbf{K}\mathrm{sys}\cml *,1\cmr$,
we could take as axioms the CMMs from the generic $\mathbf{K}\cml v,1\cmr$.
These can be built from Figure 1,
starting with $\cml\mathbf{F}\cmr$, for example:\vspace{2mm}

\qquad$\alpha_{\mathbf{S_{K}}(1,2)}
=
\mathrm{Vv_{0}+Dd_{0}+Dc_{1}+Dw_{1}+Dw_{2}}$\vspace{2mm} 

Instead, however, we shall build some equivalent, more suitable axioms.
These are constructed using the following trick:
Consider any system $\mathbf{S}$
with axiom $\cml \mathbf{S}\cmr \approx \sum_{\omega_j\in \Omega(\mathbf{S})}\omega_{j}$.
If we add any minterms to this minmatrix,
then, as long as we do not end up with new, complete prime orbits,
these extra minterms will be ``trimmed" by prime substitutions.
Hence, if we define:\vspace{-2mm}
\[
\alpha_{\mathbf{S}}
\,\,\triangleq\,\,
(\prod_{i=0}^{v-1}p_{i})\!\!\!\sum_{\omega_j\in \Omega(\mathbf{S})}\!\!\!\omega_{j}+
\mbox{!}\prod_{i=0}^{v-1}p_{i}
\,\,\approx\,\,
m_{n-1}\rightarrow m_{n-1}\!\!\!\sum_{\omega_j\in \Omega(\mathbf{S})}\!\!\!\omega_{j}
\]\vspace{-4mm}

\noindent then section $n-1$ (the positive section) of this minmatrix
includes all the minterms from the selected prime orbits, and only those minterms,
while the remaining sections include all the minterms.
Consequently, $\alpha_{\mathbf{S}}$ will collapse to $\cml \mathbf{S}\cmr$
under prime substitutions, so it is an equivalent axiom for $\mathbf{S}$.

Using this trick and Theorem 10,
we construct a similar minmatrix $\alpha_\omega$ for each prime orbit $\omega$,
where the positive section of $\alpha_\omega$ includes only the minterms from $\omega$,
while the other sections include all their minterms.
Then with our notation for a generic minterm $\mu_j=m_{s}\prod_{i=0}^{n-1}\varepsilon_{j,i}\nu_{i}$,
we get:\vspace{-5mm}
\setlength{\jot}{0.1mm}

\begin{eqnarray*}
\alpha_{\mathrm{Vv}_{0}} & \triangleq &
m_{n-1} \rightarrow
m_{n-1}\!\!\sum_{\substack{\chi(\varepsilon_j)=0 \\ \,\varepsilon_{j,n-1}=0}}
\prod_{i=0}^{n-1}\varepsilon_{j,i}\lozenge m_i
\,\,\approx\,\,
m_{n-1} \rightarrow
\!\!\sum_{\substack{\chi(\mu_j)=0 \\ \,\varepsilon_{j,n-1}=0}}
\mu_j
\\
\alpha_{\mathrm{Dd}_{0}} & \triangleq &
m_{n-1}\rightarrow
m_{n-1}\!\!\sum_{\substack{\chi(\varepsilon_j)=1 \\ \,\varepsilon_{j,n-1}=1}}
\prod_{i=0}^{n-1}\varepsilon_{j,i}\lozenge m_i
\,\,\approx\,\,
m_{n-1} \rightarrow
\!\!\sum_{\substack{\chi(\mu_j)=1 \\ \,\varepsilon_{j,n-1}=1}}
\mu_j
\end{eqnarray*}

\noindent\vspace{-5mm}and for $1\leq k<n$:

\begin{eqnarray*}
\alpha_{\mathrm{Dc}_{k}} &\triangleq &
m_{n-1}\rightarrow
m_{n-1}\sum_{\substack{\chi(\varepsilon_j)=k \\ \,\varepsilon_{j,n-1}=0}}
\,\,\prod_{i=0}^{n-1}\varepsilon_{j,i}\lozenge m_i
\,\,\approx\,\,
m_{n-1} \rightarrow
\!\!\sum_{\substack{\chi(\mu_j)=k \\ \,\varepsilon_{j,n-1}=0}}
\mu_j
\\
\alpha_{\mathrm{Dw}_{k}} & \triangleq &
m_{n-1} \rightarrow m_{n-1}
\!\!\sum_{\substack{\chi(\varepsilon_j)=k+1 \\ \,\varepsilon_{j,n-1}=1}}
\prod_{i=0}^{n-1}\varepsilon_{j,i}\lozenge m_i
\,\,\approx\,\,
m_{n-1} \rightarrow
\!\!\!\!\!\!\sum_{\substack{\chi(\mu_j)=k+1 \\ \,\varepsilon_{j,n-1}=1}}
\!\!\!\!\mu_j
\end{eqnarray*}\vspace{-3mm}

\noindent\vspace{0mm}where all the $\mu_j$ above belong to section $n-1$.

\pagebreak{}

(Note that in $\alpha_{\mathrm{Dw}_{k}}$,
$\chi(\mu_j)=k+1$ because it now includes $\varepsilon_{j,n-1}=1$.)

\vspace{2mm}Proceeding to $\mathbf{S_{K}}(x,y)$,
we add according to its CMM from Figure 1,
to obtain the following axiom:\vspace{-3mm}
\[
\alpha_{\mathbf{S_{K}}(x,y)} \,\triangleq\,
\alpha_{\mathrm{Vv_{0}}}+\sum_{k=1}^{x}\alpha_{\mathrm{Dc_{k}}}+
\alpha_{\mathrm{Dd_{0}}}+\sum_{k=1}^{y}\alpha_{\mathrm{Dw_{k}}}
\]\vspace{-3mm}

\noindent\vspace{0mm}and this minmatrix, as per the trick above,
collapses to $\cml \mathbf{S_{K}}(x,y)\cmr$.
Then, after replacing the $\alpha_\omega$ of the level 1 prime orbits $\omega$,
we combine the resulting formulas using the simple Boolean equivalence:
\vspace{-1mm}\[
(p\rightarrow q) + (p\rightarrow r) \,\approx\, p\rightarrow q+r.
\]

\noindent\vspace{-2mm}and we immediately obtain:

\[
\alpha_{\mathbf{S_{K}}(x,y)}
\,\,\approx\,\,
m_{n-1}\,\,\rightarrow
\!\!\sum_{\substack{\chi(\mu_j) \leq x \\ \,\varepsilon_{j,n-1}=0}}
\!\!\!\!\mu_j \>\>+ 
\!\!\sum_{\substack{\chi(\mu_j) \leq y+1 \\ \,\varepsilon_{j,n-1}=1}}
\!\!\!\!\mu_j
\]

\noindent\vspace{0mm}where again, all the $\mu_j$ are from section $n-1$.
This is the formula that we shall use in the next section
to determine the semantics of $\mathbf{S_{K}}(x,y)$.

As a caveat, when $y=n-1$ we call this $\alpha_{\mathbf{S_{K}}(x,*)}$,
unless $x=y=n-1$ when it is called $\alpha_{\mathbf{S_{K}}(*,*)}$.
\vspace{1mm}

Finally, for $\mathbf{S_{D}}(x,y)$ we can immediately write the axiom:\vspace{-3mm}
\[
\alpha_{\mathbf{S_{D}}(x,y)}\,\,\approx\,\,
\lozenge\mathit{1}\,\alpha_{\mathbf{S_{K}}(x,y)}
\qquad\qquad\qquad\qquad\quad
\]

\vspace{-2mm}\noindent since $\lozenge\mathit{1}\approx\mbox{!}\mathrm{Vv}_{0}$
and conjunction with it eliminates precisely the prime orbit $\mathrm{Vv}_{0}$.
Then this conjunction,
together with trimming the minmatrix to complete prime orbits,
turns $\cml \mathbf{S_{K}}(x,y)\cmr$ into $\cml \mathbf{S_{D}}(x,y)\cmr$.

\paragraph{Additional axioms}
After further calculations, it is possible to show that:

\setlength{\jot}{3mm}\vspace{-5mm}
\begin{eqnarray*}
\alpha'_{\mathbf{S_{K}}(x,y)}
& \!\!\approx & \!\!
m_{n-1}
\rightarrow
\mbox{!}\lozenge m_{n-1}\!\!\!\!
\sum_{e_{i}\in\mathrm{E}_{v}^{+}(x+1)}\!\!\!\!\square\,e_{i}
\,\,+\,\,
\lozenge m_{n-1}\!\!\!\!
\sum_{e_{i}\in\mathrm{E}_{v}^{+}(y+1)}\!\!\!\!\square\,e_{i}
\\
& \!\!\approx & \!\!
m_{n-1}
(\sum_{e_{i}\in\mathrm{E}_{v}^{+}(x+1)}\!\!\!\!\square\,e_{i}
\,\rightarrow\,\lozenge m_{n-1}\,\,)
\,\,\,\rightarrow\,\,\,
\lozenge m_{n-1}\!\!\!\!
\sum_{e_{i}\in\mathrm{E}_{v}^{+}(y+1)}\!\!\!\!\square\,e_{i}
\end{eqnarray*}

\noindent\setlength{\jot}{2mm}also works as an axiom for ${\mathbf{S_{K}}(x,y)}$.
Although we do not need to use this last formula in our paper,
we mentioned it here because it is a nice axiom with only positive subformulas.
The sums $\sum\square\,e_i$ can often be optimized,
for example they may reduce to sums that include only a subset of the $e_i$
generated by cyclic permutations of their $\langle p_{j}\rangle$.

\pagebreak{}

\section{Semantics for $\mathbf{Ksys}[\![*,1]\!]$ systems}
In the following we assume that the reader is familiar with Kripke frames
and the semantics of normal modal logics.
Let $\mathscr{F}=(W,R)$ be a Kripke frame and $w,\,w'\in W$.
If $wR\,w'$ we say that $w$ \emph{sees} $w'$.
Let ``self'' be the condition that a world $w$ is reflexive, i.e.~it sees itself.
Let ``others'' denote the number of $w'\neq w$ seen by $w$,
with ``$\mathrm{others}\leq *$" meaning $w$ sees any number of $w'\neq w$.
We define the following frame conditions on any $w$
(and since these are not modal logic formulas,
we use the classical notations):\vspace{-5mm}

\begin{eqnarray*}
\,\,\mathrm{C_{K}}(x)\triangleq\lnot``\mathrm{self}" \wedge ``0\leq\mathrm{others}\leq x";
&
\mathrm{W_{K}}(y)\triangleq``\mathrm{self}" \wedge ``0\leq\mathrm{others}\leq y";
\\
\,\,\mathrm{C_{D}}(x)\triangleq\lnot``\mathrm{ self}" \wedge ``1\leq\mathrm{others}\leq x";
&
\mathrm{W_{D}}(y)\triangleq``\mathrm{self}" \wedge ``0\leq\mathrm{others}\leq y";
\\
\mathrm{F}_{\mathrm{K}}(x,y)\triangleq\mathrm{C_{K}}(x) \vee \mathrm{W_{K}}(y);
&\!\!\!\!\!\!\!\!\!\!\!\!\!\!\!\!\!\!\!\!\!
\mathrm{F}_{\mathrm{D}}(x,y)\triangleq\mathrm{C_{D}}(x) \vee \mathrm{W_{D}}(y);
\end{eqnarray*}
\vspace{-6mm}

\begin{elabeling}{00.0000.0000.0000}
\item [{\textbf{\textsc{Theorem~11}}}]
\emph{For the $\mathbf{K}\mathrm{sys}\cml*,1\cmr$ systems defined in Section 5:}
\end{elabeling}\vspace{-4mm}

\begin{enumerate}
\item
\emph{The class of frames $\mathscr{F}$ for which
$\mathrm{F}_{\mathrm{K}}(x,y)$ holds at every world $w\in W$
corresponds to ${\mathbf{S_K}}(x,y)$}.
\item
\emph{The class of frames $\mathscr{F}$ for which
$\mathrm{F}_{\mathrm{D}}(x,y)$ holds at every world $w\in W$
corresponds to ${\mathbf{S_D}}(x,y)$}.
\end{enumerate}\vspace{1mm}

\noindent\textbf{\textsc{Proof}}\qquad{}
We organize this proof in 3 steps, as follows:\vspace{1mm}

\noindent\emph{Step 1: Preliminaries.}
Recall the $\mathbf{K}\mathrm{sys}\cml*,1\cmr$ axioms from Section 5:\vspace{-5mm}

\setlength{\jot}{2mm}
\begin{eqnarray*}
\alpha_{\mathbf{S_{K}}(x,y)}
&\!\!=\!\!&
m_{n-1}\,\,\rightarrow
\!\!\sum_{\substack{\chi(\mu_j) \leq x \\ \,\varepsilon_{j,n-1}=0}}
\!\!\!\!\mu_j \>\>+ 
\!\!\sum_{\substack{\chi(\mu_j) \leq y+1 \\ \,\varepsilon_{j,n-1}=1}}
\!\!\!\!\mu_j\,\,\,\,\mathrm{(all~\mu_\mathit{j}~in~section~\mathit{n-1})}
\\
\alpha_{\mathbf{S_{D}}(x,y)}
&\!\!=\!\!&
\mathrm{D}\,\alpha_{\mathbf{S_{K}}(x,y)}
\,=\,
\lozenge\mathit{1}\,\alpha_{\mathbf{S_{K}}(x,y)}
\end{eqnarray*}

We first prove the simple special cases:\vspace{-2mm}

\begin{itemize}
\item
When $(x,y)=(0,-1)$, $\mathrm{C_{K}}(0)$ allows only frames
consisting of \emph{blind} worlds,
i.e.~worlds that see no others,
and $\mathrm{W_{K}}(-1)$ is always false,
so the frame condition for $\mathbf{S_K}(0,-1)$ corresponds to $\mathbf{Ver}$.
But so does 
$\alpha_{\mathbf{S_K}(0,-1)}
\approx
m_{n-1} \rightarrow m_{n-1} \prod_{i=0}^{n-1}\mbox{!}\lozenge m_i
\approx
m_{n-1} \rightarrow \mbox{!}\lozenge\sum_{i=0}^{n-1}m_i
\approx
\mbox{!}\lozenge\mathit{1}$.\vspace{0.5mm}

Also, both $\mathrm{C_{D}}(0)$ and $\mathrm{W_{D}}(-1)$ are always false,
hence the frame condition for ${\mathbf{S_D}}(0,-1)$ corresponds to $\mathbf{F}$,
where
$\alpha_{\mathbf{S_D}(0,-1)}
\approx
\lozenge\mathit{1}\alpha_{\mathbf{S_K}(0,-1)}
\approx
\mathit{0}$.

\item
When $(x,y)=(n-1,n-1)=(*,*)$, ${\mathbf{S_K}}(*,*)$ is $\mathbf{K}$, because $\mathrm{F_{K}}(*,*)$ is true and
$\alpha_{\mathbf{S_K}(*,*)}
\approx
m_{n-1} \rightarrow m_{n-1} \sum\mu_j
\approx \mathit{1}$ ($\mu_j$ in all sections here).
Also, $\mathrm{C_{D}}(*)$ does not allow blind worlds,
but other than that,
both $\mathrm{C_{D}}(*)$ and $\mathrm{W_{D}}(*)$
allow worlds to see any number of other worlds,
so $\mathrm{F_{D}}(*,*)$ corresponds to serial frames,
hence ${\mathbf{S_D}}(*,*)$ is $\mathbf{D}$,
with $\alpha_{\mathbf{S_D}(*,*)}\approx\mathrm{D}$.
\end{itemize}

\vspace{-2mm}But we can still have $0\leq x\leq n-1$, $0\leq y\leq n-1$, where $(x,y)\not=(*,*)$.

\pagebreak{}

Let $V:W\times\mathcal{F}(v,d)\rightarrow\{0,1\}$ be a valuation on $\mathscr{F}$,
the other notations being as before.
To understand the significance of
$\mathrm{F}_{\mathrm{K}}(x,y)$ and $\mathrm{F}_{\mathrm{D}}(x,y)$,
the following observations are critical:\vspace{-2mm}

\begin{itemize}
\item
At any $w$, the values $V(w, p_k)$, $0\leq k<v$,
determine a unique level~0 minterm $m_i$, $0\leq i<n$, with $V(w,m_i)=1$.
But all other level 0 mintems $m'$ are disjoint with $m_i$, so $V(w,m')=0$.
Conversely, when $V(w,m_i)=1$ for some $m_i$,
this uniquely determines all $V(w, p_k)$.
Overall, choosing $V$ at $w$
is the same as choosing an $m_i$ for which $V(w,m_i)=1$,
so we can say that this $m_i$ represents \emph{the valuation at $w$}.
\item
If $w$ sees up to $z$ other $w'_i$, $0\leq i < z \leq n$,
then there are up to $z$ such $w'_i$ where valuations are pairwise distinct,
i.e.~where $V(w'_i,m_i)=1$ for all $m_i$ distinct.
Consequently, $V(w, \nu_i)=1$ for up to $z$ modal factors $\nu_i$.
In this case,
to have $V(w, \mu_j)=1$ for a minterm $\mu_j$, we need $\chi(\mu_j)\leq z$.

\end{itemize}\vspace{-1mm}

\noindent\emph{Step 2: Proof of a).}
Let $\mathscr{F}$ be an arbitrary frame.
It is well known (see for example \cite{Hughes})
that what we need to show here is the following: the axiom
$\alpha_{\mathbf{S_{K}}(x,y)}$ is valid on $\mathscr{F}$
iff $\mathrm{F}_{\mathrm{K}}(x,y)$ holds at all worlds of $\mathscr{F}$.\vspace{1mm}
Let us rewrite:\vspace{-4.5mm}

\[
\alpha_{\mathbf{S_{K}}(x,y)}
\,\approx\,
m_{n-1}\,\rightarrow \beta(x) + \gamma(y+1)
\]

\noindent where $\beta$ and $\gamma$ are the corresponding sums of minterms from section $n-1$.

\emph{Sufficiency:}
Assume $\mathrm{F}_{\mathrm{K}}(x,y)$ holds for $\mathscr{F}$.
Fixing an arbitrary valuation function $V$ on $\mathscr{F}$,
we show that $\alpha_{\mathbf{S_{K}}(x,y)}$ is valid in the model $(W,R,V)$.

If $V(w, m_{n-1})=0$ then $\alpha_{\mathbf{S_{K}}(x,y)}$ is immediately valid at $w$.
Thus, we can assume that $V(w,m_{n-1})=1$,
and it will be enough to find $V(w,\mu_j)=1$
for a single minterm $\mu_j$ from either $\beta$ or $\gamma$.
Such a minterm always has the right prefix $m_{n-1}$,
so it all depends on the states of its modal factors.

\vspace{1mm}\emph{Sufficiency case 1:}~When $0\leq x\leq n-1$ and $0\leq y\leq n-2$.

\vspace{1mm}If $\mathrm{C_K}(x)$ holds, $V(w,\nu_i)=1$
for at most $x$ modal factors $\nu_i$.
When $\nu_{n-1}$ is \emph{not} among them,
then $V(w,\mu_j)=1$ for a suitable minterm from $\beta$,
where all $\chi(\mu_j)\leq x$.
When $\nu_{n-1}$ \emph{is} among them,
then $V(w,\mu_j)=1$ for a suitable minterm from $\gamma$,
where all $\chi(\mu_j)\leq x$ because $x \leq y+1$.

If $\mathrm{W_K}(y)$ holds, $V(w,\nu_i)=1$
for at most $y+1$ modal factors $\nu_i$ in total,
one of them certainly being $\nu_{n-1}$ (because $w$ is reflexive),
plus up to $y$ others.
Then $V(w,\mu_j)=1$ for a suitable minterm from $\gamma$,
where all $\chi(\mu_j)\leq y+1$.

\vspace{1mm}\emph{Sufficiency case 2:}~When $0\leq x\leq n-2$ and $y=n-1$.

\vspace{1mm}If $\mathrm{C_K}(x)$ holds,
the argument is the same as above.\vspace{0.4mm}

If $\mathrm{W_K}(y)=\mathrm{W_K}(*)$ holds,
then $V(w, \nu_{n-1})=1$ (because $w$ is reflexive),
and $\gamma$ minterms now have
up to the maximum $n=y+1$ modal factors $\nu_i$ in state~1,
so regardless of how many distinct valuations $w$ actually sees,
one can always find a suitable minterm in $\gamma$ for which $V(w,\mu_j)=1$.

\pagebreak{}

\vspace{2mm}\emph{Necessity:}
Assume $\mathrm{F}_{\mathrm{K}}(x,y)$ does not hold for $\mathscr{F}$,
i.e. assume that there is a world $w\in W$ for which
$\mathrm{\tilde{F}}_{\mathrm{K}}(x,y) \triangleq \lnot\mathrm{F}_{\mathrm{K}}(x,y)$ holds, where:\vspace{-7mm}

\setlength{\jot}{1mm}
\begin{eqnarray*}
\mathrm{\tilde{F}_K}(x,y)
&\!\!\!\approx\!\!\!&
\lnot(\mathrm{C_K}(x)\,\vee\,\mathrm{W_K}(y))
\:\approx\:
\lnot\mathrm{C_K}(x)\,\wedge\,\lnot\mathrm{W_K}(y)
\\
&\!\!\!\approx\!\!\!&
(``\mathrm{self}"\vee\lnot``0\leq\mathrm{others}\leq x")
\wedge
(\lnot``\mathrm{self}"\vee\lnot``0\leq\mathrm{others}\leq y")
\\
&\!\!\!\approx\!\!\!&
(\lnot``\mathrm{self}"\wedge\lnot``0\leq\mathrm{others}\leq x")
\vee
(``\mathrm{self}"\wedge\lnot``0\leq\mathrm{others}\leq y")
\\
&\!\!\!\approx\!\!\!&
\mathrm{\tilde{C}_K}(x)\,\vee\,\mathrm{\tilde{W}_K}(y)
\end{eqnarray*}

\vspace{-3mm}For such a frame we choose $V(w, m_{n-1})=1$
to make the antecedent in $\alpha_{\mathbf{S_{K}}(x,y)}$ true,
and we show that, in all cases, we can find suitable valuations at the other worlds
such that the consequent in $\alpha_{\mathbf{S_{K}}(x,y)}$ is false.

\vspace{1mm}\emph{Necessity case 1:}~When $0\leq x\leq n-2$ and $0\leq y\leq n-2$.

\vspace{1mm}If $\mathrm{\tilde{C}_K}(x)$ holds,
$w$ sees at least $x+1$ other $w'_i$.
For $0\leq i\leq x$ we choose $V(w'_i, m_i)=1$ (all distinct) 
and for $i>x$, $V(w'_i, m_0)=1$ (no new valuation).
Since $x<n-1$, $V(w,\nu_{n-1})=0$,
so $V(w,\mu_j)$ could be 1 only for a $\beta$ minterm with $\chi(\mu_j)=x+1$,
which is not the case,
since our $\beta$ has only $\chi(\mu_j)\leq x$.

If $\mathrm{\tilde{W}_K}(y)$ holds,
$w$ sees itself and at least $y+1$ other $w'_i$.
But $y+1\leq n-1$, so $w$ always sees at least $y+2$ total $w'_i$.
For $0\leq i\leq y$ we choose $V(w'_i, m_i)=1$ ($y+1$ distinct valuations) 
and for $i>y$, $V(w'_i, m_{n-1})=1$ (valuation $y+2$).
Then $V(w,\mu_j)$ could be 1 only for a $\gamma$ minterm with $\chi(\mu_j)=y+2$,
which again is not the case,
since our $\gamma$ has only $\chi(\mu_j)\leq y+1$.

\vspace{1mm}\emph{Necessity case 2:}~When $x=n-1$ or $y=n-1$ (but not both).

\vspace{1mm}For $x=n-1$,
$\mathrm{\tilde{C}_K}(x)=\mathrm{\tilde{C}_K}(*)$ is false.
For $y=n-1$,
$\mathrm{\tilde{W}_K}(x)=\mathrm{\tilde{W}_K}(*)$ is false.
Either way, if the other condition holds, the argument is as above.\vspace{1mm}

This completes the proof of a).\vspace{2mm}

\noindent\emph{Step 3: Proof of b).}
For the D-plane systems,
$\alpha_{\mathbf{S_{D}}(x,y)}=
\mathrm{D}\,\alpha_{\mathbf{S_{K}}(x,y)}$,
where axiom $\mathrm{D}\triangleq\lozenge\mathit{1}$
eliminates prime orbit $\mathrm{Vv}_{0}$ from any minmatrix.
But $\mathrm{Vv}_{0}$ minterms
are precisely those that make a formula valid at blind worlds.
Thus, conjunction with $\mathrm{D}$
makes the subtle difference between $\mathrm{C_{K}}(y)$ and $\mathrm{C_{D}}(y)$
that precludes frames with blind worlds,
and Theorem 11 holds for D-plane systems too.
\textsc{\scriptsize{\hfill{}$\blacksquare$}}\vspace{3mm}

Recall that $n=2^{v}$.
Consider $v$-bit binary registers,
which can only hold unsigned numbers from 0 to $2^{v}-1$;
larger numbers cause a condition called \emph{overflow},
in which the register represents an unknown number.
In this sense, worlds of $\mathbf{K}\mathrm{sys}\cml v,1\cmr$ frames
can be said to be able to ``count'' up to $2^{v}-1$ worlds total,
themselves included;
beyond that, it is ``any number'' for them.
Together with the fact that they can only count the worlds that they see directly,
this could be interpreted as
the descriptive limitation of the semantics of non-iterative modal logics.

\pagebreak{}

We end with a result on non-finitely-axiomatizable systems
that proves that we have found all non-iterative normal logics.

\vspace{-1mm}\begin{elabeling}{0000.0000.0000.0000}
\item [{\textbf{\textsc{Corollary~12}}}]
\emph{Let $\mathbf{S}$ be a non-finitely-axiomatizable unimodal normal modal logic system.
Then $\mathbf{S}$ cannot be axiomatized
by a (countable) set of non-iterative formulas.}\vspace{-0.1mm}
\end{elabeling}

\noindent\textbf{\textsc{Proof.}}\qquad{}
Let $\mathbf{S}$ be defined by an infinite
(but necessarily countable) set of axioms $\{\alpha_{i}\}$, all non-iterative.
Every minmatrix $\alpha_{i}$ is in some $\mathbf{K}[v,1]$,
so it equals or collapses to some CMM
that has coordinates $(x_i,y_i)$ in $\mathbf{K}\mathrm{sys}\cml*,1\cmr$.
From Figure 3, $\{x_i\}$ and $\{y_i\}$ must have lower bounds,
say $x_m$ and $y_m$ respectively.
Define $\mathbf{S'}$ as follows:
if all $(x_i,y_i)$ correspond to K-plane systems
then $\mathbf{S'}\triangleq\mathbf{S_K}(x_m,y_m)$,
otherwise $\mathbf{S'}\triangleq\mathbf{S_D}(x_m,y_m)$.

Now let $\alpha\triangleq\,_v^1\cml\mathbf{S'}\cmr$
in the context where $\mathbf{S'}$ occurs first
(and then $\alpha$ is equiprovable with $_{v'}^1\cml\mathbf{S'}\cmr$ for all $v'>v$).
By construction $\mathbf{S}$ proves $\alpha$.
But by CMM inclusion $\alpha$ proves any $\alpha_i$,
making all the latter redundant.
So $\mathbf{S'}$ is really $\mathbf{S}$ and is finitely-axiomatizable.
\textsc{\scriptsize{\hfill{}$\blacksquare$}}\vspace{3mm}

In fact, we conjecture that, for any fixed positive integer $d$,
if a system \emph{$\mathbf{S}$} is not finitely-axiomatizable,
then it cannot have only axioms of modal degree $\leq d$.
Which is to claim that any non-finitely-axiomatizable system
must have non-redundant axioms of arbitrarily large modal degree.


\vspace{-0.1mm}
\AuthorAdressEmail{Adrian Soncodi}
{Lecturer, University of Texas at Dallas}
{acs151130@utdallas.edu\\
soncodi@verizon.net}

\pagebreak{}

\fontsize{11pt}{13pt}\selectfont

\section*{Appendix A}

In this appendix we show the first few
non-iterative modal contexts, lattices and systems.
These concrete examples are intended to help the reader follow
the generic calculations from the main sections of our paper.

For the systems presented below, we have already shown a method to derive defining axioms.
But axioms that are sums of minterms,
or even $\alpha_{\mathbf{S_{K}}(x,y)}$ or $\alpha_{\mathbf{S_{D}}(x,y)}$,
produce formulas that are far from minimal.
In this appendix we also show axiomatizations using a number of optimized formulas
that are more practical for inferences.

We recall some facts established in \cite{Soncodi}.
Normal minterms are a special case of canonical minterms,
namely those that satisfy the axioms of $\mathbf{K}$:
$\lozenge(p+q)\leftrightarrow\lozenge p + \lozenge q$ and $!\lozenge\mathit{0}$.
Thus, in a \emph{normal} $\mathbf{E}[v,d]$ minterm,
the states of all $\nu_i=\lozenge\mu_i$ factors
uniquely determine the state of any other $\nu_j=\lozenge\phi_j$,
where $\mu_i$ and $\phi_i$ are the minterms and formulas from $\mathbf{E}[v,d-1]$ respectively.
In $\mathbf{E}$, we assigned \emph{labels} to minterms
based on the state tuple $(\lozenge\mathit{1},\lozenge\mathit{0})$,
namely label V for (0,0), C for (0,1), D for (1,0) and W for (1,1).
Prime orbits are invariant under prime substitutions, which also preserve labels.
So prime orbits of $\mathbf{E}[v,d]$ contexts inherit labels from their minterms.
It turns out that only some D-orbits and precisely one V-orbit have normal minterms, 
which explains our notations.
And of course, prime orbits are context-dependent, so we have to rename them per context. 

\paragraph{Context {K[0,1]}}
Here, $v=0$ and $n=1$.
This context has 1 DNF factor and 2 minterms.
Since there are no substitutions, every minterm is a distinct prime orbit,
with minmatrices
$\mathrm{V} = \lozenge\mathit{1}|\mathit{0}| \approx \mbox{!}\lozenge\mathit{1}$
and
$\mathrm{D} = \lozenge\mathit{1}|\mathit{1}| \approx \lozenge\mathit{1}$.
The lattice $\mathbf{K}\mathrm{sys}\cml0,1\cmr$
is shown in Figure \ref{fig:K01-systems}.\vspace{-2mm}

\begin{wrapfigure}{r}{6.5cm}
\begin{centering}
\includegraphics[scale=0.7]{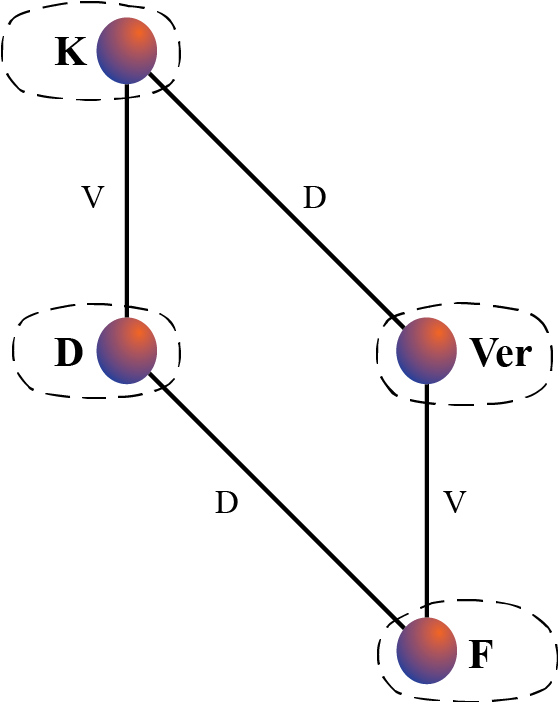}
\caption{The lattice $\mathbf{K}\mathrm{sys}\cml0,1\cmr$}\label{fig:K01-systems}
\end{centering}
\end{wrapfigure}

\setlength{\leftskip}{0.5cm}
\fontsize{9pt}{11pt}$
\\
\cml\mathbf{K}\cmr = \cml\mathbf{S_K}(*,*)\cmr = [\mathrm{V+D}]
\\
\hphantom{\cml}\mathbf{K}\hphantom{\cmr}
= \hphantom{\cml}\mathbf{K} \oplus ``\mathit{1}"
\\
\\
\cml\mathbf{D}\cmr = \cml\mathbf{S_D}(*,*)\cmr = [\mathrm{D}]
\\
\hphantom{\cml}\mathbf{D}\hphantom{\cmr}
= \hphantom{\cml}\mathbf{K} \oplus ``\lozenge\mathit{1}"
\\
\\
\cml\mathbf{Ver}\cmr =\cml\mathbf{S_K}(0,-1)\cmr = [\mathrm{V}]
\\
\hphantom{\cml}\mathbf{Ver}\hphantom{\cmr}
= \hphantom{\cml}\mathbf{K} \oplus ``\mbox{!}\lozenge\mathit{1}"
\\
\\
\cml\mathbf{F}\cmr = \cml\mathbf{S_D}(0,-1)\cmr = [\,]
\\
\hphantom{\cml}\mathbf{F}\hphantom{\cmr}
= \hphantom{\cml}\mathbf{K} \oplus ``\mathit{0}"
$

\quad\vspace{-0.5mm}\fontsize{11pt}{12pt}
\setlength{\leftskip}{0cm}

This is a limit context, because it misses the Dc and Dw prime orbits,
yet it still fits the generic pattern,
e.g. $\mathbf{S_D}(0,0) = \mathbf{S_D}(*,*)$.

\pagebreak{}

\paragraph{Context {K[1,1]}}
Here, $v=1$ and $n=2$.
This context has 3 DNF factors (1 Boolean plus 2 modal) and 8 minterms.
There are 2 prime substitutions, namely $p\mapsto p$ and $p\mapsto !p$,
which generate 4 prime orbits as follows:\vspace{3mm}

\noindent $\,\mathrm{Vv}\!=\!\begin{array}{c|cc|}
p & 1 & \!0\\
\hline \lozenge p & 0 & \!0\\
\lozenge!p & 0 & \!0
\end{array}$\,;\,
$\mathrm{Dd}\!=\!\begin{array}{c|cc|}
p & 1 & \!0\\
\hline \lozenge p & 1 & \!0\\
\lozenge!p & 0 & \!1
\end{array}$\,;\,
$\mathrm{Dc}\!=\!\begin{array}{c|cc|}
p & 1 & \!0\\
\hline \lozenge p & 0 & \!1\\
\lozenge!p & 1 & \!0
\end{array}$\,;\,
$\mathrm{Dw}\!=\!\begin{array}{c|cc|}
p & 1 & \!0\\
\hline \lozenge p & 1 & \!1\\
\lozenge!p & 1 & \!1
\end{array}$\,;\vspace{3mm}


The lattice $\mathbf{K}\mathrm{sys}\cml1,1\cmr$
is shown in Figure \ref{fig:K11-systems}.
The $\mathbf{K}[1,1]$-systems CMMs
and possible axiomatizations are as follows:

\begin{wrapfigure}{r}{6.0cm}
\begin{centering}
\vspace{2mm}\includegraphics[scale=0.6]{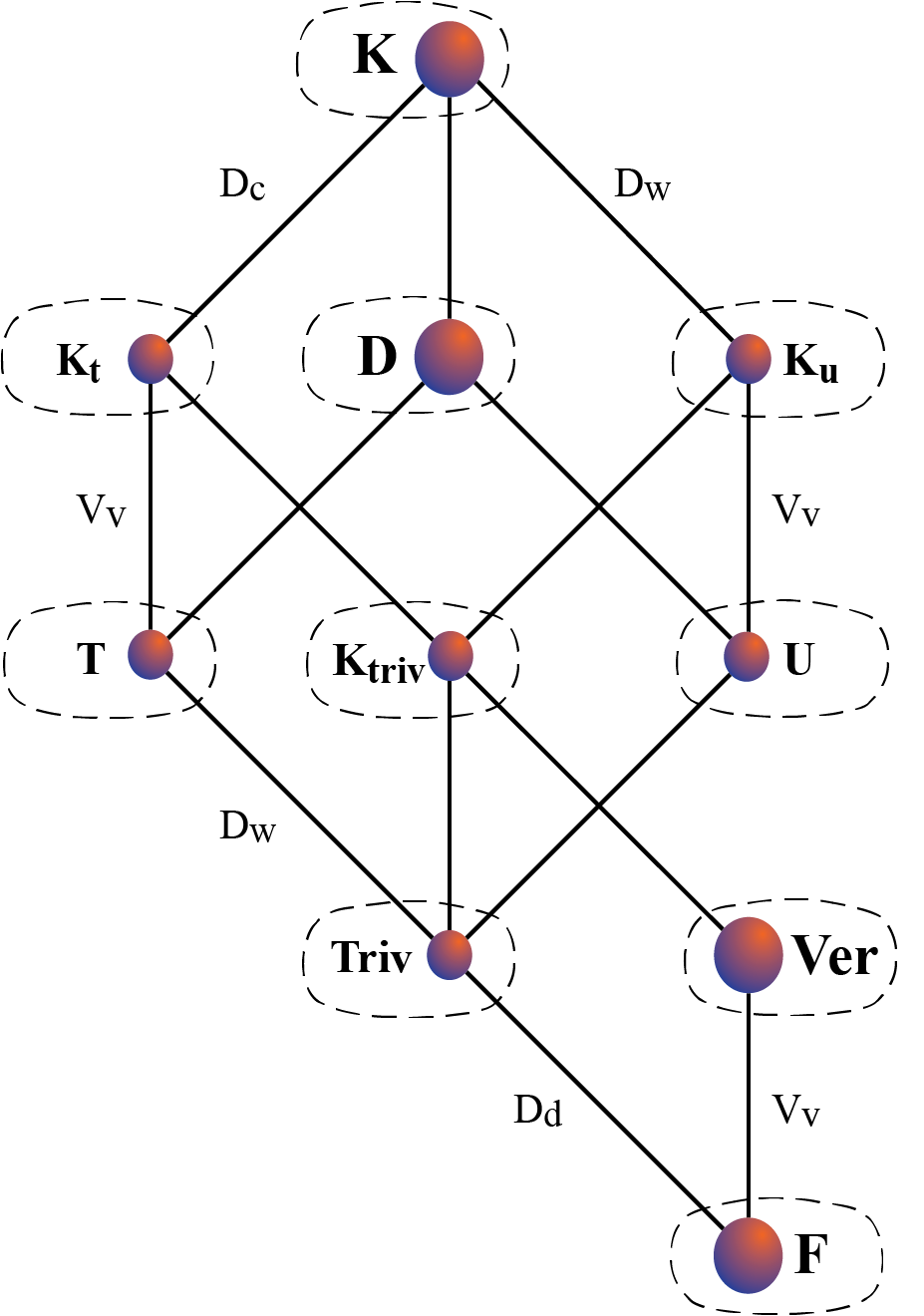}
~~~~~~~~\caption{The lattice $\mathbf{K}\mathrm{sys}\cml1,1\cmr$\newline
\newline
\newline
Note: Larger dots are systems\newline
$~~~~~~~~~\,$from predecessor contexts.
}\label{fig:K11-systems}
\end{centering}
\end{wrapfigure}

\setlength{\leftskip}{0.3cm}
\fontsize{9pt}{11pt}
\vspace{2mm}\noindent$
\cml\mathbf{K}\cmr = \cml\mathbf{S_K}(*,*)\cmr = [\mathrm{Vv+Dd+Dw+Dc}]
\\
\hphantom{\cml}\mathbf{K}\hphantom{\cmr}
= \hphantom{\cml}\mathbf{K} \oplus ``\mathit{1}"
\vspace{2mm}
\\
\cml\mathbf{D}\cmr = \cml\mathbf{S_D}(*,*)\cmr = [\mathrm{Dd+Dw+Dc}]
\\
\hphantom{\cml}\mathbf{D}\hphantom{\cmr}
= \mathbf{K} \oplus ``\lozenge\mathit{1}"
\\
\hphantom{\cml\mathbf{D}\cmr}
= \mathbf{K} \oplus ``\square p \rightarrow\lozenge p"
\vspace{2mm}
\\
\cml\mathbf{K_t}\cmr = \cml\mathbf{S_K}(0,*)\cmr = [\mathrm{Vv+Dd+Dw}]
\\
\hphantom{\cml}\mathbf{K_t}\hphantom{\cmr}
= \hphantom{\cml}\mathbf{K} \oplus ``p\rightarrow\lozenge p+\square p"
\vspace{2mm}
\\
\cml\mathbf{T}\cmr = \cml\mathbf{S_D}(0,*)\cmr = [\mathrm{Dd+Dw}]
\\
\hphantom{\cml}\mathbf{T}\hphantom{\cmr}
= \hphantom{\cml}\mathbf{K} \oplus ``p \rightarrow\lozenge p"
\\
\hphantom{\cml\mathbf{T}\cmr}
= \hphantom{\cml}\mathbf{K} \oplus ``\square p \rightarrow  p"
\vspace{2mm}
\\
\cml\mathbf{K_u}\cmr = \cml\mathbf{S_K}(0,1)\cmr = [\mathrm{Vv+Dd+Dc}]
\\
\hphantom{\cml}\mathbf{K_u}\hphantom{\cmr}
=
\hphantom{\cml}\mathbf{K} \oplus ``\lozenge p\rightarrow\square p"
\vspace{2mm}
\\
\cml\mathbf{U}\cmr = \cml\mathbf{S_D}(0,1)\cmr = [\mathrm{Dd+Dc}]
\\
\hphantom{\cml}\mathbf{U}\hphantom{\cmr}
=
\hphantom{\cml}\mathbf{K} \oplus ``\lozenge p \leftrightarrow\square p"
\vspace{2mm}
\\
\cml\mathbf{K_{triv}}\cmr = \cml\mathbf{S_K}(0,0)\cmr = [\mathrm{Vv+Dd}]
\\
\hphantom{\cml}\mathbf{K_{triv}}\hphantom{\cmr}
=\hphantom{\cml}\mathbf{K} \oplus ``\lozenge p \rightarrow  p"
\\
\hphantom{\cml\mathbf{K_{triv}}\cmr}
=\hphantom{\cml}\mathbf{K} \oplus ``p \rightarrow\square p"
\vspace{2mm}
\\
\cml\mathbf{Triv}\cmr = \cml\mathbf{S_D}(0,0)\cmr = [\mathrm{Dd}]
\\
\hphantom{\cml}\mathbf{Triv}\hphantom{\cmr}
=\hphantom{\cml}\mathbf{K} \oplus ``\lozenge p \leftrightarrow p"
\\
\hphantom{\cml\mathbf{Triv}\cmr}
=\hphantom{\cml}\mathbf{K} \oplus ``p \leftrightarrow\square p"
\vspace{2mm}
\\
\cml\mathbf{Ver}\cmr = \cml\mathbf{S_K}(0,-1)\cml = [\mathrm{Vv}]
\\
\hphantom{\cml}\mathbf{Ver}\hphantom{\cmr}
=
\hphantom{\cml}\mathbf{K} \oplus ``\mbox{!}\lozenge\mathit{1}"
\vspace{2mm}
=
\hphantom{\cml}\mathbf{K} \oplus ``\mathit{0}"
\vspace{2mm}
\\
\cml\mathbf{F}\cmr = \cml\mathbf{S_D}(0,-1)\cmr = [\,]
\\
\hphantom{\cml}\mathbf{F}\hphantom{\cmr}
=
\hphantom{\cml}\mathbf{K} \oplus ``\mathit{0}"
\vspace{2mm}
\\
$

\fontsize{11pt}{12pt}
\setlength{\leftskip}{0cm}
\vspace{-4mm}
This is one of the few contexts where calculations can be done easily by hand.
For example, one can check that other candidate CMMs
(like $\mathrm{[Dc+Dw]}$)
collapse, using the non-prime substitution $p\mapsto \mathit{0}$.

\pagebreak{}

\paragraph{Context {K[2,1]}}
Here, $v=2$ and $n=4$.
This context has 6 DNF factors (2 Boolean plus 4 modal) and 8 minterms.

Starting with this context, calculations by hand become quite tedious,
so we developed a software workbench to assist.
This is straightforward since,
after parsing formulas and converting them to minmatrices,
calculations become plain Boolean operations with bitsets.

Context $\mathbf{K}[2,1]$ has 24 prime substitutions $\varsigma(p,q)=(\varsigma_p(p,q),\varsigma_q(p,q))$,
where $(p,q)$ maps to one of the following:\vspace{1mm}
 
\noindent\begin{tabular}{ l l l l l l  } 
\,\,$(p, q)$, & $(!p, q)$  & $(p, p \leftrightarrow q)$,  & $(!p, p \leftrightarrow q)$, &
$(p \leftrightarrow q, p)$,  &  $(p \leftrightarrow !q, p)$,
\\
\,\,$(p, !q)$, & $(!p, !q)$, & $(p, p \leftrightarrow !q)$, & $(!p, p \leftrightarrow !q)$, &
$(p \leftrightarrow q, !p)$, &  $(p \leftrightarrow !q, !p)$,
\\
\,\,$(q, p)$,  & $(!q, p)$,  & $(q, p \leftrightarrow q)$,  & $(!q, p \leftrightarrow q)$,  &
$(p \leftrightarrow q, q)$,  &  $(p \leftrightarrow !q, q)$,
\\
\,\,$(q, !p)$, & $(!q, !p)$, & $(q, p \leftrightarrow !q)$, & $(!q, p \leftrightarrow !q)$, &
$(p \leftrightarrow q, !q)$, &  $(p \leftrightarrow !q, !q)$
\end{tabular}\vspace{2mm}

All the above belong to substitution class
$\mathcal{SC}_{0}(2,0) = \mathcal{S}_{p}(2,0) \subset \mathcal{S}(2,0)$.
Non-prime substitutions belong to one of 4 additional classes,
based on the dependencies between prime orbits that they reveal (see Section 4),
e.g.:\vspace{1mm}

\noindent\fontsize{10pt}{11pt}
\setlength{\leftskip}{0.15cm}
\begin{tabular}{ l l@{\hskip 0.1cm} r }
	$\mathcal{SC}_{1}(2,0)\ \mathrm{includes}\ (p,q)\mapsto
	(\mathit{0},\mathit{0}),\;(\mathit{0},\mathit{1}),
	\mathrm{\,etc.}$ & $\mathrm{(total}$ & $4)$
	\vspace{2pt}
	\\
	$\mathcal{SC}_{2}(2,0)\ \mathrm{includes}\ (p,q)\mapsto
	(\mathit{0}, pq),\;(p+q, !p!q),\;(\mathit{1}, p+q)
	\mathrm{\,etc.}$ & $\mathrm{(total}$ & $48)$
	\vspace{2pt}
	\\
	$\mathcal{SC}_{3}(2,0)\ \mathrm{includes}\ (p,q)\mapsto
	(p, \mathit{0}),\;(p, !p),\;(\mathit{0}, p\leftrightarrow q)
	\mathrm{\,etc.}$ & $\mathrm{(total}$ & $36)$
	\vspace{2pt}
	\\
	$\mathcal{SC}_{4}(2,0)\ \mathrm{includes}\ (p,q)\mapsto
	(p, pq),\;(q, p+q),\;(p\leftrightarrow q, pq)
	\mathrm{\,etc.}$ & $\mathrm{(total}$ & $144)$	
	\vspace{2pt}
\end{tabular}\vspace{1mm}

\noindent\setlength{\leftskip}{0cm}\fontsize{11pt}{12pt}where
$\mathrm{SC}_4(2,0)$ is the class of \emph{critical substitutions}.
\vspace{2mm}

The prime substitutions generate 8 prime orbits, namely:\vspace{-3mm}

\begin{center}
$\mathrm{Vvv,\>Ddd,\>Dcc1,\>Dcc2,\>Dcc3,\>Dww1,\>Dww2,\>Dww3}$
\end{center}\vspace{-3mm}

It is worth noting that promoting $\mathbf{K}[1,1]$ prime orbits does \emph{not} result in $\mathbf{K}[2,1]$ prime orbits.
Indeed the 4 promoted prime orbits from $\mathbf{K}[1,1]$ are disjoint minmatrices
and they still need to cover all $\mathbf{K}[2,1]$ minterms,
which are now partitioned between the 8 new prime orbits. We have:

\begin{itemize}
\item\vspace{-2mm}
$_1^1\mathrm{Vv}$ covers $_2^1\mathrm{Vvv}$ fully.
\item
$_1^1\mathrm{Dd}$ covers $_2^1\mathrm{Ddd}$ fully
and $_2^1\mathrm{Dcc1}$, $_2^1\mathrm{Dww1}$ partially.
\item
$_1^1\mathrm{Dc}$
covers $_2^1\mathrm{Dcc1}$ and $_2^1\mathrm{Dcc2}$ partially.
\item
$_1^1\mathrm{Dw}$
covers $_2^1\mathrm{Dww2}$, $_2^1\mathrm{Dww3}$, $_2^1\mathrm{Dcc3}$ fully
and $_2^1\mathrm{Dww1}$, $_2^1\mathrm{Dcc2}$ partially.
\end{itemize}

So to determine prime orbits, promotion does not work.
Rather, we use the following algorithm
(that works for canonical minterms from $\mathbf{E}$ too):
	
$-$ Pick any unassigned minterm, assign to a new prime orbit.

$-$ Perform all prime substitutions, assign the results to the same orbit.

$-$ Repeat as long as there are minterms unassigned to orbits.

\vspace{2mm}In the $\mathbf{K}[2,1]$ prime orbits shown below,
one can relate the highlighted bit patterns of the minterms
to the rules from Lemma 7 and Theorem 8:

\pagebreak{}

\noindent
$\quad\mathrm{Vvv}\,=\,\,\,\,
\begin{array}{c|cccc|}
    p          & 1 & 1 & 0 & 0 \\
    q          & 1 & 0 & 1 & 0 \\
\hline
\lozenge(pq)   & \!\colorbox{yellow}{\hspace{-0.15mm}0\hspace{0.15mm}}\!\! & 0 & 0 & 0 \\
\lozenge(p!q)  & 0 & \!\colorbox{yellow}{\hspace{-0.15mm}0\hspace{0.15mm}}\!\! & 0 & 0 \\
\lozenge(p!q)  & 0 & 0 & \!\colorbox{yellow}{\hspace{-0.15mm}0\hspace{0.15mm}}\!\! & 0 \\
\lozenge(!p!q) & 0 & 0 & 0 & \!\colorbox{yellow}{\hspace{-0.15mm}0\hspace{0.15mm}}\!\!
\end{array}$\,;\,
$\quad\,\,\mathrm{Ddd}\,=\,\,\,
\begin{array}{c|cccc|}
    p          & 1 & 1 & 0 & 0 \\
    q          & 1 & 0 & 1 & 0 \\
\hline
\lozenge(pq)   & \!\colorbox{yellow}{\hspace{-0.15mm}1\hspace{0.15mm}}\!\! & 0 & 0 & 0 \\
\lozenge(p!q)  & 0 & \!\colorbox{yellow}{\hspace{-0.15mm}1\hspace{0.15mm}}\!\! & 0 & 0 \\
\lozenge(p!q)  & 0 & 0 & \!\colorbox{yellow}{\hspace{-0.15mm}1\hspace{0.15mm}}\!\! & 0 \\
\lozenge(!p!q) & 0 & 0 & 0 & \!\colorbox{yellow}{\hspace{-0.15mm}1\hspace{0.15mm}}\!\!
\end{array}$\,;\vspace{2.9mm}

\noindent
$\quad\mathrm{Dcc1}\,=\,\,
\begin{array}{c|ccc:ccc:ccc:ccc|}
	p          & 1 & 1 & 1 & 1 & 1 & 1 & 0 & 0 & 0 & 0 & 0 & 0 \\
	q          & 1 & 1 & 1 & 0 & 0 & 0 & 1 & 1 & 1 & 0 & 0 & 0 \\
	\hline
	\lozenge(pq)   & \!\colorbox{yellow}{0\hspace{1.05cm}}\hspace{-1.15cm}
	                   & 0 & 0 & 1 & 0 & 0 & 1 & 0 & 0 & 1 & 0 & 0 \\
	\lozenge(p!q)  & 1 & 0 & 0 & \!\colorbox{yellow}{0\hspace{1.05cm}}\hspace{-1.15cm}
	                               & 0 & 0 & 0 & 1 & 0 & 0 & 1 & 0 \\
	\lozenge(p!q)  & 0 & 1 & 0 & 0 & 1 & 0 & \!\colorbox{yellow}{0\hspace{1.05cm}}\hspace{-1.15cm}
	                                           & 0 & 0 & 0 & 0 & 1 \\
	\lozenge(!p!q) & 0 & 0 & 1 & 0 & 0 & 1 & 0 & 0 & 1 & \!\colorbox{yellow}{0\hspace{1.05cm}}\hspace{-1.15cm}
	                                                       & 0 & 0
\end{array}$\,;\vspace{2.9mm}

\noindent
$\quad\mathrm{Dww1}\,=\!
\begin{array}{c|ccc:ccc:ccc:ccc|}
    p          & 1 & 1 & 1 & 1 & 1 & 1 & 0 & 0 & 0 & 0 & 0 & 0 \\
    q          & 1 & 1 & 1 & 0 & 0 & 0 & 1 & 1 & 1 & 0 & 0 & 0 \\
\hline
\lozenge(pq)   & \!\colorbox{yellow}{1\hspace{1.05cm}}\hspace{-1.15cm}
                   & 1 & 1 & 1 & 0 & 0 & 1 & 0 & 0 & 1 & 0 & 0 \\
\lozenge(p!q)  & 1 & 0 & 0 & \!\colorbox{yellow}{1\hspace{1.05cm}}\hspace{-1.15cm}
                               & 1 & 1 & 0 & 1 & 0 & 0 & 1 & 0 \\
\lozenge(p!q)  & 0 & 1 & 0 & 0 & 1 & 0 & \!\colorbox{yellow}{1\hspace{1.05cm}}\hspace{-1.15cm}
                                           & 1 & 1 & 0 & 0 & 1 \\
\lozenge(!p!q) & 0 & 0 & 1 & 0 & 0 & 1 & 0 & 0 & 1 & \!\colorbox{yellow}{1\hspace{1.05cm}}\hspace{-1.15cm}
                                                       & 1 & 1
\end{array}$\,;\vspace{2.9mm}

\noindent
$\quad\mathrm{Dcc2}\,=\,\,
\begin{array}{c|ccc:ccc:ccc:ccc|}
    p          & 1 & 1 & 1 & 1 & 1 & 1 & 0 & 0 & 0 & 0 & 0 & 0 \\
    q          & 1 & 1 & 1 & 0 & 0 & 0 & 1 & 1 & 1 & 0 & 0 & 0 \\
\hline
\lozenge(pq)   & \!\colorbox{yellow}{0\hspace{1.05cm}}\hspace{-1.15cm}
                   & 0 & 0 & 1 & 1 & 0 & 1 & 1 & 0 & 1 & 1 & 0 \\
\lozenge(p!q)  & 1 & 1 & 0 & \!\colorbox{yellow}{0\hspace{1.05cm}}\hspace{-1.15cm}
                               & 0 & 0 & 1 & 0 & 1 & 1 & 0 & 1 \\
\lozenge(p!q)  & 1 & 0 & 1 & 1 & 0 & 1 & \!\colorbox{yellow}{0\hspace{1.05cm}}\hspace{-1.15cm}
                                           & 0 & 0 & 0 & 1 & 1 \\
\lozenge(!p!q) & 0 & 1 & 1 & 0 & 1 & 1 & 0 & 1 & 1 & \!\colorbox{yellow}{0\hspace{1.05cm}}\hspace{-1.15cm}
                                                       & 0 & 0
\end{array}$\,;\vspace{2.9mm}

\noindent
$\quad\mathrm{Dww2}\,=\!
\begin{array}{c|ccc:ccc:ccc:ccc|}
    p          & 1 & 1 & 1 & 1 & 1 & 1 & 0 & 0 & 0 & 0 & 0 & 0 \\
    q          & 1 & 1 & 1 & 0 & 0 & 0 & 1 & 1 & 1 & 0 & 0 & 0 \\
\hline
\lozenge(pq)   & \!\colorbox{yellow}{1\hspace{1.05cm}}\hspace{-1.15cm}
                   & 1 & 1 & 1 & 1 & 0 & 1 & 1 & 0 & 1 & 1 & 0 \\
\lozenge(p!q)  & 1 & 1 & 0 & \!\colorbox{yellow}{1\hspace{1.05cm}}\hspace{-1.15cm}
                               & 1 & 1 & 1 & 0 & 1 & 1 & 0 & 1 \\
\lozenge(p!q)  & 1 & 0 & 1 & 1 & 0 & 1 & \!\colorbox{yellow}{1\hspace{1.05cm}}\hspace{-1.15cm}
                                           & 1 & 1 & 0 & 1 & 1 \\
\lozenge(!p!q) & 0 & 1 & 1 & 0 & 1 & 1 & 0 & 1 & 1 & \!\colorbox{yellow}{1\hspace{1.05cm}}\hspace{-1.15cm}
                                                       & 1 & 1
\end{array}$\,;\vspace{2.9mm}

\noindent
$\quad\mathrm{Dcc3}\,=\,\,
\begin{array}{c|cccc|}
	p          & 1 & 1 & 0 & 0 \\
	q          & 1 & 0 & 1 & 0 \\
\hline
\lozenge(pq)   & \!\colorbox{yellow}{\hspace{-0.15mm}0\hspace{0.15mm}}\!\! & 1 & 1 & 1 \\
\lozenge(p!q)  & 1 & \!\colorbox{yellow}{\hspace{-0.15mm}0\hspace{0.15mm}}\!\! & 1 & 1 \\
\lozenge(p!q)  & 1 & 1 & \!\colorbox{yellow}{\hspace{-0.15mm}0\hspace{0.15mm}}\!\! & 1 \\
\lozenge(!p!q) & 1 & 1 & 1 & \!\colorbox{yellow}{\hspace{-0.15mm}0\hspace{0.15mm}}\!\!
\end{array}$\,;\,
$\quad\!\mathrm{Dww3}\,=\,\,
\begin{array}{c|cccc|}
	p          & 1 & 1 & 0 & 0 \\
	q          & 1 & 0 & 1 & 0 \\
\hline
\lozenge(pq)   & \!\colorbox{yellow}{\hspace{-0.15mm}1\hspace{0.15mm}}\!\! & 1 & 1 & 1 \\
\lozenge(p!q)  & 1 & \!\colorbox{yellow}{\hspace{-0.15mm}1\hspace{0.15mm}}\!\! & 1 & 1 \\
\lozenge(p!q)  & 1 & 1 & \!\colorbox{yellow}{\hspace{-0.15mm}1\hspace{0.15mm}}\!\! & 1 \\
\lozenge(!p!q) & 1 & 1 & 1 & \!\colorbox{yellow}{\hspace{-0.15mm}1\hspace{0.15mm}}\!\!
\end{array}$\,;

\pagebreak{}

$\mathbf{K}\mathrm{sys}\cml2,1\cmr$ has 28 systems total:
the 10 from $\mathbf{K}\mathrm{sys}\cml1,1\cmr$ plus 18 new ones.
Here, one can already see how the systems from predecessor contexts (larger dots)
will preserve their position in the aggregate $\mathbf{K}\mathrm{sys}\cml*,1\cmr$.

\begin{figure}[H]
\begin{centering}
\includegraphics[scale=0.44]{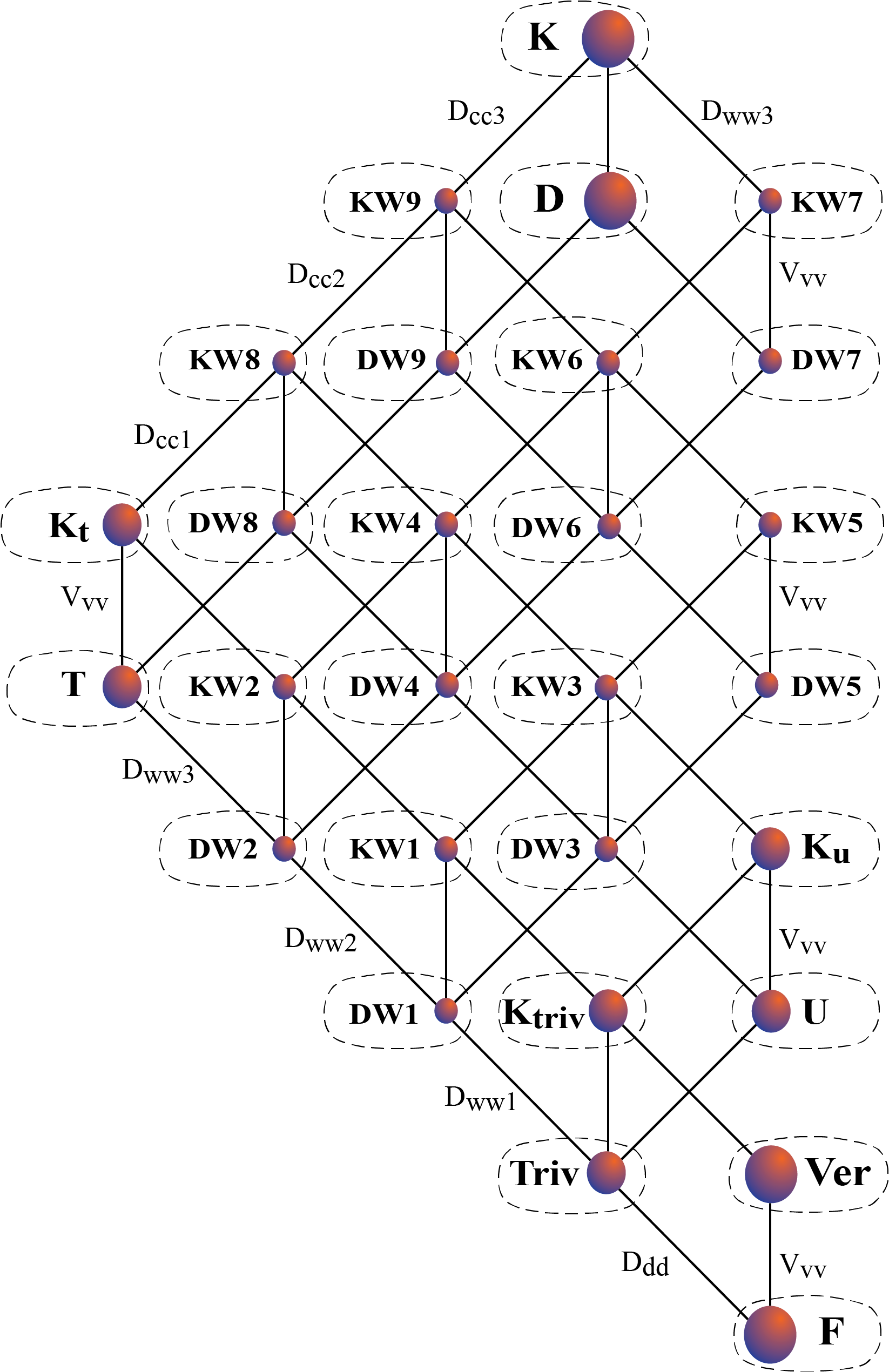}
\par\end{centering}
	
\caption{\label{fig:K21-systems}The lattice $\mathbf{K}\mathrm{sys}\cml2,1\cmr$}
\end{figure}
\vspace{-2mm}

\fontsize{11pt}{12pt}
The proposed axioms for $\mathbf{K}[2,1]$-systems
were found using our software workbench.
Although not necessarily minimal/optimal formulas,
they are often shorter than $\alpha_{\mathbf{S_K}(x,y)}$ and $\alpha_{\mathbf{S_D}(x,y)}$.
Note that, except for $\mathbf{Ver}$ and $\mathbf{F}$,
the remaining systems can be axiomatized using positive formulas.

\setlength{\leftskip}{0.3cm}
\fontsize{9pt}{11pt}
\vspace{1mm}\noindent
$
\cml\mathbf{K}\cmr = \cml\mathbf{S_K}(*,*)\cmr =
\mathrm{[Vvv+Ddd+Dcc1+Dww1+Dcc2+Dww2+Dcc3+Dww3]}
\\
\hphantom{\cml}\mathbf{K}\,
\,=\,
\mathbf{K} \oplus ``\mathit{1}"
\vspace{1.7mm}
\\
\cml\mathbf{D}\cmr = \cml\mathbf{S_D}(*,*)\cmr =
\mathrm{[Ddd+Dcc1+Dww1+Dcc2+Dww2+Dcc3+Dww3]}
\\
\hphantom{\cml}\mathbf{D}\,
\,=\,
\mathbf{K} \oplus \mathrm{D}
\,=\,
\mathbf{K} \oplus ``\square(pq) \rightarrow\lozenge(pq)"
\vspace{1.7mm}
\\
\cml\mathbf{KW9}\cmr =\cml\mathbf{S_K}(2,*)\cmr =
\mathrm{[Vvv+Ddd+Dcc1+Dww1+Dcc2+Dww2+Dww3]}
\\
\hphantom{\cml}\mathbf{KW9}\,
\,=\,
\mathbf{K} \oplus \mathrm{W9}
\,=\,
\mathbf{K} \oplus
``pq\lozenge p\lozenge q\rightarrow \lozenge(pq)+\square(p+q)"
\\
\hphantom{\cml\mathbf{KW9}\,=\,\mathbf{K} \oplus \mathrm{W9}\,}
\,=\,
\mathbf{K} \oplus
``pq\rightarrow \lozenge(pq)+\square(p\rightarrow q)+\square(q\rightarrow p)+\square(p+q)"
\vspace{1.7mm}
\\
\cml\mathbf{DW9}\cmr =\cml\mathbf{S_D}(2,*)\cmr =
\mathrm{[Ddd+Dcc1+Dww1+Dcc2+Dww2+Dww3]}
\\
\hphantom{\cml}\mathbf{DW9}\,
\,=\,
\mathbf{D} \oplus \mathrm{W9}
\,=\,
\mathbf{D} \oplus
``pq\lozenge p\lozenge q\rightarrow \lozenge(pq)+\square(p+q)"
\vspace{1.7mm}
\\
\cml\mathbf{KW8}\cmr =\cml\mathbf{S_K}(1,*)\cmr =
\mathrm{[Vvv+Ddd+Dcc1+Dww1+Dww2+Dww3]}
\\
\hphantom{\cml}\mathbf{KW8}\,
\,=\,
\mathbf{K} \oplus \mathrm{W8}
\,=\,
\mathbf{K} \oplus
``pq\lozenge p\lozenge q\rightarrow\lozenge(pq)"
\\
\hphantom{\cml\mathbf{KW8}\,=\,\mathbf{K} \oplus \mathrm{W8}\,}
\,=\,
\mathbf{K} \oplus
``p\lozenge q\rightarrow\lozenge p+\square q"
\\
\hphantom{\cml\mathbf{KW8}\,=\,\mathbf{K} \oplus \mathrm{W8}\,}
\,=\,
\mathbf{K} \oplus
``pq\lozenge p\lozenge (p\rightarrow q)\rightarrow\lozenge q"
\\
\hphantom{\cml\mathbf{KW8}\,=\,\mathbf{K} \oplus \mathrm{W8}\,}
\,=\,
\mathbf{K} \oplus
``pq\lozenge q\rightarrow\lozenge(pq)+\square(p+q)"
\\
\hphantom{\cml\mathbf{KW8}\,=\,\mathbf{K} \oplus \mathrm{W8}\,}
\,=\,
\mathbf{K} \oplus
``pq\lozenge(p\rightarrow q)\lozenge(q\rightarrow p)\rightarrow\lozenge(p\leftrightarrow q)"
\\
\hphantom{\cml\mathbf{KW8}\,=\,\mathbf{K} \oplus \mathrm{W8}\,}
\,=\,
\mathbf{K} \oplus
``pq\rightarrow\lozenge(pq)+\square(p\rightarrow q)+\square(q\rightarrow p)"
\vspace{1.7mm}
\\
\cml\mathbf{DW8}\cmr =\cml\mathbf{S_D}(1,*)\cmr =
\mathrm{[Ddd+Dcc1+Dww1+Dww2+Dww3]}
\\
\hphantom{\cml}\mathbf{DW8}\,
\,=\,
\mathbf{D} \oplus \mathrm{W8}
\,=\,
\mathbf{D} \oplus
``pq\lozenge p\lozenge q\rightarrow\lozenge(pq)"
\\
\hphantom{\cml\mathbf{DW8}\,=\,\mathbf{K} \oplus \mathrm{W8}\,}
\,=\,
\mathbf{K} \oplus
``pq\rightarrow \lozenge(pq)+(\lozenge p\leftrightarrow \square p)"
\\
\hphantom{\cml\mathbf{DW8}\,=\,\mathbf{K} \oplus \mathrm{W8}\,}
\,=\,
\mathbf{K} \oplus
``pq(\square p\leftrightarrow\square q)\rightarrow\lozenge(p\leftrightarrow q)"
\\
\hphantom{\cml\mathbf{DW8}\,=\,\mathbf{K} \oplus \mathrm{W8}\,}
\,=\,
\mathbf{K} \oplus
``pq\square(p\leftrightarrow q)\rightarrow(\lozenge p\leftrightarrow\lozenge q)"
\vspace{1.7mm}
\\
\cml\mathbf{K_t}\cmr =\cml\mathbf{S_K}(0,*)\cmr =
\mathrm{[Vvv+Ddd+Dww1+Dww2+Dww3]}
\\
\hphantom{\cml}\mathbf{K_t}\,
\,=\,
\mathbf{K} \oplus \mathrm{K_t}
\,=\,
\mathbf{K} \oplus
``pq\rightarrow\lozenge(pq)+\square(pq)"
\vspace{1.7mm}
\\
\cml\mathbf{T}\cmr =\cml\mathbf{S_D}(0,*)\cmr =
\mathrm{[Ddd+Dww1+Dww2+Dww3]}
\\
\hphantom{\cml}\mathbf{T}\,
\,=\,
\mathbf{K} \oplus \mathrm{T}
\,=\,
\mathbf{K} \oplus
``\square(pq)\rightarrow pq"
\vspace{1.7mm}
\\
\cml\mathbf{KW7}\cmr =\cml\mathbf{S_K}(3,2)\cmr =
\mathrm{[Vvv+Ddd+Dcc1+Dww1+Dcc2+Dww2+Dcc3]}
\\
\hphantom{\cml}\mathbf{KW7}\,
\,=\,
\mathbf{K} \oplus \mathrm{W7}
\,=\,
\mathbf{K} \oplus
``\lozenge(pq)\rightarrow\square(p\rightarrow q)+\square(q\rightarrow p)+\square(p+q)"
\vspace{1.7mm}
\\
\cml\mathbf{DW7}\cmr = \cml\mathbf{S_D}(3,2)\cmr =
\mathrm{[Ddd+Dcc1+Dww1+Dcc2+Dww2+Dcc3]}
\\
\hphantom{\cml}\mathbf{DW7}\,
\,=\,
\mathbf{D} \oplus \mathrm{W7}
\,=\,
\mathbf{D} \oplus
``\lozenge(pq)\rightarrow\square(p\rightarrow q)+\square(q\rightarrow p)+\square(p+q)"
\vspace{1.7mm}
\\
\cml\mathbf{KW6}\cmr =\cml\mathbf{S_K}(2,2)\cmr =
\mathrm{[Vvv+Ddd+Dcc1+Dww1+Dcc2+Dww2]}
\\
\hphantom{\cml}\mathbf{KW6}\,
\,=\,
\mathbf{K} \oplus \mathrm{W6}
\,=\,
\mathbf{K} \oplus
``pq\rightarrow\square(p\rightarrow q)+\square(q\rightarrow p)+\square(p+q)"
\vspace{1.7mm}
\\
\cml\mathbf{DW6}\cmr =\cml\mathbf{S_D}(2,2)\cmr =
\mathrm{[Ddd+Dcc1+Dww1+Dcc2+Dww2]}
\\
\hphantom{\cml}\mathbf{DW6}\,
\,=\,
\mathbf{D} \oplus \mathrm{W6}
\,=\,
\mathbf{D} \oplus
``pq\rightarrow\square(p\rightarrow q)+\square(q\rightarrow p)+\square(p+q)"
\vspace{1.7mm}
\\
\cml\mathbf{KW5}\cmr =\cml\mathbf{S_K}(2,1)\cmr =
\mathrm{[Vvv+Ddd+Dcc1+Dww1+Dcc2]}
\\
\hphantom{\cml}\mathbf{KW5}\,
\,=\,
\mathbf{K} \oplus \mathrm{W5}
\,=\,
\mathbf{K} \oplus
``\lozenge p\lozenge q\rightarrow\lozenge(pq)+\square(p+q)"
\\
\hphantom{\cml\mathbf{KW5}\,=\,\mathbf{K} \oplus \mathrm{W5}\,}
\,=\,
\mathbf{K} \oplus
``\lozenge(pq)\rightarrow\square(p\rightarrow q)+\square(q\rightarrow p)"
\\
\hphantom{\cml\mathbf{KW5}\,=\,\mathbf{K} \oplus \mathrm{W5}\,}
\,=\,
\mathbf{K} \oplus
``\lozenge(pq)(\square p\rightarrow\square q)\rightarrow\square(p\rightarrow q)"
\\
\hphantom{\cml\mathbf{KW5}\,=\,\mathbf{K} \oplus \mathrm{W5}\,}
\,=\,
\mathbf{K} \oplus
``\lozenge(pq)\square(p+q)\rightarrow\lozenge p\square q+\lozenge q\square p"
\\
\hphantom{\cml\mathbf{KW5}\,=\,\mathbf{K} \oplus \mathrm{W5}\,}
\,=\,
\mathbf{K} \oplus
``\square(p\rightarrow q)+\square(q\rightarrow p)+\square(p+q)"
\vspace{1.7mm}
\\
\cml\mathbf{DW5}\cmr =\cml\mathbf{S_D}(2,1)\cmr =
\mathrm{[Ddd+Dcc1+Dww1+Dcc2]}
\\
\hphantom{\cml}\mathbf{DW5}\,
\,=\,
\mathbf{D} \oplus \mathrm{W5}
\,=\,
\mathbf{D} \oplus
``\lozenge p\lozenge q\rightarrow\lozenge(pq)+\square(p+q)"
\vspace{1.7mm}
\\
\cml\mathbf{KW4}\cmr =\cml\mathbf{S_K}(1,2)\cmr =
\mathrm{[Vvv+Ddd+Dcc1+Dww1+Dww2]}
\\
\hphantom{\cml}\mathbf{KW4}\,
\,=\,
\mathbf{K} \oplus \mathrm{W4}
\,=\,
\mathbf{K} \oplus
``pq\rightarrow\square(p\rightarrow q)+\square(q\rightarrow p)+\lozenge(pq)\square(p+q)"
\\
\hphantom{\cml\mathbf{KW4}\,=\,\mathbf{K} \oplus \mathrm{W4}\,}
\,=\,
\mathbf{K} \oplus
``pq\lozenge p\rightarrow\lozenge q(\square(q\rightarrow p)+\square(p\rightarrow q))+\square(p+q)"
\vspace{1.7mm}
\\
\cml\mathbf{DW4}\cmr =\cml\mathbf{S_D}(1,2)\cmr =
\mathrm{[Ddd+Dcc1+Dww1+Dww2]}
\\
\hphantom{\cml}\mathbf{DW4}\,
\,=\,
\mathbf{D} \oplus \mathrm{W4}
\,=\,
\mathbf{D} \oplus
``pq\rightarrow\square(p\rightarrow q)+\square(q\rightarrow p)+\lozenge(pq)\square(p+q)"
\vspace{1.7mm}
\\
\cml\mathbf{KW3}\cmr =\cml\mathbf{S_K}(1,1)\cmr =
\mathrm{[Vvv+Ddd+Dcc1+Dww1]}
\\
\hphantom{\cml}\mathbf{KW3}\,
\,=\,
\mathbf{K} \oplus \mathrm{W3}
\,=\,
\mathbf{K} \oplus
``pq\rightarrow \square(p\rightarrow q)+\square(q\rightarrow p)"
\\
\hphantom{\cml\mathbf{KW3}\,=\,\mathbf{K} \oplus \mathrm{W3}\,}
\,=\,
\mathbf{K} \oplus
``pq\rightarrow \square(p\leftrightarrow q)+\square(p+q)"
\\
\hphantom{\cml\mathbf{KW3}\,=\,\mathbf{K} \oplus \mathrm{W3}\,}
\,=\,
\mathbf{K} \oplus
``pq\square(p+q)\rightarrow\square p+\square q"
\\
\hphantom{\cml\mathbf{KW3}\,=\,\mathbf{K} \oplus \mathrm{W3}\,}
\,=\,
\mathbf{K} \oplus
``pq(\square p\rightarrow\square q)\rightarrow\square(p\rightarrow q)"
\\
\hphantom{\cml\mathbf{KW3}\,=\,\mathbf{K} \oplus \mathrm{W3}\,}
\,=\,
\mathbf{K} \oplus
``p\lozenge q\rightarrow\lozenge(pq)+\square(p+q)"
\\
\hphantom{\cml\mathbf{KW3}\,=\,\mathbf{K} \oplus \mathrm{W3}\,}
\,=\,
\mathbf{K} \oplus
``p+q+\square(p\rightarrow q)+\square(q\rightarrow p)"
\vspace{1.7mm}
\\
\cml\mathbf{DW3}\cmr =\cml\mathbf{S_D}(1,1)\cmr =
\mathrm{[Ddd+Dcc1+Dww1]}
\\
\hphantom{\cml}\mathbf{DW3}\,
\,=\,
\mathbf{D} \oplus \mathrm{W3}
\,=\,
\mathbf{D} \oplus
``pq\rightarrow \square(p\rightarrow q)+\square(q\rightarrow p)"
\\
\hphantom{\cml\mathbf{DW3}\,=\,\mathbf{D} \oplus \mathrm{W3}\,}
\,=\,
\mathbf{K} \oplus
``pq\square(p+q)\rightarrow\lozenge p\square p+\lozenge q\square q"
\\
\hphantom{\cml\mathbf{DW3}\,=\,\mathbf{D} \oplus \mathrm{W3}\,}
\,=\,
\mathbf{K} \oplus
``pq\rightarrow\lozenge(p\rightarrow q)\square(p\rightarrow q)+\lozenge(q\rightarrow p)\square(q\rightarrow p)"
\vspace{1.7mm}
\\
\cml\mathbf{KW2}\cmr =\cml\mathbf{S_K}(0,2)\cmr =
\mathrm{[Vvv+Ddd+Dww1+Dww2]}
\\
\hphantom{\cml}\mathbf{KW2}\,
\,=\,
\mathbf{K} \oplus \mathrm{W2}
\,=\,
\mathbf{K} \oplus
``pq\rightarrow\lozenge(pq)(\square(p\rightarrow q)+\square(q\rightarrow p))+\square(p+q)"
\\
\hphantom{\cml\mathbf{KW2}\,=\,\mathbf{K} \oplus \mathrm{W2}\,}
\,=\,
\mathbf{K} \oplus
``pq\rightarrow\lozenge p\square(p\rightarrow q)+\lozenge q\square(q\rightarrow p)+\square(p+q)"
\vspace{1.7mm}
\\
\cml\mathbf{DW2}\cmr =\cml\mathbf{S_D}(0,2)\cmr =
\mathrm{[Ddd+Dww1+Dww2]}
\\
\hphantom{\cml}\mathbf{DW2}\,
\,=\,
\mathbf{D} \oplus \mathrm{W2}
\,=\,
\mathbf{D} \oplus
``pq\rightarrow\lozenge(pq)(\square(p\rightarrow q)+\square(q\rightarrow p))+\square(p+q)"
\\
\hphantom{\cml\mathbf{DW2}}\,
\,=\,\mathbf{T} \oplus \mathrm{W7}
\,=\,
\mathbf{T} \oplus
``pq\rightarrow\square(p\rightarrow q)+\square(q\rightarrow p)+\square(p+q)"
\\
\hphantom{\cml\mathbf{DW2}\,=\,\mathbf{D} \oplus \mathrm{W2}\,}
\,=\,
\mathbf{K} \oplus
``pq\rightarrow\lozenge(pq)(\square(p\rightarrow q)+\square(q\rightarrow p)+\square(p+q))"
\vspace{1.7mm}
\\
\cml\mathbf{K_u}\cmr =\cml\mathbf{S_K}(0,1)\cmr =
\mathrm{[Vvv+Ddd+Dcc1]}
\\
\hphantom{\cml}\mathbf{K_u}\,
\,=\,
\mathbf{K} \oplus \mathrm{K_u}
\,=\,
\mathbf{K} \oplus
``\lozenge p\rightarrow\square p"
\vspace{1.7mm}
\\
\cml\mathbf{U}\cmr =\cml\mathbf{S_D}(0,1)\cmr =
\mathrm{[Ddd+Dcc1]}
\\
\hphantom{\cml}\mathbf{U}\,
\,=\,
\mathbf{D} \oplus \mathrm{U}
\,=\,
\mathbf{D} \oplus
``\lozenge p\leftrightarrow\square p"
\vspace{1.7mm}
\\
\cml\mathbf{KW1}\cmr =\cml\mathbf{S_K}(0,1)\cmr =
\mathrm{[Vvv+Ddd+Dww1]}
\\
\hphantom{\cml}\mathbf{KW1}\,
\,=\,
\mathbf{K} \oplus \mathrm{W1}
\,=\,
\mathbf{K} \oplus
``pq\rightarrow\lozenge(pq)\square(p\leftrightarrow q)+\square(p+q)"
\\
\hphantom{\cml\mathbf{KW1}\,=\,\mathbf{K} \oplus \mathrm{W1}\,}
\,=\,
\mathbf{K} \oplus
``pq\rightarrow\square(p\leftrightarrow q)+\lozenge p\lozenge q\square(p+q)"
\vspace{1.7mm}
\\
\cml\mathbf{DW1}\cmr =\cml\mathbf{S_D}(0,1)\cmr =
\mathrm{[Ddd+Dww1]}
\\
\hphantom{\cml}\mathbf{DW1}\,
\,=\,
\mathbf{D} \oplus \mathrm{W1}
\,=\,
\mathbf{D} \oplus
``pq\rightarrow\lozenge(pq)\square(p\leftrightarrow q)+\square(p+q)"
\\
\hphantom{\cml\mathbf{DW1}}\,
\,=\,\mathbf{T} \oplus \mathrm{W3}
\,=\,
\mathbf{T} \oplus
``pq\rightarrow\square(p\rightarrow q)+\square(q\rightarrow p)"
\\
\hphantom{\cml\mathbf{DW1}\,=\,\mathbf{D} \oplus \mathrm{W1}\,}
\,=\,
\mathbf{K} \oplus
``pq\rightarrow\lozenge p\square(p\rightarrow q)+\lozenge q\square(q\rightarrow p)"
\vspace{1.7mm}
\\
\cml\mathbf{K_{triv}}\cmr =\cml\mathbf{S_K}(0,0)\cmr =
\mathrm{[Vvv+Ddd0]}
\\
\hphantom{\cml}\mathbf{K_{triv}}\,
\,=\,
\mathbf{K} \oplus \mathrm{K_{triv}}
\,=\,
\mathbf{K} \oplus
``p\rightarrow\square p"
\vspace{1.7mm}
\\
\cml\mathbf{Triv}\cmr =\cml\mathbf{S_D}(0,0)\cmr =
\mathrm{[Ddd0]}
\\
\hphantom{\cml}\mathbf{Triv}\,
\,=\,
\mathbf{D} \oplus \mathrm{Triv}
\,=\,
\mathbf{D} \oplus
``p\leftrightarrow\square p"
\vspace{1.7mm}
\\
\cml\mathbf{Ver}\cmr =\cml\mathbf{S_K}(0,-1)\cmr =
\mathrm{[Vvv+Ddd0]}
\\
\hphantom{\cml}\mathbf{Ver}\,
\,=\,
\mathbf{K} \oplus
``!\lozenge\mathit{1}"
\vspace{1.7mm}
\\
\cml\mathbf{F}\cmr =\cml\mathbf{S_D}(0,-1)\cmr =
\mathrm{[\,]}
\\
\hphantom{\cml}\mathbf{F}\,
\,=\,
\mathbf{K} \oplus
``\mathit{0}"
\vspace{1.7mm}
$

\setlength{\leftskip}{0.2cm}
\fontsize{11pt}{12pt}
\paragraph{Context {K[3,1]}}
Here, $v=3$ and $n=8$.
This context has 11 DNF factors (3 Boolean plus 8 modal) and 2048 minterms.
There are now 22 substitution classes,
starting with class $\mathcal{SC}_{00}(3,0) = \mathcal{S}_{p}(3,0)$ of  prime substitutions
and up to class $\mathcal{SC}_{21}(3,0)$ of critical substitutions.
Due to the large number ($n^n$) of level 0 substitutions, we only provide a few examples:
\vspace{2mm}

\noindent\fontsize{9pt}{11pt}
\setlength{\leftskip}{0cm}
\addtolength{\tabcolsep}{-0.15cm}
\begin{tabular}{ l l r }
	$\mathcal{SC}_{00}(3,0)\ \mathrm{includes}\ (p,q,r)\mapsto
	(p, q, q\!\leftrightarrow\!r),\;(p, q\!\leftrightarrow\!pr, r),
	\mathrm{\,etc.}$ & $\mathrm{(total}$ & $40,320)$
	\vspace{2pt}
	\\
	$\mathcal{SC}_{01}(3,0)\ \mathrm{includes}\ (p,q,r)\mapsto
	(\mathit{0}, \mathit{0}, \mathit{0}),\;(\mathit{0}, \mathit{1}, \mathit{1}),
	\mathrm{\,etc.}$ & $\mathrm{(total}$ & $8)$
	\vspace{2pt}
	\\
	$\mathcal{SC}_{02}(3,0)\ \mathrm{includes}\ (p,q,r)\mapsto
	(\mathit{0}, \mathit{0}, pqr),\;(\mathit{1}, pqr, pqr),
	\mathrm{\,etc.}$ & $\mathrm{(total}$ & $448)$
	\vspace{2pt}
	\\
	$\mathcal{SC}_{03}(3,0)\ \mathrm{includes}\ (p,q,r)\mapsto
	(\mathit{0}, \mathit{0}, pq),\;(\mathit{1}, pr, pr),
	\mathrm{\,etc.}$ & $\mathrm{(total}$ & $1,568)$	
	\vspace{2pt}
	\\
	$\mathcal{SC}_{04}(3,0)\ \mathrm{includes}\ (p,q,r)\mapsto
	(\mathit{0}, \mathit{0}, pq\!+\!r),\;(\mathit{0}, q+pr, \mathit{1}),
	\mathrm{\,etc.}$ & $\mathrm{(total}$ & $3,136)$
	\vspace{2pt}
	\\
	$\mathcal{SC}_{05}(3,0)\ \mathrm{includes}\ (p,q,r)\mapsto
	(\mathit{0}, \mathit{0}, r),\;(\mathit{0}, \mathit{1}, p\!\leftrightarrow\!q),
	\mathrm{\,etc.}$ & $\mathrm{(total}$ & $1,960)$
	\vspace{2pt}
	\\
	$\mathcal{SC}_{06}(3,0)\ \mathrm{includes}\ (p,q,r)\mapsto
	(\mathit{0}, pr, pqr),\;(\mathit{1}, p\!+\!q\!+\!r, pqr),
	\mathrm{\,etc.}$ & $\mathrm{(total}$ & $9,408)$
	\vspace{2pt}
	\\
	$\mathcal{SC}_{07}(3,0)\ \mathrm{includes}\ (p,q,r)\mapsto
	(\mathit{0}, p\!+\!r, pqr),\;(pqr, \mathit{1}, p+r),
	\mathrm{\,etc.}$ & $\mathrm{(total}$ & $56,448)$
	\vspace{2pt}
	\\
	$\mathcal{SC}_{08}(3,0)\ \mathrm{includes}\ (p,q,r)\mapsto
	(\mathit{0}, q, pqr),\;(\mathit{0}, q(p\!+\!r), q),
	\mathrm{\,etc.}$ & $\mathrm{(total}$ & $94,080)$
	\vspace{2pt}
	\\
	$\mathcal{SC}_{09}(3,0)\ \mathrm{includes}\ (p,q,r)\mapsto
	(\mathit{0}, pr, pq),\ (pqr, pq, pr),
	\mathrm{\,etc.}$ & $\mathrm{(total}$ & $94,080)$
	\vspace{2pt}
	\\
	$\mathcal{SC}_{10}(3,0)\ \mathrm{includes}\ (p,q,r)\mapsto
	(\mathit{0}, q, pq),\ (r(p\!\leftrightarrow\!q), r, r),
	\mathrm{\,etc.}$ & $\mathrm{(total}$ & $70,560)$
	\vspace{2pt}
	\\
	$\mathcal{SC}_{11}(3,0)\ \mathrm{includes}\ (p,q,r)\mapsto
	(\mathit{0}, pq, r(p\!+\!q)),\ (pqr, qr, q),
	\mathrm{\,etc.}$ & $\mathrm{(total}$ & $705,600)$
	\vspace{2pt}
	\\
	$\mathcal{SC}_{12}(3,0)\ \mathrm{includes}\ (p,q,r)\mapsto
	(\mathit{0}, pq, r\!+\!pq),\ (qr, qr, p\!+\!qr),
	\mathrm{\,etc.}$ & $\mathrm{(total}$ & $94,080)$
	\vspace{2pt}
	\\
	$\mathcal{SC}_{13}(3,0)\ \mathrm{includes}\ (p,q,r)\mapsto
	(\mathit{0}, q, pr),\ (qr, p, p),
	\mathrm{\,etc.}$ & $\mathrm{(total}$ & $470,400)$
	\vspace{2pt}
	\\
	$\mathcal{SC}_{14}(3,0)\ \mathrm{includes}\ (p,q,r)\mapsto
	(\mathit{0}, p(q+r), q),\ (pqr, q, q\!+\!r),
	\mathrm{\,etc.}$ & $\mathrm{(total}$ & $1,411,200)$
	\vspace{2pt}
	\\
	$\mathcal{SC}_{15}(3,0)\ \mathrm{includes}\ (p,q,r)\mapsto
	(\mathit{0}, p\!\leftrightarrow\!qr, r),\ (qr, r, q),
	\mathrm{\,etc.}$ & $\mathrm{(total}$ & $176,400)$
	\vspace{2pt}
	\\
	$\mathcal{SC}_{16}(3,0)\ \mathrm{includes}\ (p,q,r)\mapsto
	(p, pq, pr),\ (p\!+\!r, q\!+\!r, r),
	\mathrm{\,etc.}$ & $\mathrm{(total}$ & $470,400)$
	\vspace{2pt}
	\\
	$\mathcal{SC}_{17}(3,0)\ \mathrm{includes}\ (p,q,r)\mapsto
	(p, qr, r(p\!\leftrightarrow\!q)),\ (p, qr, p\!+\!q),
	\mathrm{\,etc.}$ & $\mathrm{(total}$ & $3,763,200)$
	\vspace{2pt}
	\\
	$\mathcal{SC}_{18}(3,0)\ \mathrm{includes}\ (p,q,r)\mapsto
	(p, q\!\leftrightarrow\!r, pqr),\ (p\!+\!qr, pr, r),
	\mathrm{\,etc.}$ & $\mathrm{(total}$ & $2,822,400)$
	\vspace{2pt}
	\\
	$\mathcal{SC}_{19}(3,0)\ \mathrm{includes}\ (p,q,r)\mapsto
	(p, pr, q(p\!+\!r)),\ (p(q\!+\!r), qr, r),
	\mathrm{\,etc.}$ & $\mathrm{(total}$ & $1,128,960)$
	\vspace{2pt}
	\\
	$\mathcal{SC}_{20}(3,0)\ \mathrm{includes}\ (p,q,r)\mapsto
	(pr, p\!\leftrightarrow\!r, q),\ (p, q\!\leftrightarrow\!r, p\!+\!q),
	\mathrm{\,etc.}$ & $\mathrm{(total}$ & $4,233,600)$
	\vspace{2pt}
	\\
	$\mathcal{SC}_{21}(3,0)\ \mathrm{includes}\ (p,q,r)\mapsto
	(p\!+\!qr, q, r),\ (p, p\!\leftrightarrow\!q, r\!+\!pq),
	\mathrm{\,etc.}$ & $\mathrm{(total}$ & $1,128,960)$
	\vspace{2pt}
\end{tabular}\vspace{2mm}

\setlength{\leftskip}{0cm}
\fontsize{11pt}{12pt}
\addtolength{\tabcolsep}{0.15cm}

The exact patterns that determine which substitutions are in which class are for further study.
Nevertheless, $\mathcal{SC}_{00}(3,0)$ and $\mathcal{SC}_{21}(3,0)$ are the most important classes.
The labeling of the others is not very relevant,
because (see Section 4) they determine prime orbit dependencies that are weaker
than those revealed by the critical substitutions.\vspace{2mm}

The prime substitutions determine 16 prime orbits, which conform to the patterns from Lemma 7 and Theorem 8, namely:\vspace{1mm}

\begin{tabular}{l l}
Vvvv and Dddd & have 8 minterms each, and\vspace{2pt}
\\
$\mathrm{Dccc}_i$ and $\mathrm{Dwww}_i$ & have $8{n-1 \choose i}$ minterms each, $1\leq i \leq 7$.\vspace{2pt}
\\
\end{tabular}

\pagebreak{}

The lattice $\mathbf{K}\mathrm{sys}\cml3,1\cmr$} 
is shown in Figure \ref{fig:K31-systems}.
It has 88 systems total: the 28 from $\mathbf{K}\mathrm{sys}\cml2,1\cmr$ plus 60 new ones.
Larger dots are systems from predecessor contexts.
Note that, as for the other contexts,
we present the Hasse diagram by CMM \emph{rank} (i.e. number of prime orbits)
rather than by CMM \emph{count} (i.e. number of minterms).
\vspace{2mm}

\begin{figure}[H]
\begin{centering}
\includegraphics[scale=0.71]{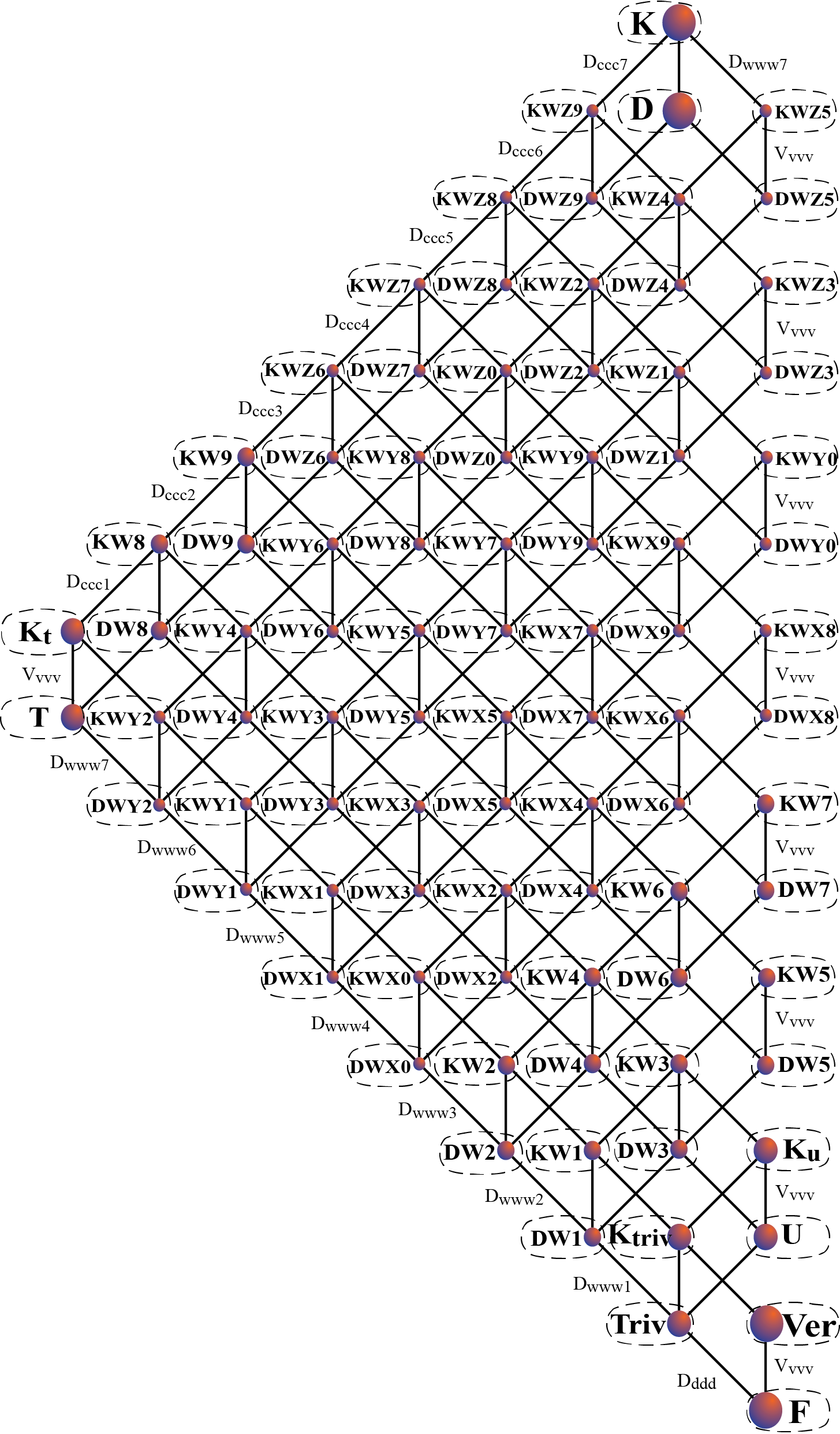}
\par\end{centering}
	
\caption{\label{fig:K31-systems}The lattice $\mathbf{K}\mathrm{sys}\cml3,1\cmr$}
\end{figure}

\pagebreak{}

The $\mathbf{K}\mathrm[3,1]$-systems are similar to those in the predecessors contexts.
We present them here for completness, but only the new K-plane systems.
For the corresponding D-plane system one just needs to add axiom D.

In their axioms we note the repeated occurrence of $\sum\square e_i$ terms,
which reflect the sums from  $\alpha_{\mathbf{S_K}(x,y)}$ and $\alpha_{\mathbf{S_D}(x,y)}$
(although in many cases these sums can be optimized).
Also, to write more compact formulas, we define the following notations
for \emph{cyclic sums} and \emph{cyclic products}:
\vspace{-7mm}

\begin{eqnarray*}
\sum^{\circ}e(p,q,r) & \triangleq & e(p,q,r)+e(q,r,p)+e(r,p,q)
\\
\prod^{\circ}e(p,q,r) & \triangleq & e(p,q,r)\,e(q,r,p)\,e(r,p,q)
\end{eqnarray*}

\vspace{-1.5mm}Then the systems can be described as follows:

\setlength{\leftskip}{0.1cm}
\fontsize{8pt}{11pt}
\vspace{1mm}\noindent
$
\cml\mathbf{KWZ9}\cmr =\cml\mathbf{S_K}(6,*)\cmr =
\mathrm{[Vvvv+Dddd+Dcc1+...+Dccc6+Dwww1+...+Dwww7]}
\\
\hphantom{\cml}\mathbf{KWZ9}\phantom{\cmr}
=\,
\mathbf{K} \oplus \mathrm{WZ9}
=\, \mathbf{K} \oplus
``pqr\rightarrow\lozenge(pqr)
+ \sum^\circ(\square(p\rightarrow q+r) + \square(qr\rightarrow p)) + \square(p+q+r)"
\vspace{1.7mm}
\\
\cml\mathbf{KWZ8}\cmr =\cml\mathbf{S_K}(5,*)\cmr =
\mathrm{[Vvvv+Dddd+Dcc1+...+Dccc5+Dwww1+...+Dwww7]}
\\
\hphantom{\cml}\mathbf{KWZ8}\hphantom{\cmr}=\,\mathbf{K} \oplus \mathrm{WZ8}
=\, \mathbf{K} \oplus
``pqr\rightarrow\lozenge(pqr)
+ \sum^\circ(\square(p\rightarrow q+r) + \square(qr\rightarrow p))" \hphantom{+ \square(p+q+r)"}
\\
\hphantom{\cml\mathbf{KWZ8}\cmr=\,\mathbf{K} \oplus \mathrm{WZ8}}
=\, \mathbf{K} \oplus
``pqr\rightarrow\lozenge(pqr)
+ \sum^\circ(\square(p+(q\leftrightarrow r)) + \square(qr\rightarrow p)) + \square(p+q+r)"
\vspace{1.7mm}
\\
\cml\mathbf{KWZ7}\cmr =\cml\mathbf{S_K}(4,*)\cmr =
\mathrm{[Vvvv+Dddd+Dcc1+...+Dccc4+Dwww1+...+Dwww7]}
\\
\hphantom{\cml}\mathbf{KWZ7}\hphantom{\cmr}=\,\mathbf{K} \oplus \mathrm{WZ7}
=\, \mathbf{K} \oplus
``pqr\rightarrow\lozenge(pqr)
+ \sum^\circ(\square(p+(q\leftrightarrow r)) + \square(qr\rightarrow p))"
\vspace{1.7mm}
\\
\cml\mathbf{KWZ6}\cmr =\cml\mathbf{S_K}(3,*)\cmr =
\mathrm{[Vvvv+Dddd+Dcc1+...+Dccc3+Dwww1+...+Dwww7]}
\\
\hphantom{\cml}\mathbf{KWZ6}\hphantom{\cmr}=\,\mathbf{K} \oplus \mathrm{WZ6}
=\, \mathbf{K} \oplus
``pqr\rightarrow\lozenge(pqr)
+ \sum^\circ\square(p\rightarrow q) + \square(p+q+r)"
\vspace{1.7mm}
\\
\cml\mathbf{KWZ5}\cmr =\cml\mathbf{S_K}(7,6)\cmr =
\mathrm{[Vvvv+Dddd+Dcc1+...+Dccc7+Dwww1+...+Dwww6]}
\\
\hphantom{\cml}\mathbf{KWZ5}\hphantom{\cmr}=\,\mathbf{K} \oplus \mathrm{WZ5}
=\, \mathbf{K} \oplus
``\lozenge(pqr)\rightarrow
\sum^\circ(\square(p\rightarrow q+r) + \square(qr\rightarrow p)) + \square(p+q+r)"
\vspace{1.7mm}
\\
\cml\mathbf{KWZ4}\cmr =\cml\mathbf{S_K}(6,6)\cmr =
\mathrm{[Vvvv+Dddd+Dcc1+...+Dccc6+Dwww1+...+Dwww6]}
\\
\hphantom{\cml}\mathbf{KWZ4}\hphantom{\cmr}=\,\mathbf{K} \oplus \mathrm{WZ4}
=\, \mathbf{K} \oplus
``pqr\rightarrow
\sum^\circ(\square(p\rightarrow q+r) + \square(qr\rightarrow p)) + \square(p+q+r)"
\vspace{1.7mm}
\\
\cml\mathbf{KWZ3}\cmr =\cml\mathbf{S_K}(6,5)\cmr =
\mathrm{[Vvvv+Dddd+Dcc1+...+Dccc6+Dwww1+...+Dwww5]}
\\
\hphantom{\cml}\mathbf{KWZ3}\hphantom{\cmr}=\,\mathbf{K} \oplus \mathrm{WZ3}
=\, \mathbf{K} \oplus
``\lozenge(pqr)\rightarrow
\sum^\circ(\square(p\rightarrow q+r) + \square(qr\rightarrow p))"
\\
\hphantom{\cml\mathbf{KWZ3}\hphantom{\cmr}=\,\mathbf{K} \oplus \mathrm{WZ3}}
=\, \mathbf{K} \oplus
``\sum^\circ(\square(p\rightarrow q+r) + \square(qr\rightarrow p)) +\square(p+q+r)"
\vspace{1.7mm}
\\
\cml\mathbf{KWZ2}\cmr =\cml\mathbf{S_K}(5,6)\cmr =
\mathrm{[Vvvv+Dddd+Dcc1+...+Dccc5+Dwww1+...+Dwww6]}
\\
\hphantom{\cml}\mathbf{KWZ2}\hphantom{\cmr}=\,\mathbf{K} \oplus \mathrm{WZ2}
=\, \mathbf{K} \oplus
``pqr\rightarrow\lozenge(pqr)\square(p+q+r)+
\sum^\circ(\square(p\rightarrow q+r) + \square(qr\rightarrow p))"
\vspace{1.7mm}
\\
\cml\mathbf{KWZ1}\cmr =\cml\mathbf{S_K}(5,5)\cmr =
\mathrm{[Vvvv+Dddd+Dcc1+...+Dccc5+Dwww1+...+Dwww5]}
\\
\hphantom{\cml}\mathbf{KWZ1}\hphantom{\cmr}=\,\mathbf{K} \oplus \mathrm{WZ1}
=\, \mathbf{K} \oplus
``pqr\rightarrow
\sum^\circ(\square(p\rightarrow q+r) + \square(qr\rightarrow p))"
\\
\hphantom{\cml\mathbf{KWZ1}\hphantom{\cmr}=\,\mathbf{K} \oplus \mathrm{WZ1}}
=\, \mathbf{K} \oplus
``pqr\rightarrow
\sum^\circ( \square(p+(q\leftrightarrow r)) +\square(qr\rightarrow p)) + \square(p+q+r)"
\vspace{1.7mm}
\\
\cml\mathbf{KWZ0}\cmr =\cml\mathbf{S_K}(4,6)\cmr =
\mathrm{[Vvvv+Dddd+Dcc1+...+Dccc4+Dwww1+...+Dwww6]}
\\
\hphantom{\cml}\mathbf{KWZ0}\hphantom{\cmr}=\,\mathbf{K} \oplus \mathrm{WZ0}
=\, \mathbf{K} \oplus
``pqr\rightarrow
\lozenge(pqr)(\sum^\circ\square(qr\rightarrow p) + \square(p+q+r))
\\
\phantom{\qquad\qquad\qquad\qquad\qquad\qquad\quad}
+\sum^\circ(\square(p+(q\leftrightarrow r)) + \square(qr\rightarrow p))"
\vspace{1.7mm}
\\
\cml\mathbf{KWY9}\cmr =\cml\mathbf{S_K}(4,5)\cmr =
\mathrm{[Vvvv+Dddd+Dcc1+...+Dccc4+Dwww1+...+Dwww5]}
\\
\hphantom{\cml}\mathbf{KWY9}\hphantom{\,}=\,\mathbf{K} \oplus \mathrm{WY9}
=\, \mathbf{K} \oplus
``pqr\rightarrow
\lozenge(pqr)\square(p+q+r)
+\sum^\circ(\square(p+(q\leftrightarrow r)) + \square(qr\rightarrow p))" 
\vspace{1.7mm}
\\
\cml\mathbf{KWY8}\cmr =\cml\mathbf{S_K}(3,6)\cmr =
\mathrm{[Vvvv+Dddd+Dcc1+...+Dccc3+Dwww1+...+Dwww6]}
\\
\hphantom{\cml}\mathbf{KWY8}\hphantom{\,}=\,\mathbf{K} \oplus \mathrm{WY8}
=\, \mathbf{K} \oplus
``pqr\rightarrow
\lozenge(pqr)\sum^\circ\square(p\rightarrow q+r) + \sum^\circ\square(qr\rightarrow p) + \square(p+q+r)" 
\\
\hphantom{\cml\mathbf{KWY8}\hphantom{\,}=\,\mathbf{K} \oplus \mathrm{WY8}}
=\, \mathbf{K} \oplus
``pqr\rightarrow
\lozenge(pqr)
\sum^\circ\square(qr\rightarrow p) + \sum^\circ\square(p\rightarrow q+r) + \square(p+q+r)" 
\vspace{1.7mm}
\\
\cml\mathbf{KWY7}\cmr =\cml\mathbf{S_K}(3,5)\cmr =
\mathrm{[Vvvv+Dddd+Dcc1+...+Dccc4+Dwww1+...+Dwww5]}
\\
\hphantom{\cml}\mathbf{KWY7}\hphantom{\,}=\,\mathbf{K} \oplus \mathrm{WY7}
=\, \mathbf{K} \oplus
``pqr\rightarrow
\lozenge(pqr)\sum^\circ\square(p\rightarrow q+r) + \sum^\circ(\square(qr\rightarrow p)+\square(q+r))" 
\vspace{1.7mm}
\\
\cml\mathbf{KWY6}\cmr =\cml\mathbf{S_K}(2,6)\cmr =
\mathrm{[Vvvv+Dddd+Dcc1+Dccc2+Dwww1+...+Dwww6]}
\\
\hphantom{\cml}\mathbf{KWY6}\hphantom{\cmr}=\,\mathbf{K} \oplus \mathrm{WY6}
=\, \mathbf{K} \oplus
``\prod^\circ \!p\lozenge p\rightarrow
\lozenge(pqr)
\sum^\circ\square(qr\rightarrow p) + \sum^\circ\square(p\rightarrow q+r)+\square(p+q+r)" 
\\
\hphantom{\cml\mathbf{KWY6}=\,\mathbf{K} \oplus \mathrm{WY6}}
=\, \mathbf{K} \oplus
``pqr\rightarrow
\lozenge(pqr)\sum^\circ(\square(p\rightarrow q+r) + \square(qr\rightarrow p))"
\\
\phantom{\qquad\qquad\qquad\qquad\qquad\qquad\quad}
+ \sum^\circ\square(p\leftrightarrow q) + \square(p+q+r)" 
\vspace{1.7mm}
\\
\cml\mathbf{KWY5}\cmr =\cml\mathbf{S_K}(2,5)\cmr =
\mathrm{[Vvvv+Dddd+Dcc1+Dccc2+Dwww1+...+Dwww5]}
\\
\hphantom{\cml}\mathbf{KWY5}\hphantom{\,}=\,\mathbf{K} \oplus \mathrm{WY5}
=\, \mathbf{K} \oplus
``pqr\rightarrow
\lozenge(pqr)\sum^\circ\square(p\rightarrow q+r) + \sum^\circ\square(qr\rightarrow p)" 
\\
\hphantom{\cml\mathbf{KWY5}=\,\mathbf{K} \oplus \mathrm{WY5}}
\,=\, \mathbf{K} \oplus
``pqr\rightarrow
\lozenge(pqr)\sum^\circ\square(qr\rightarrow p) + \sum^\circ\square(p\rightarrow q+r)"
\vspace{1.7mm}
\\
\cml\mathbf{KWY4}\cmr =\cml\mathbf{S_K}(1,6)\cmr =
\mathrm{[Vvvv+Dddd+Dcc1+Dwww1+...+Dwww6]}
\\
\hphantom{\cml}\mathbf{KWY4}\hphantom{\,}=\,\mathbf{K} \oplus \mathrm{WY4}
=\, \mathbf{K} \oplus
``\prod^\circ \!p\lozenge p\rightarrow
\lozenge(pqr)
(\sum^\circ(\square(p\rightarrow q+r) + \square(qr\rightarrow p))
+\square(p+q+r))"
\\
\hphantom{\cml\mathbf{KWY6}=\,\mathbf{K} \oplus \mathrm{WY6}}
=\, \mathbf{K} \oplus
``pqr\rightarrow
\lozenge(pqr)(\sum^\circ(\square(p\rightarrow q+r) + \square(qr\rightarrow p)) + \square(p+q+r))"
\\
\phantom{\qquad\qquad\qquad\qquad\qquad\qquad\quad}
+ \sum^\circ\square(p\leftrightarrow q)" 
\vspace{1.7mm}
\\
\cml\mathbf{KWY3}\cmr =\cml\mathbf{S_K}(1,5)\cmr =
\mathrm{[Vvvv+Dddd+Dcc1+Dwww1+...+Dwww5]}
\\
\hphantom{\cml}\mathbf{KWY3}\hphantom{\,}=\,\mathbf{K} \oplus \mathrm{WY3}
=\, \mathbf{K} \oplus
``\prod^\circ \!p\lozenge p\rightarrow
\lozenge(pqr)\sum^\circ(\square(p\rightarrow q+r) + \square(qr\rightarrow p))" 
\\
\hphantom{\cml\mathbf{KWY6}=\,\mathbf{K} \oplus \mathrm{WY6}}
=\, \mathbf{K} \oplus
``pqr\rightarrow
\lozenge(pqr)\sum^\circ(\square(p\rightarrow q+r) + \square(qr\rightarrow p))
+ \sum^\circ\square(p\leftrightarrow q)" 
\vspace{1.7mm}
\\
\cml\mathbf{KWY2}\cmr =\cml\mathbf{S_K}(0,6)\cmr =
\mathrm{[Vvvv+Dddd+Dwww1+...+Dwww6]}
\\
\hphantom{\cml}\mathbf{KWY2}\hphantom{\,}=\,\mathbf{K} \oplus \mathrm{WY2}
=\, \mathbf{K} \oplus
``pqr\rightarrow
\lozenge(pqr)\sum^\circ(\square(p\rightarrow q+r) + \square(qr\rightarrow p)) + \square(p+q+r)" 
\\
\hphantom{\cml\mathbf{KWY2}\hphantom{\,}=\,\mathbf{K} \oplus \mathrm{WY2}}
=\, \mathbf{K} \oplus
``pqr\rightarrow
\sum^\circ\!\lozenge p\,(\square(p\rightarrow q+r) + \square(qr\rightarrow p))+\square(p+q+r)" 
\vspace{1.7mm}
\\
\cml\mathbf{KWY1}\cmr =\cml\mathbf{S_K}(0,5)\cmr =
\mathrm{[Vvvv+Dddd+Dwww1+...+Dwww5]}
\\
\hphantom{\cml}\mathbf{KWY1}\hphantom{\,}=\,\mathbf{K} \oplus \mathrm{WY1}
=\, \mathbf{K} \oplus
``pqr\rightarrow
\lozenge(pqr)\sum^\circ(\square(p\rightarrow q+r) + \square(qr\rightarrow p)) + \square(pqr)" 
\\
\hphantom{\cml\mathbf{KWY1}\hphantom{\,}=\,\mathbf{K} \oplus \mathrm{WY1}}
=\, \mathbf{K} \oplus
``pqr\rightarrow
\sum^\circ\!\lozenge p\,(\square(p\rightarrow q+r) + \square(qr\rightarrow p)) + \square(pqr)" 
\\
\hphantom{\cml\mathbf{KWY1}\hphantom{\,}=\,\mathbf{K} \oplus \mathrm{WY1}}
=\, \mathbf{K} \oplus
``\lozenge(pqr)\rightarrow
\sum^\circ\lozenge p(\square(p\rightarrow q+r) + \square(qr\rightarrow p))" 
\vspace{1.7mm}
\\
\cml\mathbf{KWY0}\cmr =\cml\mathbf{S_K}(5,4)\cmr =
\mathrm{[Vvvv+Dddd+Dccc1+...+Dccc5+Dwww1+...+Dwww4]}
\\
\hphantom{\cml}\mathbf{KWY0}\hphantom{\,}=\,\mathbf{K} \oplus \mathrm{WY0}
=\, \mathbf{K} \oplus
``\sum^\circ(\square(p\rightarrow q+r) + \square(qr\rightarrow p))" 
\\
\hphantom{\cml\mathbf{KWY0}\hphantom{\,}=\,\mathbf{K} \oplus \mathrm{WY0}}
=\, \mathbf{K} \oplus
``\sum^\circ(\square(p+(q\leftrightarrow r)) + \square(qr\rightarrow p)) + \square(p+q+r)" 
\\
\hphantom{\cml\mathbf{KWY0}\hphantom{\,}=\,\mathbf{K} \oplus \mathrm{WY0}}
=\, \mathbf{K} \oplus
``\lozenge(pqr)\rightarrow
\sum^\circ(\square(p\rightarrow q+r) + \square(p+(q\leftrightarrow r)))" 
\vspace{1.7mm}
\\
\cml\mathbf{KWX9}\cmr =\cml\mathbf{S_K}(4,4)\cmr =
\mathrm{[Vvvv+Dddd+Dccc1+...+Dccc4+Dwww1+...+Dwww4]}
\\
\hphantom{\cml}\mathbf{KWX9}\hphantom{\,}=\,\mathbf{K} \oplus \mathrm{WX9}
=\, \mathbf{K} \oplus
``pqr\rightarrow
\sum^\circ(\square(p\rightarrow q+r) + \square(p+(q\leftrightarrow r)))" 
\\
\hphantom{\cml\mathbf{KWX9}\hphantom{\,}=\,\mathbf{K} \oplus \mathrm{WX9}}
=\, \mathbf{K} \oplus
``pqr\rightarrow
\sum^\circ(\square(p+(q\leftrightarrow r)) + \square(p\rightarrow(q\leftrightarrow r))) + \square(p+q+r)" 
\vspace{1.7mm}
\\
\pagebreak{}
\\
\cml\mathbf{KWX8}\cmr =\cml\mathbf{S_K}(4,3)\cmr =
\mathrm{[Vvvv+Dddd+Dccc1+...+Dccc4+Dwww1+...+Dwww3]}
\\
\hphantom{\cml}\mathbf{KWX8}\hphantom{\,}=\,\mathbf{K} \oplus \mathrm{WX8}
=\, \mathbf{K} \oplus
``\lozenge(pqr)\rightarrow
\sum^\circ\square(p\rightarrow q) + \square(p+q+r)" 
\\
\hphantom{\cml\mathbf{KWX8}\hphantom{\,}=\,\mathbf{K} \oplus \mathrm{WX8}}
=\, \mathbf{K} \oplus
``\lozenge(pqr)\rightarrow
\sum^\circ\square(p\rightarrow q+r) + \square(p+q+r)" 
\\
\hphantom{\cml\mathbf{KWX8}\hphantom{\,}=\,\mathbf{K} \oplus \mathrm{WX8}}
=\, \mathbf{K} \oplus
``\lozenge(pqr)\rightarrow
\sum^\circ\square(qr\rightarrow p) + \square(p+q+r)" 
\\
\hphantom{\cml\mathbf{KWX8}\hphantom{\,}=\,\mathbf{K} \oplus \mathrm{WX8}}
=\, \mathbf{K} \oplus
``\lozenge(pqr)\rightarrow
\sum^\circ(\square(qr\rightarrow p) + \square(q+r))" 
\\
\hphantom{\cml\mathbf{KWX8}\hphantom{\,}=\,\mathbf{K} \oplus \mathrm{WX8}}
=\, \mathbf{K} \oplus
``\lozenge(pqr)\rightarrow
\sum^\circ(\square(p\rightarrow q+r) + \square(q+r\rightarrow p)) + \square(p+q+r)" 
\\
\hphantom{\cml\mathbf{KWX8}\hphantom{\,}=\,\mathbf{K} \oplus \mathrm{WX8}}
=\, \mathbf{K} \oplus
``\sum^\circ(\square(p\rightarrow q+r) + \square(p\rightarrow(q\leftrightarrow r)))" 
\vspace{1.7mm}
\\
\cml\mathbf{KWX7}\cmr =\cml\mathbf{S_K}(3,4)\cmr =
\mathrm{[Vvvv+Dddd+Dccc1+...+Dccc3+Dwww1+...+Dwww4]}
\\
\hphantom{\cml}\mathbf{KWX7}\hphantom{\,}=\,\mathbf{K} \oplus \mathrm{WX7}
=\, \mathbf{K} \oplus
``pqr\rightarrow
\!\lozenge(pqr)\!\sum^\circ\!\square(p\!\rightarrow\!q+r)
 + \sum^\circ(\square(pq\!\leftrightarrow\!pr) + \square(p+(q\!\leftrightarrow\!r)))" 
\\
\hphantom{\cml\mathbf{KWX7}\hphantom{\,}=\,\mathbf{K} \oplus \mathrm{WX7}}
=\, \mathbf{K} \oplus
``pqr\rightarrow
(\lozenge(pqr)+\square(p+q+r))\sum^\circ\square(p\rightarrow q+r)
\\
\phantom{\qquad\qquad\qquad\qquad\qquad\qquad\quad}
+ \sum^\circ(\square(pq\leftrightarrow pr) + \square(p\rightarrow q))" 
\vspace{1.7mm}
\\
\cml\mathbf{KWX6}\cmr =\cml\mathbf{S_K}(3,3)\cmr =
\mathrm{[Vvvv+Dddd+Dccc1+...+Dccc3+Dwww1+...+Dwww3]}
\\
\hphantom{\cml}\mathbf{KWX6}\hphantom{\,}=\,\mathbf{K} \oplus \mathrm{WX6}
=\, \mathbf{K} \oplus
``\sum^\circ\square(p\rightarrow q) + \square(p+q+r)" 
\\
\hphantom{\cml\mathbf{KWX6}\hphantom{\,}=\,\mathbf{K} \oplus \mathrm{WX6}}
=\, \mathbf{K} \oplus
``\sum^\circ\square(p\rightarrow q+r) + \square(p+q+r)" 
\\
\hphantom{\cml\mathbf{KWX6}\hphantom{\,}=\,\mathbf{K} \oplus \mathrm{WX6}}
=\, \mathbf{K} \oplus
``\sum^\circ\square(qr\rightarrow p) + \square(p+q+r)" 
\vspace{1.7mm}
\\
\cml\mathbf{KWX5}\cmr =\cml\mathbf{S_K}(2,4)\cmr =
\mathrm{[Vvvv+Dddd+Dccc1+Dccc2+Dwww1+...+Dwww4]}
\\
\hphantom{\cml}\mathbf{KWX5}\hphantom{\,}=\,\mathbf{K} \oplus \mathrm{WX5}
=\, \mathbf{K} \oplus
``pqr\rightarrow
\lozenge(pqr)\sum^\circ\square(p\rightarrow q+r)
+ \sum^\circ(\square(pq\leftrightarrow pr) + \square p)"
\vspace{1.7mm}
\\
\cml\mathbf{KWX4}\cmr =\cml\mathbf{S_K}(2,3)\cmr =
\mathrm{[Vvvv+Dddd+Dccc1+Dccc2+Dwww1+...+Dwww3]}
\\
\hphantom{\cml}\mathbf{KWX4}\hphantom{\,}=\,\mathbf{K} \oplus \mathrm{WX4}
=\, \mathbf{K} \oplus
``pqr\rightarrow
\lozenge(pqr)\square(p+q+r)
+ \sum^\circ\square(p\rightarrow q)"
\\
\hphantom{\cml\mathbf{KWX4}\hphantom{\,}=\,\mathbf{K} \oplus \mathrm{WX4}}
=\, \mathbf{K} \oplus
``pqr\rightarrow
\lozenge(pqr)\sum^\circ\square(p\rightarrow q) + \sum^\circ\square(p\leftrightarrow q) + \square(p+q+r)" 
\\
\hphantom{\cml\mathbf{KWX4}\hphantom{\,}=\,\mathbf{K} \oplus \mathrm{WX4}}
=\, \mathbf{K} \oplus
``\prod^\circ\!p\lozenge p\rightarrow
\lozenge(pqr)\sum^\circ\square(p\rightarrow q) + \square(p+q+r)" 
\vspace{1.7mm}
\\
\cml\mathbf{KWX3}\cmr =\cml\mathbf{S_K}(1,4)\cmr =
\mathrm{[Vvvv+Dddd+Dccc1+Dwww1+...+Dwww4]}
\\
\hphantom{\cml}\mathbf{KWX3}\hphantom{\,}=\,\mathbf{K} \oplus \mathrm{WX3}
=\, \mathbf{K} \oplus
``pqr\rightarrow
\lozenge(pqr)\sum^\circ\square(qr\rightarrow p)
+ \sum^\circ\square(p+(q\leftrightarrow r))"
\\
\hphantom{\cml\mathbf{KWX3}\hphantom{\,}=\,\mathbf{K} \oplus \mathrm{WX3}}
=\, \mathbf{K} \oplus
``pqr\rightarrow
\lozenge(pqr)\sum^\circ\square(p\rightarrow q+r)
+ \sum^\circ\square(pq\leftrightarrow pr)"
\vspace{1.7mm}
\\
\cml\mathbf{KWX2}\cmr =\cml\mathbf{S_K}(1,3)\cmr =
\mathrm{[Vvvv+Dddd+Dccc1+Dwww1+...+Dwww3]}
\\
\hphantom{\cml}\mathbf{KWX2}\hphantom{\,}=\,\mathbf{K} \oplus \mathrm{WX2}
=\, \mathbf{K} \oplus
``pqr\rightarrow
\lozenge(pqr)(\sum^\circ\square(p\rightarrow q) + \square(p+q+r))
+ \sum^\circ\square(p\leftrightarrow q)"
\\
\hphantom{\cml\mathbf{KWX2}\hphantom{\,}=\,\mathbf{K} \oplus \mathrm{WX2}}
=\, \mathbf{K} \oplus
``\prod^\circ\!p\lozenge p\rightarrow
\lozenge(pqr)(\sum^\circ\square(p\rightarrow q) + \square(p+q+r))"
\vspace{1.7mm}
\\
\cml\mathbf{KWX1}\cmr =\cml\mathbf{S_K}(0,4)\cmr =
\mathrm{[Vvvv+Dddd+Dwww1+...+Dwww4]}
\\
\hphantom{\cml}\mathbf{KWX1}\hphantom{\,}=\,\mathbf{K} \oplus \mathrm{WX1}
=\, \mathbf{K} \oplus
``pqr\rightarrow
\lozenge(pqr)\sum^\circ(\square(qr\rightarrow p) + \square(p+(q\leftrightarrow r)))
+ \square(pqr)"
\\
\hphantom{\cml\mathbf{KWX1}\hphantom{\,}=\,\mathbf{K} \oplus \mathrm{WX1}}
=\, \mathbf{K} \oplus
``pqr\rightarrow
\lozenge(pqr)\sum^\circ(\square(p\rightarrow q+r) + \square(p\rightarrow(q\leftrightarrow r)))
+ \square(pqr)"
\\
\hphantom{\cml\mathbf{KWX1}\hphantom{\,}=\,\mathbf{K} \oplus \mathrm{WX1}}
=\, \mathbf{K} \oplus
``\sum^\circ(\lozenge p\square p\rightarrow p)( \square(p\rightarrow q+r) + \square(qr\rightarrow p) )"
\vspace{1.7mm}
\\
\cml\mathbf{KWX0}\cmr =\cml\mathbf{S_K}(0,3)\cmr =
\mathrm{[Vvvv+Dddd+Dwww1+...+Dwww3]}
\\
\hphantom{\cml}\mathbf{KWX0}\hphantom{\,}=\,\mathbf{K} \oplus \mathrm{WX0}
=\, \mathbf{K} \oplus
``pqr\rightarrow
\lozenge(pqr)\sum^\circ\square(q\rightarrow r) + \square(p+q+r)"
\\
\hphantom{\cml\mathbf{KWX0}\hphantom{\,}=\,\mathbf{K} \oplus \mathrm{WX0}}
=\, \mathbf{K} \oplus
``pqr\rightarrow
\sum^\circ\!\lozenge p\,\square(q\rightarrow r) + \square(p+q+r)"
\\
\hphantom{\cml\mathbf{KWX0}\hphantom{\,}=\,\mathbf{K} \oplus \mathrm{WX0}}
=\, \mathbf{K} \oplus
``pqr\rightarrow
\sum^\circ\!\lozenge p\,\square(qr\rightarrow p) + \square(p+q+r)"
\\
\hphantom{\cml\mathbf{KWX0}\hphantom{\,}=\,\mathbf{K} \oplus \mathrm{WX0}}
=\, \mathbf{K} \oplus
``pqr\rightarrow
\sum^\circ\!\lozenge p\,\square(p\rightarrow q+r) + \square(p+q+r)"
\vspace{1.7mm}
$

\paragraph{The pre-normal system R}
\fontsize{11pt}{12pt}\noindent
\setlength{\leftskip}{0cm}The system commonly denoted by $\mathbf{R}$
is the extension of the classical modal logic $\mathbf{E}$
with axiom $\mathrm{R}=\lozenge(p+q)\leftrightarrow\lozenge p+\lozenge q$.
(Note: This is in fact an $\mathbf{E}[2,1]$-system
that we called $\mathbf{R_{cw}}$ in \cite{Soncodi},
as it only sinks into another, ``genuine" system $\mathbf{R}$
that we defined as an $\mathbf{E}[0,1]$-system.
But here we continue to refer to this system $\mathbf{R}$ by its popular name.)

We call $\mathbf{R}$ \emph{pre-normal} because it is similar to $\mathbf{K}$
except that it lacks the normality axiom $\mathrm{N}=\square\mathit{1}$.
But axiom $\mathrm{R}$ still allows us to use
the previously-defined $\nu_j=\lozenge\mu_j, 0\leq j<n$
to determine the state of all other $\lozenge e_i$,
except for the modal factor $\nu_{\varnothing}=\lozenge\mathit{0}$,
which can be in any state,
so we need to add it to every $\mathbf{R}$-minmatrix.
Thus, in $\mathbf{R}[1,1]$ for example:\vspace{-2mm}

\begin{center}
$_1^1\cml \mathbf{R}\cmr = [\mathit{1}]=
\begin{array}{c|cccccccccc|}
         p & 1 & 1 & 1 & 1 & 1 & 0 & 0 & 0 & 0 & 0 \\
\hline
\lozenge p & 1 & 1 & 1 & 0 & 0 & 1 & 1 & 1 & 0 & 0 \\
\lozenge!p & 1 & 1 & 0 & 1 & 0 & 1 & 1 & 0 & 1 & 0 \\
\lozenge
\mathit{0} & 1 & 0 & 0 & 0 & 0 & 1 & 0 & 0 & 0 & 0
\end{array}$
\end{center}\vspace{-2mm}

\noindent\vspace{0mm}and we observe that, in addition to the normal minterms,
some pre-normal minterms show up, where the state of $\lozenge\mathit{0}$ is 1.
Then, by axiom R, these pre-normal minterms must have all modal factors in state 1, 
so they end up in the same prime orbit, which we label $\mathrm{Ww_0}$.

We note that $\mathbf{R}$ is not a \emph{perfect base},
since $_1^1\cml \mathbf{R}\cmr$ does not have the full $2^4=16$ minterms.
Nevertheless, the theory from this paper holds with only minor adjustments.
Then it turns out that prime orbit $\mathrm{Ww_0}$ is similar to $\mathrm{Vv_0}$,
in the sense that it covers only itself under substitutions
and therefore can be added to any $\mathbf{S_D}(x,y)$ and $\mathbf{S_K}(x,y)$ CMM
without collapsing it.
So $\mathbf{R}\mathrm{sys}\cml *,1\cmr$
is similar to $\mathbf{K}\mathrm{sys}\cml *,1\cmr$,
but adds a J-plane and an R-plane:

\begin{figure}[H]\vspace{-2mm}
\begin{centering}
\setlength{\leftskip}{0.15cm}\includegraphics[scale=0.63]{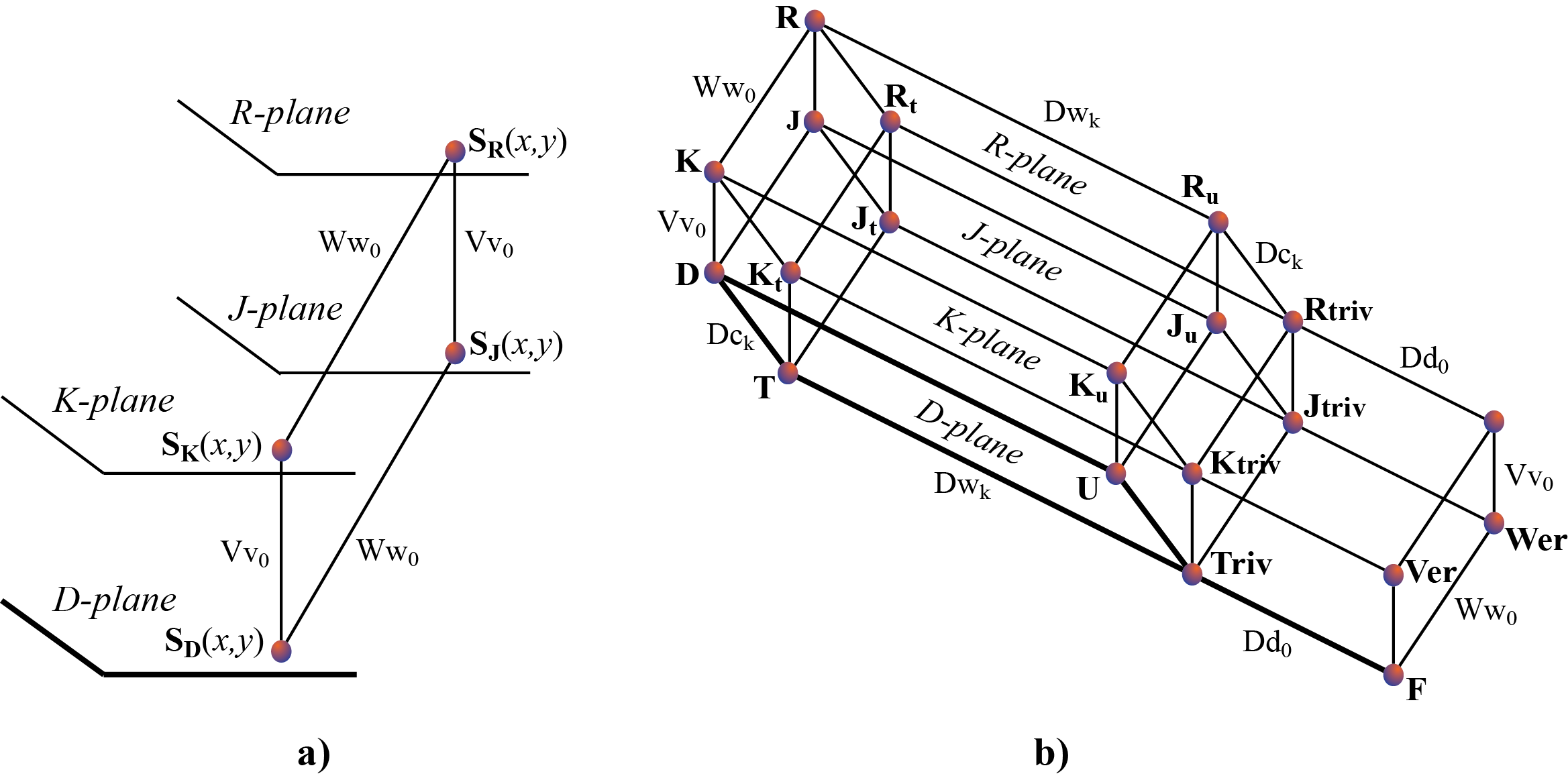}
\par\end{centering}
 	
\vspace{-1mm}\caption{\label{fig:R-systems}The lattice $\mathbf{R}\mathrm{sys}\cml *,1\cmr$}
\end{figure}

\pagebreak{}

\setlength{\leftskip}{0cm}

\vspace{5mm}

\fontsize{9pt}{11pt}\noindent
\emph{~~~This is a revised version of a paper submitted to Studia Logica in January 2016.}\vspace{1mm}


\AuthorAdressEmail{Adrian Soncodi}
{Lecturer, University of Texas at Dallas}
{acs151130@utdallas.edu\\
soncodi@verizon.net}


\vspace{4mm}

\begin{thebibliography}{1}

\bibitem{Blackburn}
{\sc Blackburn},~P., M.~{\sc de~Rijke}, and Y.~{\sc Venema},
{\em Modal Logic},
Cambridge University Press, 2001.

\bibitem{Fine}
{\sc Fine},~K.,
{\em Normal Forms in Modal Logic},
Notre Dame Journal of Formal Logic, Volume XVI, Number 2, April 1975.

\bibitem{Hughes}
{\sc Hughes},~G.W., M.J.~{\sc Cresswell},
{\em A New Introduction to Modal Logic},
Routledge 1996

\bibitem{Lewis}
{\sc Lewis},~D.,
{\em Intensional Logics without Iterative Axioms},
Journal of Philosophical Logic, Volume 3, 1974

\bibitem{Soncodi}
{\sc Soncodi},~A.,
{\em Automorphisms of the Lattice of Classical Modal Logics},
Studia Logica, Volume 104, 2015

\bibitem{Surendonk}
{\sc Surendonk},~T.J.
{\em Canonicity for Intensional Logics without Iterative Axioms},
Journal of Philosophical Logic, Volume 26, 1996


\end{thebibliography}

\begin{thebibliography}{1}
	
	\bibitem{bg1}
	{\sc Blackburn},~P., M.~{\sc de~Rijke}, and Y.~{\sc Venema},
	{\em Modal Logic},
	Cambridge University Press, 2001.
	
	\bibitem{bg2}
	{\sc Fine},~K.,
	{\em Normal Forms in Modal Logic},
	Notre Dame Journal of Formal Logic, Volume XVI, Number 2, April 1975.
	
	\bibitem{bg3}
	{\sc Hughes},~G.W., M.J.~{\sc Cresswell},
	{\em A New Introduction to Modal Logic},
	Routledge 1996
	
	\bibitem{bg4}
	{\sc Lewis},~D.,
	{\em Intensional Logics without Iterative Axioms},
	Journal of Philosophical Logic, Volume 3, 1974
	
	\bibitem{bg5}
	{\sc Soncodi},~A.,
	{\em Automorphisms of the Lattice of Classical Modal Logics},
	Studia Logica, Volume 104, 2015
	
	\bibitem{bg6}
	{\sc Surendonk},~T.J.
	{\em Canonicity for Intensional Logics without Iterative Axioms},
	Journal of Philosophical Logic, Volume 26, 1996
	
\end{thebibliography}
\end{document}
